\documentclass{jams-l}

\usepackage{amssymb}
\usepackage{amsmath}

\usepackage{graphicx}

\usepackage[cmtip,all]{xy}

\usepackage{stmaryrd}
\usepackage{comment}
\usepackage{tikz}
\usepackage{caption}
\usepackage{subcaption}
\usepackage{float}
\usepackage{wrapfig}
\usepackage{bm}
\usepackage{dashbox}
\usepackage{multirow}
\usepackage{xcolor}
\usepackage{hyperref}
\usepackage{cancel}

\usepackage{todonotes}

\usetikzlibrary{quotes,angles}

\newtheorem{theorem}{Theorem}[section]
\newtheorem{lemma}[theorem]{Lemma}
\newtheorem{proposition}[theorem]{Proposition}
\newtheorem{proposition-definition}[theorem]{Proposition-Definition}
\newtheorem{corollary}[theorem]{Corollary}

\theoremstyle{definition}
\newtheorem{definition}[theorem]{Definition}

\theoremstyle{remark}
\newtheorem{remark}[theorem]{Remark}

\numberwithin{equation}{section}

\DeclareMathOperator{\CuspEnd}{\mathcal{C}_a}
\DeclareMathOperator{\gammaAW}{\gamma_{\scriptscriptstyle{AW}}}
\DeclareMathOperator{\dom}{Dom}
\DeclareMathOperator{\Fp}{Fp}

\newcommand{\CuspidalEnd}[1]{\mathcal{C}_{#1}}

\begin{document}

\title[Pseudo-Laplacian on a cuspidal end: Alvarez--Wentworth conditions]{Pseudo-Laplacian on a cuspidal end with a flat unitary line bundle: Alvarez--Wentworth boundary conditions}


\author{Mathieu Dutour}
\address{University of Alberta, Department of Mathematical and Statistical Sciences}
\curraddr{}
\email{dutour@ualberta.ca}
\thanks{}



\date{}

\dedicatory{}

\begin{abstract}

    A cuspidal end is a type of metric singularity, described as a product $S^1 \times \left] a, +\infty \right[$ with the Poincaré metric. The underlying set can also be seen as $\mathbb{R} \times \left] a, +\infty \right[$ subject to the action of the translation $T : \left( x,y \right) \longrightarrow \left( x+1, y \right)$. On it, one may consider a holomorphic line bundle $L$, coming from a unitary character of the group generated by $T$. The complex modulus induces a flat metric on $L$, and a pseudo-Laplacian $\Delta_{L,0}$ acting on functions can be associated to the Chern connection. One needs to specify boundary conditions, and they are here chosen to be the Alvarez--Wentworth boundary conditions, which are a combination of Dirichlet and Neumann boundary conditions. The aim of this paper is to find the asymptotic behavior of the zeta-regularized determinant $\det \left( \Delta_{L,0} + \mu \right)$, as $\mu > 0$ goes to infinity for any $a$, and also as $a$ goes to infinity for $\mu = 0$. \vspace{10pt}
\end{abstract}

\maketitle

\tableofcontents


\section{Introduction}

    This paper is a follow-up to \cite{dutour:cuspDirichlet}, where a closely related situation was studied, using similar ideas and techniques. As such, the reader is referred to the introduction of \cite{dutour:cuspDirichlet} for a more detailed survey of an analogous layout. Nevertheless, parts of the introduction below will be the same as in \cite{dutour:cuspDirichlet}, with a focus on the differences between the two papers.

    \subsection{Description of the situation}
    
        In this paper, as in \cite{dutour:cuspDirichlet}, we study a type of metric singularity on Riemann surfaces called \textit{cusps}, with flat unitary holomorphic line bundles on them. These occur in a natural way when studying modular curves induced by a Fuchsian group of the first kind $\Gamma$, together with a flat unitary holomorphic vector bundle associated with a unitary representation of $\Gamma$. The computations presented here will eventually be used to derive a \textit{Deligne--Riemann--Roch isometry} which extends the work of Freixas i Montplet and von Pippich discussed in \cite{freixas-vonPippich:RROrbifold}. It is expected that this isometry will lead to a formula relating the first non-zero derivative at $1$ of the Selberg zeta function to some arithmetic intersection numbers. In particular, when the representation of $\Gamma$ is a non-trivial character, this relation will provide information on the Weil--Peterson norms of weight $1$ modular forms with Nebentypus.
        
        \subsubsection{Metric singularities}
        
            A \textit{cuspidal end} $\CuspEnd$ of height $a > 0$ is defined as the quotient of $\mathbb{R} \times \left] a, + \infty \right[$ by the action of the translation $\left( x, y \right) \mapsto \left( x + 1, y \right)$, which can be identified to $S^1 \times \left] a, +\infty \right[$, endowed with the Poincaré metric
            \begin{equation}
            \label{eq:intro:defPoincareMetric}
                \begin{array}{lll}
                    \mathrm{d}s_{\mathrm{hyp}}^2 & = & \frac{\mathrm{d}x^2 + \mathrm{d}y^2}{y^2}.
                \end{array}
            \end{equation}
            The coordinate $z = \exp \left( 2 i \pi \left( x + i y \right) \right)$ further identifies the cuspidal end to a punctured disk of radius $\varepsilon = \exp \left( - 2 \pi a \right)$, on which the Poincaré metric is
            \begin{equation}
            \label{eq:intro:defPoincareMetricPuncturedDisk}
                \begin{array}{lll}
                    \mathrm{d}s_{\mathrm{hyp}}^2 & = & \frac{\left \vert \mathrm{d}z \right \vert^2}{\left( \left \vert z \right \vert \log \left \vert z \right \vert \right)^2}.
                \end{array}
            \end{equation}
            The metric \eqref{eq:intro:defPoincareMetricPuncturedDisk} cannot be smoothly extended at $z=0$, and is thus referred to as a ``singular metric''. Let us also consider a unitary character $\chi : \mathbb{Z} \longrightarrow \mathbb{C}^{\ast}$, which provides an action of $\mathbb{Z}$ onto the trivial line bundle $\left( \mathbb{R} \times \left] a, +\infty \right[ \right) \times \mathbb{C}$ by
            \begin{equation}
            \label{eq:intro:actionUnitaryCharacter}
                \begin{array}{lll}
                    k \cdot \left( \left( x, y \right), \lambda \right) & = & \left( \left( x + k, y \right), \chi \left( k \right) \lambda \right).
                \end{array}
            \end{equation}
            The quotient of $\left( \mathbb{R} \times \left] a, +\infty \right[ \right) \times \mathbb{C}$ by $\mathbb{Z}$ under the action \eqref{eq:intro:actionUnitaryCharacter} yields a flat unitary holomorphic line bundle $L$ on the cuspidal end, which is entirely determined by the complex number of modulus $1$
            \begin{equation}
            \label{eq:intro:chiOf1}
                \begin{array}{lll}
                    \chi \left( 1 \right) & = & e^{2 i \pi \alpha},
                \end{array}
            \end{equation}
            where $\alpha$ is a real number well-defined modulo $1$. The line bundle $L$ is endowed with a \textit{canonical metric}, induced by the complex modulus on $\mathbb{C}$, which is fully characterized by the choice of a representative modulo $1$ of $\alpha$, which we assume here belongs to $\left[ 0,1 \right[$. Although $L$ can be extended over the cusp, by a process called Deligne's canonical extension, the metric cannot, in general, be smoothly continued.

        \subsubsection{Pseudo-Laplacian}
        
            In this paper, we consider the \textit{pseudo-Laplacian} with \textit{Alvarez--Wentworth boundary conditions}. The change in boundary conditions marks the difference between the present paper and \cite{dutour:cuspDirichlet}. These can be though of as a blend between the Dirichlet and the Neumann boundary conditions, having been introduced by Alvarez in \cite[Sec. 4]{alvarez:theoryStrings}, and extended to (metrically) non-trivial bundles by Wentworth in \cite{wentworth:gluingFormulasAlvarez}. Let us briefly describe them. The line bundle $L$ on the cuspidal end can be trivialized by the section induced by
            \begin{equation}
                \begin{array}{lllll}
                    s_{\alpha} & : & \mathbb{R} \times \left] a, + \infty \right[ & \longrightarrow & \mathbb{C} \\[0.5em]
                    && \left( x, y \right) & \longmapsto & \exp \left( 2 i \pi \alpha \left( x + i y \right) \right)
                \end{array}
            \end{equation}
            Any smooth section of $L$ on $\CuspEnd$ can then be written as a product $f s_{\alpha}$, where $f$ is a smooth complex function on $\CuspEnd$. The Alvarez--Wentworth boundary conditions are then defined by
            \begin{equation}
                \begin{array}{lllllll}
                    \Im f \left( x, a \right) & = & 0 & \text{ and } & \frac{\partial \Re f}{\partial y} \left( x, a \right) & = & 0
                \end{array}
            \end{equation}
            which is to say the Dirichlet condition on the imaginary part and the Neumann condition on the real part. As noted in \cite{dutour:cuspDirichlet}, the notion of pseudo-Laplacians has been studied by Colin de Verdière in \cite{colinDeVerdiere:pseudoLaplacian1, colinDeVerdiere:pseudoLaplacian2}. The value of $\alpha \in \left[ 0, 1 \right[$ splits the discussion into two parts.
            
            \medskip
            
            \hspace{10pt} $\bullet$ If $L$ is (metrically) non-trivial, \textit{i.e.} if we have $\alpha > 0$, the Chern Laplacian acting on smooth sections satisfying the Alvarez--Wentworth boundary conditions is a symmetric operator, and its Friedrichs extension is a self-adjoint operator, called the \textit{Chern Laplacian with Alvarez--Wentworth boundary conditions}. For the purpose of this paper, in order to have a uniform terminology with the case of the trivial line bundle presented below, this Laplacian is renamed the \textit{pseudo-Laplacian with Alvarez--Wentworth boundary conditions}, and is denoted by $\Delta_{L,0}$.
            
            \medskip
            
            \hspace{10pt} $\bullet$ If $L$ is (metrically) trivial, \textit{i.e.} if we have $\alpha = 0$, the Chern Laplacian has an essential spectrum, which we must remove before computing the determinant. This is done by considering the orthogonal decomposition
            \begin{equation}
            \label{eq:intro:orthDecompPseudoLaplacian}
                \begin{array}{lll}
                    \scriptstyle L^2 \left( \CuspEnd, \frac{\mathrm{d}x^2 + \mathrm{d}y^2}{y^2} \right) & \scriptstyle = & \scriptstyle L^2 \left( \CuspEnd, \frac{\mathrm{d}x^2 + \mathrm{d}y^2}{y^2} \right)_0 \; \; \oplus \; \; L^2 \left( \left] a, + \infty \right[, \frac{\mathrm{d}y^2}{y^2} \right),
                \end{array}
            \end{equation}
            where the index $0$ means ``with vanishing constant Fourier coefficient''. The \textit{pseudo-Laplacian with Alvarez--Wentworth boundary conditions}, denoted by $\Delta_{L,0}$, is defined as the Chern Laplacian with Alvarez--Wentworth boundary conditions induced by restriction to the first space on the right-hand side of \eqref{eq:intro:orthDecompPseudoLaplacian}. Its determinant can be seen as a \textit{relative determinant}, following the notion introduced by M\"{u}ller in \cite{muller:relativeDeterminants}, of the Chern Laplacian with Alvarez--Wentworth boundary conditions $\Delta_L$ and of the Laplacian $-y^2 \mathrm{d}^2/\mathrm{d}y^2$ on $\left] a, + \infty \right[$ with Alvarez--Wentworth boundary conditions.

    \subsection{Statement of the main result}
    
        Similarly to \cite{dutour:cuspDirichlet}, we prove, this paper, two types of results related to the \textit{zeta-regularized determinant} of the pseudo-Laplacian with Alvarez--Wentworth boundary conditions. Both use explicitely the results obtained in \cite{dutour:cuspDirichlet}.
        
        \medskip
        
        \hspace{10pt} $\bullet$ Our first result, in theorems \ref{thm:asymptMuDeterminantAlphaNeq0} and \ref{thm:asymptMuDeterminantAlpha0}, is a $\mu$-aymptotic expansion
            \begin{equation}
            \label{AsymptoticStudyMu}
                \begin{array}{lll}
                    \log \det \left( \Delta_{L,0} + \mu \right) & = & \fbox{$\mu$-divergent part} + \fbox{$\mu$-constant term} + o \left( 1 \right)
                \end{array}
            \end{equation}
            
            \noindent for the logarithm of the zeta-regularized determinant of the pseudo-Laplacian (with Alvarez--Wentworth boundary conditions), as $\mu$ goes to infinity through strictly positive real values. This type of evaluation can be used to compute the constant in Mayer--Vietoris type formulas with parameter, in a way similar to what Burghelea, Friedlander, and Kappeler do in \cite[Sec. 3.19 \& 4.8]{bfk:mayerVietoris} for smooth metrics on compact manifolds, to what Wentworth does in \cite[Sec. 3.4]{wentworth:gluingFormulasAlvarez} for the boundary conditions we study here, and to what Carron does in \cite{carron:determinantRelatifFonctionXi} for potentially singular metrics.
        
        \medskip
        
        \hspace{10pt} $\bullet$ Our second result, in theorems \ref{thm:asymptADeterminantAlphaNeq0} and \ref{thm:asymptADeterminantAlpha0}, is an $a$-asymptotic expansion
            \begin{equation}
            \label{AsymptoticStudyA}
                \begin{array}{lll}
                    \log \det \Delta_{L,0} & = & \fbox{$a$-divergent part} + \fbox{$a$-constant term} + o \left( 1 \right)
                \end{array}
            \end{equation}
            
            \noindent for the logarithm of the zeta-regularized determinant of the pseudo-Laplacian (with Alvarez--Wentworth boundary conditions), as the height $a$ of the cuspidal end goes to infinity, \textit{i.e.} as the cusp shrinks, without parameter $\mu$.

    \subsection{Presentation of the paper}
    
        As was already indicated, the general layout of this paper is similar to that of \cite{dutour:cuspDirichlet}, and the reader is referred to the introduction of this earlier paper to get a more technical presentation of an analogous situation. Let us focus instead on the differences with \cite{dutour:cuspDirichlet}.
        
        \subsubsection{Step $1$: preliminary work on the pseudo-Laplacian}
        
            To begin this step, we briefly recall what cuspidal ends are, how to construct flat unitary holomorphic line bundles on them, and, in subsection \ref{SubSec:SobolevSpaces} the definitions of some of the classical Sobolev spaces we need. Then, in subsection \ref{SubSec:AWBoundaryConditions}, we see the precise definition of the Alvarez--Wentworth boundary conditions, and introduce other Sobolev type spaces, which are necessary to understand the Friedrichs extension process, used in subsection \ref{SubSec:pseudoLaplacian} to define the pseudo-Laplacian studied here. The main purpose of this first step is achieved in subsection \ref{SubSec:WeylLawSpectralZeta}, by obtaining a Weyl type law
            \begin{equation}
            \label{eq:intro:WeylTypeLaw}
                \begin{array}{lll}
                    N \left( \Delta_{L,0}, \lambda \right) & \leqslant & C \lambda
                \end{array},
            \end{equation}
            on the spectral counting function. To get this result, we have to assume, as in \cite{dutour:cuspDirichlet} that $\alpha$ is rational, which allows us to find an inequality which effectively removes the line bundle $L$. The resulting auxiliary Laplacian with Alvarez--Wentworth boundary conditions is then related to its Neumann counterpart, and we then refer to \cite{dutour:cuspDirichlet} for the remainder of the proof of the Weyl type law. As a consequence, the spectral zeta function of $\Delta_{L,0}$ exists and is holomorphic on the half-plane $\Re s > 1$.

        \subsubsection{Step $2$: localizing the eigenvalues}
        
            In subsections \ref{SubSec:spectralProblem} and \ref{SubSec:localizationEigenvalues}, we study the eigenvalues of the pseudo-Laplacian with Alvarez--Wentworth boundary conditions, by solving the spectral problem
            \begin{equation}
            \label{eq:intro:spectralProblem1}
                \left \{ \begin{array}{llll}
                    - y^2 \left( \frac{\partial^2}{\partial x^2} + \frac{\partial^2}{\partial y^2} \right) \psi \left( x, y \right) & = & \multicolumn{2}{l}{\lambda \psi \left( x, y \right)} \\[0.5em]
                
                    \psi \left( x+1, y \right) & = & \multicolumn{2}{l}{e^{2 i \pi \alpha} \psi \left( x, y \right)} \\[0.5em]
                
                    \int_{S^1 \times \left] a, +\infty \right[} \; \left \vert \psi \right \vert^2 & < & + \infty & \text{(Integrability condition)} \\[0.7em]
                
                    \Im \left( \psi s_{- \alpha} \right) \left( x, a \right) & = & 0 & \multirow{2}{*}{\text{(AW boundary condition)}} \\[0.5em]
                
                    \frac{\partial}{\partial y} \Re \left( \psi s_{- \alpha} \right) \left( x, a \right) & = & 0 \\[0.5em]
                
                    \int_{S^1} \; \psi \left( x, y \right) \; \mathrm{d}x & = & 0 &  \text{for almost all } y \geqslant a \text{ if } \alpha = 0
                \end{array} \right.
            \end{equation}
            in the same way as the analogous spectral problem in \cite{dutour:cuspDirichlet}. We find that the eigenvalues in \eqref{eq:intro:spectralProblem1} are characterized by the equations
            \begin{equation}
            \label{eq:intro:AWBoundaryCondition4}
                \begin{array}{lll}
                    \left( 1 + 4 \pi \alpha a \right) K_{s-1/2} \left( 2 \pi \left \vert k + \alpha \right \vert a \right) K_{s-1/2} \left( 2 \pi \left \vert k - \alpha \right \vert a \right) \\[0.5em]
                \quad + 2 \pi \left \vert k - \alpha \right \vert a K_{s-1/2} \left( 2 \pi \left \vert k + \alpha \right \vert a \right) K_{s-1/2} ' \left( 2 \pi \left \vert k - \alpha \right \vert a \right) \\[0.5em]
                \quad \quad \quad + 2 \pi \left \vert k + \alpha \right \vert a K_{s-1/2} \left( 2 \pi \left \vert k - \alpha \right \vert a \right) K_{s-1/2} ' \left( 2 \pi \left \vert k + \alpha \right \vert a \right) & = & 0
                \end{array}
            \end{equation}
            for $k \in \mathbb{Z}$, except $k = 0$ if we have $\alpha = 0$. Assuming we have $a > 1/\left( 4 \pi \left( 1 - \alpha \right) \right)$, the solutions of \eqref{eq:intro:AWBoundaryCondition4} are precisely located in corollary \ref{cor:locZerosKNeq0} for $k \neq 0$, in the same way as in \cite{dutour:cuspDirichlet} using an adaptation of Saharian's argument from \cite[App. A]{saharian}. The case $k = 0$, which must be considered if we have $\alpha \neq 0$, is different here, owing to the fact that the pseudo-Laplacian with ALvarez--Wentworth boundary conditions has a kernel. Nevertheless, proposition \ref{prop:zerosFunctionF0} localizes the associated contribution to the eigenvalues.

        \subsubsection{Step $3$: asymptotic studies}
        
            Using the argument principle on suitable contours with the left-hand side of \eqref{eq:intro:AWBoundaryCondition4}, we recover in subsection \ref{SubSec:contoursOfIntegration} the spectral zeta function of the pseudo-Laplacian with Alvarez--Wentworth boundary conditions. The contours can be seen in figures \ref{fig:integrationContourPM} and \ref{fig:integrationContourK0}. Note that the contour we need to compute the contribution from the case $k=0$ is more complicated because of the presence of a kernel for the pseudo-Laplacian. Thus, we obtain two ``partial'' spectral zeta functions $\zeta_{L,\mu}^{\pm} \left( s \right)$ and $\zeta_{L,\mu}^{0}$. The factorization \eqref{eq:functionf} then allows us to extract a Dirichlet contribution in subsection \ref{subsec:factorizationDirichletContribution}, thereby giving decompositions
            \begin{equation}
                \begin{array}{llll}
                    \zeta_{L,\mu}^{\pm} \left( s \right) & = & \zeta_{L,\mu}^{\pm, \, D} \left( s \right) + \zeta_{L,\mu}^{\pm, \, AW} \left( s \right) \\[1em]
                    \zeta_{L,\mu}^0 \left( s \right) & = & \zeta_{L,\mu}^{0, \, D} \left( s \right) + \zeta_{L,\mu}^{0, \, AW} \left( s \right) & \text{if we have } \alpha \neq 0
                \end{array}
            \end{equation}
            The study of the two Dirichlet functions $\zeta_{L,\mu}^{\pm, \, D}$ and $\zeta_{L,\mu}^{0, \, D}$ having been done in \cite{dutour:cuspDirichlet}, we are left with contributions coming purely from the Alvarez--Wentworth boundary conditions. We study these functions by successive splittings, a process similar to what is used in \cite{dutour:cuspDirichlet}.

        \subsubsection{Acknowledgements}
        
            The work in this paper originated from my doctoral thesis, and was realized as a postdoctoral fellow at the University of Alberta. Funding coming from the endowment of M. V. Subbarao Professorship in number theory and NSERC Discovery Grant RGPIN-2019-06112 should be acknowledged for that period. I would like to thank Gerard Freixas i Montplet, with whom I had many discussions on topics related to the computations presented in this paper, and whose advice has been invaluable. I would also like to thank Richard Wentworth for some exchanges we had during my PhD on the boundary conditions considered here, which helped me better understand them. Finally, I want to thank Manish Patnaik for our mathematical exchanges.

\section{Description of the spectral problem}
\label{Sec:DescriptionSpectralProblem}

    The situation considered here being the same as in \cite{dutour:cuspDirichlet}, the reader is referred to sections $2.1$ and $2.2$ of this paper, whose notations are adopted here, unless explicitely redefined, for details. Consider the cuspidal end $\CuspEnd = S^1 \times \left] a,+\infty \right[$ of height $a > 0$, endowed with the Poincaré metric 
    \begin{equation}
        \begin{array}{ccc}
            \mathrm{d}s^2_{\text{hyp}} & = & \frac{\mathrm{d}x^2  + \mathrm{d}y^2}{y^2}    
        \end{array} ,
    \end{equation}
    \noindent and a flat holomorphic unitary line bundle $L$ coming from a unitary character
    \begin{equation}
        \begin{array}{ccccccc}
            \chi & : & \mathbb{Z} & \longrightarrow & \mathbb{U} & \subset & \mathbb{C}^{\ast}
        \end{array} .
    \end{equation}
    \noindent The line bundle $L$ is then characterized by $\chi \left( 1 \right) = e^{2 i \pi \alpha}$, with $\alpha \in \left[ 0,1 \right[$, and a smooth section $s$ of $L$ over $\CuspEnd$ is identified to a smooth function
    \begin{equation}
    \label{identification:section:function}
        \begin{array}{ccccc}
            s & : & \mathbb{R} \times \left] a, + \infty \right[ & \longrightarrow & \mathbb{C}
        \end{array}
    \end{equation}
    \noindent such that we have $s \left( x+1, y \right) = e^{2 i \pi \alpha} s \left( x, y \right)$. Denote by $\nabla^L$ the Chern connection associated to the canonical metric on $L$, and by $\Delta_L$ the Laplacian associated to the Chern connection, acting on smooth sections of $L$. Furthermore, denote by $\nabla_{LC}$ the Levi--Civita connection. These two connections then induce connections, still denoted by $\nabla^L$, on each bundle of differential forms with values in $L$. Let us recall some definitions related to Sobolev spaces.

    \subsection{Sobolev spaces}
    \label{SubSec:SobolevSpaces}
    
        There are two Sobolev spaces we need to define, which are based on the space $L^2 \left( \CuspEnd, L \right)$ of $L^2$ sections of $L$.
        
        \begin{definition}
        
            The space of $L_1^2$ sections of $L$ is defined by
            \begin{equation}
                \begin{array}{ccc}
                    L_1^2 \left( \CuspEnd, L \right) & = & \left \{ s \in L^2 \left( \CuspEnd, L \right), \; \nabla^L s \in L^2 \left( \CuspEnd, \Omega \otimes L \right) \right \}.
                \end{array}
            \end{equation}
            \noindent In other words, it is the space of $L^2$ sections whose covariant derivative, in the distributional sense, is an $L^2$ section of $\Omega \otimes L$, where $\Omega$ is the cotangent bundle.
        
        \end{definition}
        
        \begin{proposition}
                
            The space of $L_1^2$ sections of $L$ is a Hilbert space for the norm
            \begin{equation}
                \begin{array}{ccc}
                    \left \Vert s \right \Vert_{L_1^2}^2 & = & \left \Vert s \right \Vert_{L^2}^2 + \left \Vert \nabla^L s \right \Vert_{L^2}^2
                \end{array} .
            \end{equation}
            \noindent Furthermore, this norm is called the $L_1^2$-norm.
            
        \end{proposition}
        
        \begin{proof}
        
            This is a classical result.
        
        \end{proof}
        
        \begin{definition}
        
            The space of $L_2^2$ sections of $L$ is defined by
            \begin{equation}
                \begin{array}{ccc}
                    L_2^2 \left( \CuspEnd, L \right) & = & \left \{ s \in L^2 \left( \CuspEnd, L \right), \; \Delta_L s \in L^2 \left( \CuspEnd, L \right) \right \}.
                \end{array}
            \end{equation}
            \noindent In other words, it is the space of $L^2$ sections whose Chern Laplacian, in the distributional sense, is again an $L^2$ section.
        
        \end{definition}
        
        \begin{proposition}
                
            The space of $L_2^2$ sections of $L$ is a Hilbert space for the norm
            \begin{equation}
                \begin{array}{ccc}
                    \left \Vert s \right \Vert_{L_2^2}^2 & = & \left \Vert s \right \Vert_{L^2}^2 + \left \Vert \Delta_L s \right \Vert_{L^2}^2
                \end{array} .
            \end{equation}
            \noindent Furthermore, this norm is called the $L_2^2$-norm.
                
        \end{proposition}
        
        \begin{proof}
        
            This is a classical result.
        
        \end{proof}
        
        \begin{remark}
        
            Some authors define these Sobolev spaces using alternatively iterations of the Chern connection or fractional powers of the Chern Laplacians. In our case, these yield the same spaces. For more details on the various definitions of Sobolev spaces and their comparisons, the reader is referred to \cite{eichhorn:openManifolds}.
        
        \end{remark}

    \subsection{The Alvarez--Wentworth boundary conditions}
    \label{SubSec:AWBoundaryConditions}
    
        The cuspidal end $\CuspEnd$ having a boundary, we need to consider a set of boundary conditions in order to study the Laplacian $\Delta_L$. In this paper, we choose the \textit{Alvarez--Wentworth} boundary conditions, first introduced by Alvarez in \cite[Sec. 4]{alvarez:theoryStrings} for functions on a surface, and extended by Wentworth in \cite{wentworth:gluingFormulasAlvarez} to sections of a line bundle on a compact Riemann surface. The non-compactness of $\CuspEnd$ does not change anything for this part.

        \subsubsection{Dirichlet boundary conditions}
        \label{SubSub:DirichletBoundaryConditions}
        
            There are two elementary sets of boundary conditions which play a role in defining the Alvarez--Wentworth boundary conditions. Let us begin by recalling the first one, which is of Dirichlet type.
        
            \begin{definition}
        
                The \textit{boundary trace operator} $\gamma_1$ is defined as
                \begin{equation}
                \label{boundaryTraceOperator}
                    \begin{array}{ccccc}
                        \gamma_1 & : & \mathcal{C}^{\infty} \left( \CuspEnd, L \right) & \longrightarrow & \mathcal{C}^{\infty} \left( S^1, L_{\vert y = a} \right) \\[0.5em]
                    
                        && s & \longmapsto & s \left( \cdot, \; a \right)
                    \end{array} .
                \end{equation}
            \end{definition}
        
            \begin{proposition}
                
                The map $\gamma_1$ can be extended into a continuous operator
                \begin{equation}
                    \begin{array}{ccccc}
                        \gamma_1 & : & L_1^2 \left( \CuspEnd, L \right) & \longrightarrow & L^2 \left( S^1, L_{\vert y = a} \right)
                    \end{array} .
                \end{equation}
            \end{proposition}
        
            \begin{proof}
        
                It is a consequence of the integration by parts formula.
        
            \end{proof}
            
            \begin{definition}
            
                A section $s \in L_1^2 \left( \CuspEnd, L \right)$ is said to satisfy the \textit{Dirichlet boundary conditions} if it belongs to
                \begin{equation}
                    \begin{array}{ccccc}
                        L_{1, \text{Dir}}^2 \left( \CuspEnd, L \right) & = & \ker \gamma_1 & = & \left \{ s \in L_1^2 \left( \CuspEnd, L \right), \; \; \gamma_1 s = 0 \right \}
                     \end{array} .
                \end{equation}
            \end{definition}
            
            \begin{remark}
            
                The Sobolev space $L_{1, \text{Dir}}^2$ is the closure, for the $L_1^2$-norm, of the space of smooth, compactly supported sections of $L$ with support in the interior of $\CuspEnd$.
            
            \end{remark}

        \subsubsection{Neumann boundary conditions}
        \label{SubSub:NeumannBoundaryConditions}
        
            The second set of boundary conditions which intervenes in the definition of the Alvarez--Wentworth boundary conditions is of Neumann type. Let us present those.
            
            \begin{definition}
            
                The \textit{normal trace operator} $\gamma_2$ is defined as
                \begin{equation}
                \label{normalTraceOperator}
                    \begin{array}{ccccc}
                        \gamma_1 & : & \mathcal{C}^{\infty} \left( \CuspEnd, L \right) & \longrightarrow & \mathcal{C}^{\infty} \left( S^1, L_{\vert y = a} \right) \\[0.5em]
                    
                        && s & \longmapsto & \frac{\partial s}{\partial y} \left( \cdot, a \right)
                    \end{array} .
                \end{equation}
            \end{definition}
            
            \begin{proposition}
                    
                The map $\gamma_2$ can be extended into a continuous operator
                \begin{equation}
                    \begin{array}{ccccc}
                        \gamma_2 & : & L_2^2 \left( \CuspEnd, L \right) & \longrightarrow & L^2 \left( S^1, L_{\vert y = a} \right)
                    \end{array} .
                \end{equation}
            \end{proposition}
            
            \begin{proof}
        
                It is a consequence of the integration by parts formula.
        
            \end{proof}
            
            \begin{definition}
            
                A section $s \in L_2^2 \left( \CuspEnd, L \right)$ is said to satisfy the \textit{Neumann boundary conditions} if it belongs to
                \begin{equation}
                    \begin{array}{ccccc}
                        L_{2, \text{Neu}}^2 \left( \CuspEnd, L \right) & = & \ker \gamma_2 & = & \left \{ s \in L_2^2 \left( \CuspEnd, L \right), \; \; \gamma_2 s = 0 \right \}
                     \end{array} .
                \end{equation}
            \end{definition}
            
            \begin{remark}
            
                The Sobolev space $L_{2, \text{Neu}}^2$ is the closure, for the $L_2^2$-norm, of the space of smooth sections of $L$ over $\CuspEnd$ whose normal derivative vanishes at the boundary.
            
            \end{remark}

        \subsubsection{Alvarez--Wentworth boundary conditions}
        \label{SubSub:AWBoundaryConditions}
        
            These boundary conditions are of mixed type. In other words, we will need to split a section of $L$ into two parts, impose the Dirichlet conditions on the first part, and the Neumann ones on the second. Let us see that more precisely. For clarity, we will repeat some of the content of \cite[Sec. 2.3]{dutour:cuspDirichlet}.
            
            \begin{proposition}
                
                The function
                \begin{equation}
                \label{trivializationL}
                    \begin{array}{ccccl}
                        s_{\alpha} & : & \mathbb{R} \times \left] a, + \infty \right[ & \longrightarrow & \mathbb{C} \\[0.5em]
                        
                        && \left( x,y \right) & \longmapsto & \exp \left( 2 i \pi \alpha \left( x + i y \right) \right)
                    \end{array}
                \end{equation}
                \noindent is a smooth, trivializing section of $L$.
            \end{proposition}
            
            \begin{proof}
            
                The only point to prove is that $s_{\alpha}$ trivializes $L$ over $\CuspEnd$. Let $s$ be a section of $L$ over the cuspidal end. Following \eqref{identification:section:function}, we identify $s$ to a function
                \begin{equation}
                    \begin{array}{lllll}
                        s & : & \mathbb{R} \times \left] a, + \infty \right[ & \longrightarrow & \mathbb{C}
                    \end{array}
                \end{equation}
                such that we have $s \left( x+1, y \right) = e^{2 i \pi \alpha} s \left( x, y \right)$ for every $\left( x, y \right) \in \mathbb{R} \times \left] a, + \infty \right[$. The function $s_{\alpha}$ having no zero, we can form the quotient
                \begin{equation}
                    \begin{array}{lllll}
                        \frac{s}{s_{\alpha}} & : & \mathbb{R} \times \left] a, + \infty \right[ & \longrightarrow & \mathbb{C}
                    \end{array}
                \end{equation}
                which is $1$-periodic in the first variable, thus inducing a (smooth) function on $\CuspEnd$, and proving the proposition.
            
            \end{proof}
            
            \begin{definition}
            
                Let $s$ be a smooth section of $L$, written as $s = f s_{\alpha}$, where $f$ is a smooth function on $\CuspEnd$. The \textit{real and imaginary parts of $s$} are defined by
                \begin{equation}
                    \begin{array}{ccccccc}
                        \Re s & = & \Re f & \text{ and } & \Im s & = & \Im f
                    \end{array},
                \end{equation}
                \noindent under the identification \eqref{identification:section:function}. These are not sections of $L$, but smooth functions on $\CuspEnd$. However, after multiplying them by $s_{\alpha}$, they do induce sections of $L$. Note that we have $f = s s_{-\alpha}$.
                
            \end{definition}
            
            \begin{definition}
                The \textit{Alvarez--Wentworth trace operator} $\gammaAW$ is defined as
                \begin{equation}
                \label{AWTrace}
                    \begin{array}{lllll}
                        \gamma_{\scriptscriptstyle{AW}} & : & \mathcal{C}^{\infty} \left( \CuspEnd, L \right) & \longrightarrow & \mathcal{C}^{\infty} \left( S^1, L_{\vert y = a} \right) \; \oplus \; \mathcal{C}^{\infty} \left( S^1, L_{\vert y = a} \right) \\[0.5em]
                        
                        && s & \longmapsto & \left( \left( \Im s \right) \left( x, a \right) s_{\alpha} \left( x, a \right), \; \frac{\partial \Re s}{\partial y} \left( x, a \right) s_{\alpha} \left( x, a \right) \right)
                    \end{array},
                \end{equation}
                under the identification \eqref{identification:section:function}.
            \end{definition}
            
            \begin{remark}
                The trace operator $\gammaAW$ is a blend between the boundary trace operator \eqref{boundaryTraceOperator} and the normal trace operator \eqref{normalTraceOperator}. It can also be defined in terms of the Dolbeault operator $\overline{\partial}_L$, as explained in \cite[Def. 2.1]{wentworth:gluingFormulasAlvarez}.
            \end{remark}
            
            \begin{proposition}
                The map $\gamma_{\scriptscriptstyle{AW}}$ can be extended into a continuous operator
                \begin{equation}
                    \begin{array}{lllll}
                        \gammaAW & : & L_2^2 \left( \CuspEnd, L \right) & \longrightarrow & L^2 \left( S^1, L_{\vert y = a} \right) \; \oplus \; L^2 \left( S^1, L_{\vert y = a} \right)
                    \end{array}.
                \end{equation}
            \end{proposition}
            
            \begin{proof}
                It is a consequence of the integration by parts formula.
            \end{proof}
            
            \begin{remark}
                Note that the operator $\gammaAW$ is only real linear, even though the vectors spaces involved are naturally complex.
            \end{remark}
            
            \begin{definition}
                A section $s \in L_2^2 \left( \CuspEnd, L \right)$ is said to satisfy the \textit{Alvarez--Wentworth boundary conditions} if it belongs to
                \begin{equation}
                \label{eq:SobolevAWL22}
                    \begin{array}{lllll}
                        L_{2, \text{AW}}^2 \left( \CuspEnd, L \right) & = & \ker \gammaAW & = & \left \{ s \in L_2^2 \left( \CuspEnd, L \right), \; \; \gammaAW s = 0 \right \}
                     \end{array} .
                \end{equation}
            \end{definition}

    \subsection{Pseudo-Laplacian}
    \label{SubSec:pseudoLaplacian}
    
        As in \cite[Sec. 2.3]{dutour:cuspDirichlet}, we need to introduce a pseudo-Laplacian, which will be the operator we study in this paper.
        
        \begin{proposition}
        \label{prop:friedrichsLaplacianAW}
            The Chern Laplacian $\Delta_L$, acting on smooth sections $s$ of $L$ satisfying $\gammaAW s = 0$, is a symmetric positive operator. Its Friedrichs extension is a positive $L^2$ self-adjoint operator
            \begin{equation}
                \begin{array}{lllll}
                    \Delta_L & : & L_{2, \text{AW}}^2 \left( \CuspEnd, L \right) & \longrightarrow & L^2 \left( \CuspEnd, L \right)
                \end{array}
            \end{equation}
            called the Chern Laplacian with Alvarez--Wentworth boundary conditions.
        \end{proposition}
        
        \begin{proof}
            The proof follows an argument similar to \cite[Prop. 2.13]{dutour:cuspDirichlet}, using the expression of the domain of the Friedrichs extension. Let us describe part of the argument. Denote by $Q_{\Delta_L}$ the quadratic form associated to $\Delta_L$. Its closure is defined on
            \begin{equation}
            \label{eq:domainClosureQuadForm}
                \begin{array}{lllll}
                    \dom \overline{Q_{\Delta_L}} & = & L_{1, \text{AW}}^2 \left( \CuspEnd, L \right) & = & \left \{ s \in L^2 \left( \CuspEnd, L \right), \; \gamma_1 \left( \left( \Re s \right) s_{\alpha} \right) = 0 \right \}
                \end{array},
            \end{equation}
            where $\gamma_1$ is the boundary trace operator \eqref{boundaryTraceOperator}. This space will be important later. It should be noted that the symmetry of $\Delta_L$ is a consequence of the integration by parts formula \cite[Eq. (2.13)]{wentworth:gluingFormulasAlvarez}.
        \end{proof}
        As explained in \cite[Sec. 2.3]{dutour:cuspDirichlet}, the eigenvalues of $\Delta_L$ are only well-behaved when the line bundle $L$ is (metrically) non-trivial, \textit{i.e.} when we have $\alpha \neq 0$. In the case of the (metrically) trivial line bundle, we need to modify $\Delta_L$ to study its determinant. As we did in \eqref{identification:section:function}, a smooth function $f$ on $\CuspEnd$ can be identified to a smooth function on $\mathbb{R} \times \left] a, +\infty \right[$ which is $1$-periodic in the first variable. Such a function can be written as sum of its Fourier series
        \begin{equation}
            \begin{array}{lll}
                f \left( x, y \right) & = & a_0 \left( y \right) + \sum\limits_{k \neq 0} \; a_k \left( y \right) e^{2 i \pi k x}
            \end{array}.
        \end{equation}
        The hyperbolic laplacian of $f$ is then given by
        \begin{equation}
            \begin{array}{lll}
                \Delta f & = & \left( -y^2 \frac{\mathrm{d}^2}{\mathrm{d}y^2} \right) a_0 \left( y \right) + \left( -y^2 \left( \frac{\partial^2}{\partial x^2} + \frac{\partial^2}{\partial y^2} \right) \right) \left( \sum\limits_{k \neq 0} \; a_k \left( y \right) e^{2 i \pi k x} \right)
            \end{array}
        \end{equation}
        and equals the Laplacian $\Delta_L$ of the section of the trivial bundle induced by $f$. We further consider the projection onto the constant term in the Fourier expansion, denoted by
        \begin{equation}
            \begin{array}{ccccc}
                p & : & L^2 \left( \CuspEnd \right) & \longrightarrow & L^2 \left( \left] a, + \infty \right[, L \right) \\[0.5em]
                
                && f & \longmapsto & a_0 \left( y \right)
            \end{array},
        \end{equation}
        which is a surjective map.
        
        \begin{definition}
        \label{def:constantFourierCoeff}
            The kernel of $p$ is denoted by $L^2 \left( \CuspEnd \right)_0$. The Sobolev spaces~corresponding to those considered in \eqref{eq:SobolevAWL22} and \eqref{eq:domainClosureQuadForm} are denoted in a similar way.
        \end{definition}
        
        \begin{proposition}
            We have the following orthogonal decomposition
            \begin{equation}
                \begin{array}{lll}
                    L^2 \left( \CuspEnd \right) & = & L^2 \left( \CuspEnd \right)_0 \oplus L^2 \left( \left] a, +\infty \right[ \right)
                \end{array} .
            \end{equation}
            Furthermore, the Chern Laplacian with Dirichlet boundary condition splits
            \begin{equation}
                \begin{array}{lll}
                    \Delta & = & \Delta \oplus \left( - y^2 \frac{\mathrm{d}^2}{\mathrm{d}y^2} \right)
                \end{array} ,
            \end{equation}
            where the first Laplacian on the right-hand side acts on
            \begin{equation}
                \begin{array}{lll}
                    L^2_{2, \mathrm{AW}} \left( \CuspEnd \right)_0 & = & L^2_{2, \mathrm{AW}} \left( \CuspEnd \right) \cap L^2 \left( \CuspEnd \right)_0
                \end{array} .
            \end{equation}
        \end{proposition}
        
        For the last definition, let us go back to a general flat line bundle $L$.
        
        \begin{definition}
        \label{def:pseudoLaplacianAW}
            The \textit{pseudo-Laplacian with Alvarez--Wentworth boundary condition} $\Delta_{L,0}$ is defined as the Chern Laplacian:
            \begin{itemize}
                \item acting on $L^2_{2, \mathrm{AW}} \left( \CuspEnd, L \right)$, if the character $\chi$ is non-trivial;
                
                \item acting on $L^2_{2, \mathrm{AW}} \left( \CuspEnd \right)_0$ if $\chi$, and thus $L$, is trivial.
            \end{itemize}
        \end{definition}
        
        \begin{remark}
            When $\chi$ is trivial, there is an added condition of vanishing constant Fourier coefficient. This will be important in subsection \ref{SubSec:spectralProblem} to make sense of the \textit{spectral zeta function} of the pseudo-Laplacian.
        \end{remark}

    \subsection{Weyl type law and the spectral zeta function}
    \label{SubSec:WeylLawSpectralZeta}
    
        To define the spectral zeta function of the pseudo-Laplacian introduced in definition \ref{def:pseudoLaplacianAW}, we need an asymptotic control on its eigenvalues, which will be achieved by obtaining a \textit{Weyl type law}, taking the form of an inequality on the \textit{spectral counting function} $N \left( \Delta_{L,0}, \cdot \right)$, as defined in \cite[Def. A.14]{dutour:cuspDirichlet}. Note that this function differs \textit{a priori} from the more classical \textit{eigenvalue counting function}, insofar as it counts the elements in the sequence obtained in the inf-sup theorem \cite[Thm. A.13]{dutour:cuspDirichlet}, and therefore includes a potentially non-empty essential spectrum.
        
        \begin{theorem}[Weyl type law]
        \label{thm:WeylTypeLaw}
            There exists a constant $C > 0$ such that we have
            \begin{equation}
            \label{eq:WeylTypeLaw}
                \begin{array}{lll}
                    N \left( \Delta_{L,0}, \lambda \right) & \leqslant & C \lambda
                \end{array},
            \end{equation}
            for any strictly positive real number $\lambda$.
        \end{theorem}
        
        \begin{corollary}
        \label{cor:WeylTypeLaw}
            The pseudo-Laplacian with Alvarez--Wentworth boundary condition $\Delta_{L,0}$ has no essential spectrum. Its eigenvalues are isolated and have finite multiplicity. Denoting them by $\left( \lambda_j \right)$ in ascending order with multiplicity, we have
            \begin{equation}
            \label{eq:CorWeylTypeLaw}
                \begin{array}{lllll}
                    N \left( \Delta_{L,0}, \lambda \right) & = & \# \left \{ j, \; \lambda_j \leqslant \lambda \right \} & \leqslant & C \lambda
                \end{array} .
            \end{equation}
        \end{corollary}
        
        \begin{proof}
            This is similar to \cite[Cor. 2.21]{dutour:cuspDirichlet}.
        \end{proof}
        
        \begin{proposition-definition}
            For any real number $\mu \geqslant 0$, the function
            \begin{equation}
            \label{eq:spectralZetaFunction}
                \begin{array}{lllll}
                    \zeta_{L, \mu} & : & s & \longmapsto & \sum\limits_j \; \left( \lambda_j + \mu \right)^{-s}
                \end{array}
            \end{equation}
            is holomorphic on the half-plane $\Re s > 1$, and called the spectral zeta function associated to the pseudo-Laplacian with Alvarez--Wentworth boundary condition. The sum in \eqref{eq:spectralZetaFunction} defining the zeta function ranges over $\mathbb{Z}_{\geqslant 0}$ if we have $\mu > 0$, and over $\mathbb{Z}_{\geqslant 1}$ if we have $\mu = 0$, so as to avoid the eigenvalue $0$.
        \end{proposition-definition}
        
        \begin{proof}
            This is the same argument as \cite[Prop.-Def. 2.22]{dutour:cuspDirichlet}.
        \end{proof}
        
        We will now dedicate the rest of this subsection to the proof of theorem \ref{thm:WeylTypeLaw}.

        \subsubsection{Getting rid of the line bundle}
        
            The first step towards provign the Weyl type law stated in theorem \ref{thm:WeylTypeLaw} is to find an upper-bound of the spectral counting function which does not involve the line bundle $L$, but rather the trivial line bundle. We need to introduce a variation of the Alvarez--Wentworth boundary conditions for the trivial line bundle which corresponds, in the language of Wentworth \cite{wentworth:gluingFormulasAlvarez}, to the boundary conditions for a different \textit{framing}.
            
            \begin{definition}
                Let $k \in \mathbb{Z}$ be an integer. We set
                \begin{equation}
                    \begin{array}{lllll}
                        s_k & : & \mathbb{R} \times \left] a, +\infty \right[ & \longrightarrow & \mathbb{C} \\
                        && \left( x, y \right) & \longmapsto & e^{2 i k \pi \left( x + i y \right)}
                    \end{array}
                \end{equation}
                which induces a smooth function on $\CuspEnd$.
            \end{definition}
            
            \begin{remark}
                Every function $s_k$ is a trivialization of the trivial line bundle over $\CuspEnd$.
            \end{remark}
            
            \begin{definition}
            \label{def:AWFramingsk}
                A function $f$ on $\CuspEnd$ is said to satisfy the \textit{Alvarez--Wentworth boundary conditions relative to the framing} $s_k$ if we have
                \begin{equation}
                    \begin{array}{lll@{\qquad}l@{\qquad}lll}
                        \Im f_k & = & 0 & \text{ and } & \frac{\partial}{\partial y} \Re f_k & = & 0
                    \end{array},
                \end{equation}
                where we wrote $f = f_k s_k$ and we identify functions on $\CuspEnd$ to $1$-periodic in the first variable functions on $\mathbb{R} \times \left] a,  +\infty \right[$.
            \end{definition}
            
            \begin{definition}
                For any reals $a > 0$ and $\lambda > 0$, we denote by $N_{a,k} \left( \lambda \right)$ the spectral counting function of the pseudo-laplacian with Alvarez--Wentworth boundary condition for the (metrically) trivial line bundle over $\CuspidalEnd{a}$, relative to the framing $s_k$.
            \end{definition}
        
            \begin{lemma}
            \label{lem:getRidOfLineBundle}
                Assume the character $\chi$ has finite order $n$. For any $\lambda > 0$, we have
                \begin{equation}
                    \begin{array}{lll}
                        N \left( \Delta_{L,0} \right) & \leqslant & N_{a/n, n \alpha} \left( \lambda \right)
                    \end{array},
                \end{equation}
                where $a>0$ is fixed.
            \end{lemma}
            
            \begin{proof}
                The argument follows the one used in \cite[Lem. 2.24]{dutour:cuspDirichlet}. Nevertheless, let us present the argument, to highlight the importance of the framing. First, we note that the result is empty if we have $\alpha = 0$. We thus assume that we have $\alpha \neq 0$. The character $\chi$ having finite order $n$, we have $n \alpha \in \mathbb{Z}$. Using \eqref{identification:section:function}, we set
                \begin{equation}
                    \begin{array}{lllllll}
                        m_n & : & \mathcal{C}^{\infty} \left( \CuspEnd, L \right) & \longrightarrow & \multicolumn{3}{l}{\mathcal{C}^{\infty} \left( \CuspidalEnd{a/n} \right)} \\[0.5em]
                        && s & \longmapsto & \left( x, y \right) & \longmapsto & s \left( nx, ny \right)
                    \end{array},
                \end{equation}
                which is injective. Consider a smooth section $s$ of $L$ over $\CuspEnd$ and write $s = f s_{\alpha}$. Assume $s$ satisfies the Alvarez--Wentworth boundary condition, \textit{i.e.} that we have
                \begin{equation}
                    \begin{array}{lllllll}
                        \left( \Im s \right) \left( x, a \right) & = & 0 & \text{ and } & \left( \frac{\partial}{\partial y} \Re f \right) \left( x, a \right) & = & 0
                    \end{array}.
                \end{equation}
                Taking the image of $s$ by $m_n$, we get
                \begin{equation}
                    \begin{array}{lllll}
                        m_n \left( s \right) & = & f \left( nx, ny \right) s_{\alpha} \left( nx, ny \right) & = & f \left( nx, ny \right) s_{n \alpha} \left( x, y \right)
                    \end{array}
                \end{equation}
                and we have
                \begin{equation}
                    \begin{array}{lllllll}
                        f \left( nx, n \frac{a}{n} \right) & = & 0 & \text{ and } & n \left( \frac{\partial}{\partial y} f \right) \left( nx, n \frac{a}{n} \right) & = & 0
                    \end{array}.
                \end{equation}
                Thus, the function $m_n \left( s \right)$ satisfies the Alvarez--Wentworth boundary conditions relative to the framing $s_{n \alpha}$. The rest of the argument is as in \cite[Lem. 2.24]{dutour:cuspDirichlet}.
            \end{proof}

        \subsubsection{Relation to the Neumann Laplacian}
        
            In this part of the argument, we will prove theorem \ref{thm:WeylTypeLaw} by relating the spectral counting function therein to that of the Neumann Laplacian, and using the argument detailed in \cite[Sec. 2.4]{dutour:cuspDirichlet}.
            
            \begin{definition}
                Let $\Lambda$ be a union of open intervals contained in $\left] a, + \infty \right[$. The \textit{pseudo-Laplacian with Neumann boundary conditions} $\Delta^N_{\Lambda}$ is defined as the Laplacian associated to the closed quadratic form
                \begin{equation}
                \label{eq:quadFormLaplacianDirichlet}
                    \begin{array}{lllll}
                        Q_N & : & L_1^2 \left( S^1 \times \Lambda \right)_0 \times L_1^2 \left( S^1 \times \Lambda \right)_0 & \longrightarrow & \mathbb{C} \\[0.5em]
                    
                        && \left( u,v \right) & \longmapsto & \int_{S^1 \times \Lambda} \; \nabla u \wedge \overline{\nabla v}
                    \end{array}
                \end{equation}
                on the Sobolev space $L_1^2$ with vanishing constant Fourier coefficient. Here $\nabla$ is the gradient, \textit{i.e.} the Chern connection for the trivial line bundle.
            \end{definition}
            
            \begin{definition}
                Let $\Lambda$ be a union of open intervals contained in $\left] a, + \infty \right[$ and $k \geqslant 0$ be an integer. The \textit{pseudo-Laplacian with Alvarez--Wentworth boundary conditions} is defined as the Laplacian $\Delta^{AW}_{k, \Lambda}$ associated to the closed quadratic form
                \begin{equation}
                \label{eq:quadFormLaplacianAW}
                    \begin{array}{lllll}
                        Q_{AW,k} & : & L_1^2 \left( S^1 \times \Lambda \right)_{AW,0} \times L_1^2 \left( S^1 \times \Lambda \right)_{AW,0} & \longrightarrow & \mathbb{C} \\[0.5em]
                    
                        && \left( u,v \right) & \longmapsto & \int_{S^1 \times \Lambda} \; \overline{\partial} u \wedge \partial v
                    \end{array}
                \end{equation}
                on the Sobolev space $L_1^2$ with vanishing constant Fourier coefficient and Alvarez--Wentworth boundary condition relative to $s_k$ from definition \ref{def:AWFramingsk}.
            \end{definition}
            
            We will now compare these two Laplacians for the order on self-adjoint operators introduced in \cite[Def. A.16]{dutour:cuspDirichlet} and deduce from that a comparison of the associated spectral counting functions, using \cite[Prop. A.20]{dutour:cuspDirichlet}.
            
            \begin{lemma}
            \label{lem:compareNeumannAW}
                Let $\Lambda$ be an open interval in $\left] a, + \infty \right[$. We have $\Delta_{\Lambda}^N \preccurlyeq \Delta_{\Lambda}^{AW}$. Consequently, we have $N \left( \Delta_{\Lambda}^{AW}, \lambda \right) \leqslant N \left( \Delta_{\Lambda}^N, \lambda \right)$ for any $\lambda > 0$.
            \end{lemma}
            
            \begin{proof}
                This result stems from a comparison of the domain of the quadratic forms defined in \eqref{eq:quadFormLaplacianDirichlet} and \eqref{eq:quadFormLaplacianAW}, together with the observation that we have
                \begin{equation}
                    \begin{array}{lll}
                        \left \Vert \overline{\partial} u \right \Vert^2 & \leqslant & \left \Vert \nabla u \right \Vert^2
                    \end{array}
                \end{equation}
                for any smooth complex valued function $u$ on $\CuspEnd$, with respect to the $L^2$-product, the inequality steming from the orthogonality of $\left( 0,1 \right)$-forms and $\left( 1,0 \right)$-forms.
            \end{proof}

        \subsubsection{Proof of the Weyl type estimate}    
            We now have all the required tools to prove theorem \ref{thm:WeylTypeLaw}. We will rely on \cite[Sec. 2.4]{dutour:cuspDirichlet}, which itself is based on an argument by Colin de Verdière from \cite{colinDeVerdiere:pseudoLaplacian2}.
            
            \begin{proof}[Proof of theorem \ref{thm:WeylTypeLaw}]
                Consider a real number $\lambda > 0$. We first note that we have
                \begin{equation}
                    \begin{array}{lll}
                        N_{a/n, n \alpha} \left( \lambda \right) & = & N \left( \Delta_{n\alpha, \left] a/n, + \infty \right[}^{AW}, \lambda \right)
                    \end{array}
                \end{equation}
                and we deduce from this observation the inequalities
                \begin{equation}
                    \begin{array}{llll}
                        N \left( \Delta_{L,0}, \lambda \right) & \leqslant & N \left( \Delta_{n\alpha, \left] a/n, + \infty \right[}^{AW}, \lambda \right) & \text{using lemma \ref{lem:getRidOfLineBundle}} \\[1em]
                        
                        & \leqslant & N \left( \Delta_{\Lambda}^N, \lambda \right) & \text{using lemma \ref{lem:compareNeumannAW}}.
                    \end{array}
                \end{equation}
                Hence, it is enough to prove the analogue of the Weyl type law for the Neumann Laplacian. This is done in \cite[Sec. 2.4]{dutour:cuspDirichlet}.
            \end{proof}

    \subsection{Spectral problem}
    \label{SubSec:spectralProblem}
    
        Similarly to \cite[Sec. 2.5]{dutour:cuspDirichlet}, we will now state and study a spectral problem related to the pseudo-Laplacian with Alvarez--Wentworth boundary condition (sometimes shortened as ``AW''). It is given by:
        
        \begin{equation}
        \label{eq:spectralProblem1}
            \left \{ \begin{array}{llll}
                - y^2 \left( \frac{\partial^2}{\partial x^2} + \frac{\partial^2}{\partial y^2} \right) \psi \left( x, y \right) & = & \multicolumn{2}{l}{\lambda \psi \left( x, y \right)} \\[0.5em]
                
                \psi \left( x+1, y \right) & = & \multicolumn{2}{l}{e^{2 i \pi \alpha} \psi \left( x, y \right)} \\[0.5em]
                
                \int_{S^1 \times \left] a, +\infty \right[} \; \left \vert \psi \right \vert^2 & < & + \infty & \text{(Integrability condition)} \\[0.7em]
                
                \Im \left( \psi s_{- \alpha} \right) \left( x, a \right) & = & 0 & \multirow{2}{*}{\text{(AW boundary condition)}} \\[0.5em]
                
                \frac{\partial}{\partial y} \Re \left( \psi s_{- \alpha} \right) \left( x, a \right) & = & 0 \\[0.5em]
                
                \int_{S^1} \; \psi \left( x, y \right) \; \mathrm{d}x & = & 0 &  \text{for almost all } y \geqslant a \text{ if } \alpha = 0
            \end{array} \right.
        \end{equation}
        
        Using the change of function $\varphi \left( x,y \right) = e^{-2 i \pi \alpha x} \psi \left( x, y \right)$, the spectral problem given in \eqref{eq:spectralProblem1} is equivalent to
        
        \begin{equation}
        \label{eq:spectralProblem2}
            \left \{ \begin{array}{llll}
                - y^2 \left( \frac{\partial^2}{\partial x^2} + \frac{\partial^2}{\partial y^2} \right) \varphi \left( x, y \right) & = & \multicolumn{2}{l}{\left( \lambda - 4 \pi^2 \alpha^2 \right) \varphi \left( x, y \right) + 4 i \pi \alpha y^2 \frac{\partial}{\partial x} \varphi} \\[0.5em]
                
                \varphi \left( x+1, y \right) & = & \multicolumn{2}{l}{\varphi \left( x, y \right)} \\[0.5em]
                
                \int_{S^1 \times \left] a, +\infty \right[} \; \left \vert \varphi \right \vert^2 & < & + \infty & \text{(Integrability condition)} \\[0.7em]
                
                \Im \left( \varphi e^{2 \pi \alpha y} \right) \left( x, a \right) & = & 0 & \multirow{2}{*}{\text{(AW boundary condition)}} \\[0.5em]
                
                \frac{\partial}{\partial y} \Re \left( \varphi e^{2 \pi \alpha y} \right) \left( x, a \right) & = & 0 \\[0.5em]
                
                \int_{S^1} \; \varphi \left( x, y \right) \; \mathrm{d}x & = & 0 &  \text{for almost all } y \geqslant a \text{ if } \alpha = 0
            \end{array} \right.
        \end{equation}
        
        To make the Alvarez--Wentworth boundary condition appear more easily, we decompose $\varphi = u + i v$ into its real and imaginary parts in \eqref{eq:spectralProblem2}, which becomes equivalent to
        
        \begin{equation}
        \label{eq:spectralProblem3}
            \left \{ \begin{array}{llll}
                - y^2 \left( \frac{\partial^2}{\partial x^2} + \frac{\partial^2}{\partial y^2} \right) u \left( x, y \right) & = & \multicolumn{2}{l}{\left( \lambda - 4 \pi^2 \alpha^2 \right) u \left( x, y \right) - 4 \pi \alpha y^2 \frac{\partial}{\partial x} v} \\[0.5em]
                
                - y^2 \left( \frac{\partial^2}{\partial x^2} + \frac{\partial^2}{\partial y^2} \right) v \left( x, y \right) & = & \multicolumn{2}{l}{\left( \lambda - 4 \pi^2 \alpha^2 \right) v \left( x, y \right) + 4 \pi \alpha y^2 \frac{\partial}{\partial x} u} \\[0.5em]
                
                u \left( x+1, y \right) & = & \multicolumn{2}{l}{u \left( x, y \right)} \\[0.5em]
                
                v \left( x+1, y \right) & = & \multicolumn{2}{l}{v \left( x, y \right)} \\[0.5em]
                
                \int_{S^1 \times \left] a, +\infty \right[} \; \left \vert u \right \vert^2 & < & + \infty & \multirow{2}{*}{\text{(Integrability conditions)}} \\[0.7em]
                
                \int_{S^1 \times \left] a, +\infty \right[} \; \left \vert v \right \vert^2 & < & + \infty & \\[0.7em]
                
                v \left( x, a \right) & = & 0 & \multirow{2}{*}{\text{(AW boundary condition)}} \\[0.5em]
                
                \frac{\partial u}{\partial y} \left( x, a \right) + 2 \pi \alpha u \left( x, a \right) & = & 0 \\[0.5em]
                
                \int_{S^1} \; u \left( x, y \right) \; \mathrm{d}x & = & 0 &  \multirow{2}{*}{\text{for almost all $y \geqslant a$ if $\alpha = 0$}} \\[0.5em]
                
                \int_{S^1} \; v \left( x, y \right) \; \mathrm{d}x & = & 0 &
            \end{array} \right.
        \end{equation}
        
        The functions $u$ and $v$ are smooth and $1$-periodic in the first variable, and can thus be written as sums of their respective Fourier series. We have
        \begin{equation}
        \label{eq:fourierDecompositions}
            \begin{array}[t]{lllllll}
                u \left( x, y \right) & = & \sum\limits_k \; a_k \left( y \right) e^{2 i k \pi x} & \text{ and } & v \left( x, y \right) & = & \sum\limits_k \; b_k \left( y \right) e^{2 i k \pi x}
            \end{array}
        \end{equation}
        where each sum ranges over $k \in \mathbb{Z}$ if we have $\alpha \neq 0$, and over $k \in \mathbb{Z}^{\ast} = \mathbb{Z} \setminus \left \{ 0 \right \}$ if we have $\alpha = 0$. Using the unicity of Fourier decompositions, the partial differential equations from \eqref{eq:spectralProblem3} become
        \begin{equation}
        \label{eq:odeFourier1}
            \left \{ \begin{array}{llllr}
                y^2 \left( a_k'' \left( y \right) - 4 \pi^2 k^2 a_k \left( y \right) \right) \hspace{-6pt} & + & \hspace{-6pt} \left( \lambda - 4 \pi^2 \alpha^2 y^2 \right) a_k \left( y \right) & \hspace{-6pt} = \hspace{-6pt} & 8 i \pi^2 \alpha k y^2 b_k \left( y \right) \\[0.5em]
                
                y^2 \left( b_k'' \left( y \right) - 4 \pi^2 k^2 b_k \left( y \right) \right) \hspace{-6pt} & + & \hspace{-6pt} \left( \lambda - 4 \pi^2 \alpha^2 y^2 \right) b_k \left( y \right) & \hspace{-6pt} = \hspace{-6pt} & - 8 i \pi^2 \alpha k y^2 a_k \left( y \right)
            \end{array} \right.
        \end{equation}
        and is equivalent to
        \begin{equation}
        \label{eq:odeFourier2}
            \left \{ \begin{array}{llllr}
                y^2 c_k'' \left( y \right) & + & \left( \lambda - 4 \pi^2 \left( k + \alpha \right)^2 y^2 \right) c_k \left( y \right) & = & 0 \\[0.5em]
                
                y^2 d_k'' \left( y \right) & + & \left( \lambda - 4 \pi^2 \left( k + \alpha \right)^2 y^2 \right) d_k \left( y \right) & = & 0
            \end{array} \right.
        \end{equation}
        by setting $c_k \left( y \right) = a_k \left( y \right) + i b_k \left( y \right)$ and $d_k \left( y \right) = a_k \left( y \right) - i b_k \left( y \right)$. Taking into account the integrability conditions from \eqref{eq:spectralProblem3}, we have
        \begin{equation}
        \label{eq:solOdeFourier2}
            \begin{array}{lll}
                c_k \left( y \right) & = & \alpha_k \sqrt{y} K_{s-1/2} \left( 2 \pi \left \vert k + \alpha \right \vert y \right) \\[0.5em]
                d_k \left( y \right) & = & \beta_k \sqrt{y} K_{s-1/2} \left( 2 \pi \left \vert k - \alpha \right \vert y \right)
            \end{array}
        \end{equation}
        for some complex constants $\alpha_k$ and $\beta_k$ depending only on $k$. In \eqref{eq:solOdeFourier2}, we denoted by $K$ the modified Bessel functions of the second kind, and we set $\lambda = s \left( 1-s \right)$. The essential properties of these functions is recalled in \cite[App. C]{dutour:cuspDirichlet}, which is based on \cite{olver, nist}. Returning to $a_k \left( y \right)$ and $b_k \left( y \right)$, the solutions of $\eqref{eq:odeFourier1}$ satisfying the integrability conditions are given by
        \begin{equation}
        \label{eq:solOdeFourier1}
            \begin{array}{lll}
                a_k \left( y \right) & = & \frac{1}{2} \sqrt{y} \left[ \alpha_k K_{s-1/2} \left( 2 \pi \left \vert k + \alpha \right \vert y \right) + \beta_k K_{s-1/2} \left( 2 \pi \left \vert k - \alpha \right \vert y \right) \right] \\[0.7em]
                
                b_k \left( y \right) & = & \frac{1}{2i} \sqrt{y} \left[ \alpha_k K_{s-1/2} \left( 2 \pi \left \vert k + \alpha \right \vert y \right) - \beta_k K_{s-1/2} \left( 2 \pi \left \vert k - \alpha \right \vert y \right) \right]
            \end{array}
        \end{equation}
        We now need to understand which of these solutions satisfy the Alvarez--Wentworth boundary conditions extracted from \eqref{eq:spectralProblem3} on the Fourier coefficients, given by
        \begin{equation}
        \label{eq:AWBoundaryCondition1}
            \begin{array}{lllllll}
                b_k \left( a \right) & = 0 & \text{ and } & a_k' \left( a \right) + 2 \pi \alpha a_k \left( a \right) & = & 0
            \end{array}.
        \end{equation}
        Using \eqref{eq:solOdeFourier1}, we can rewrite \eqref{eq:AWBoundaryCondition1} as
        \begin{equation}
        \label{eq:AWBoundaryCondition2}
            \left \{ \begin{array}{lll}
                \alpha_k K_{s-1/2} \left( 2 \pi \left \vert k + \alpha \right \vert a \right) + \beta_k \left( - K_{s-1/2} \left( 2 \pi \left \vert k - \alpha \right \vert a \right) \right) & = & 0 \\[1em]
                
                \alpha_k \left[ \left( 1 + 4 \pi \alpha a \right) K_{s-1/2} \left( 2 \pi \left \vert k + \alpha \right \vert a \right) \right. \\[0.3em]
                \qquad \qquad \qquad \qquad \qquad + \left. 4 \pi \left \vert k + \alpha \right \vert a K_{s-1/2}' \left( 2 \pi \left \vert k + \alpha \right \vert a \right) \right] \\[0.3em]
                + \beta_k \left[ \left( 1 + 4 \pi \alpha a \right) K_{s-1/2} \left( 2 \pi \left \vert k - \alpha \right \vert a \right) \right. \\[0.3em]
                \qquad \qquad \qquad \qquad \qquad + \left. 4 \pi \left \vert k - \alpha \right \vert a K_{s-1/2}' \left( 2 \pi \left \vert k - \alpha \right \vert a \right) \right] & = & 0
            \end{array} \right.
        \end{equation}
        The system \eqref{eq:AWBoundaryCondition2} has a non-trivial solution if and only if its determinant vanishes, meaning if we have
        \begin{equation}
        \label{eq:AWBoundaryCondition3}
            \begin{array}{lll}
                K_{s-1/2} \left( 2 \pi \left \vert k + \alpha \right \vert a \right) \left[ \left( 1 + 4 \pi \alpha a \right) K_{s-1/2} \left( 2 \pi \left \vert k - \alpha \right \vert a \right) \right. \\
                \qquad \qquad \qquad \qquad \qquad \qquad \left. + 4 \pi \left \vert k - \alpha \right \vert a K_{s-1/2}' \left( 2 \pi \left \vert k - \alpha \right \vert a \right) \right] \\
                + K_{s-1/2} \left( 2 \pi \left \vert k - \alpha \right \vert a \right) \left[ \left( 1 + 4 \pi \alpha a \right) K_{s-1/2} \left( 2 \pi \left \vert k + \alpha \right \vert a \right) \right. \\
                \qquad \qquad \qquad \qquad \qquad \qquad \left. + 4 \pi \left \vert k - \alpha \right \vert a K_{s-1/2}' \left( 2 \pi \left \vert k - \alpha \right \vert a \right) \right] & = & 0.
            \end{array}
        \end{equation}
        This equation can be rewritten as
        \begin{equation}
        \label{eq:AWBoundaryCondition4}
            \begin{array}{lll}
                \left( 1 + 4 \pi \alpha a \right) K_{s-1/2} \left( 2 \pi \left \vert k + \alpha \right \vert a \right) K_{s-1/2} \left( 2 \pi \left \vert k - \alpha \right \vert a \right) \\[0.5em]
                \quad + 2 \pi \left \vert k - \alpha \right \vert a K_{s-1/2} \left( 2 \pi \left \vert k + \alpha \right \vert a \right) K_{s-1/2} ' \left( 2 \pi \left \vert k - \alpha \right \vert a \right) \\[0.5em]
                \quad \quad \quad + 2 \pi \left \vert k + \alpha \right \vert a K_{s-1/2} \left( 2 \pi \left \vert k - \alpha \right \vert a \right) K_{s-1/2} ' \left( 2 \pi \left \vert k + \alpha \right \vert a \right) & = & 0.
            \end{array}
        \end{equation}
        Recall that we previously set $\lambda = s \left( 1-s \right)$. Hence, solving \eqref{eq:AWBoundaryCondition4} in $s$ is equivalent to finding the eigenvalues of the pseudo-Laplacian with Alvarez--Wentworth boundary conditions, which were initally characterized by the spectral problem \eqref{eq:spectralProblem1}.

    \subsection{Localization of the eigenvalues}
    \label{SubSec:localizationEigenvalues}
        
        In this subsection, we will find more information on the solutions in $s$ of equation \eqref{eq:AWBoundaryCondition4}, in order to localize the eigenvalues of the pseudo-Laplacian with Alvarez--Wentworth boundary conditions. The argument below is similar to the one used in \cite[Prop. C.9]{dutour:cuspDirichlet}, and is also based on an argument presented by Saharian in \cite[App. 9]{saharian}. However, the case $k=0$, which must be considered if we have $\alpha \neq 0$, involves some substantial changes. As such, the localization of eigenvalues in this case, and the subsequent computation of determinants, will be performed separately from the case $k \neq 0$.
        
        \begin{definition}
            For any non-zero integer $k \in \mathbb{Z} \setminus \left \{ 0 \right \}$, we set
            \begin{equation}
            \label{eq:functionf}
                \begin{array}{lll}
                    f_k \left( \nu \right) & = & K_{i\nu} \left( 2 \pi \left \vert k + \alpha \right \vert a \right) K_{i\nu} \left( 2 \pi \left \vert k - \alpha \right \vert a \right) g_k \left( \nu \right)
                \end{array}
            \end{equation}
            where the function $g_k$ is defined by
            \begin{equation}
            \label{eq:functiong}
                \begin{array}{lll}
                    g_k \left( \nu \right) & = & 1 + 4 \pi \alpha a + 2 \pi \left \vert k + \alpha \right \vert a \frac{K_{i\nu}' \left( 2 \pi \left \vert k + \alpha \right \vert a \right)}{K_{i\nu} \left( 2 \pi \left \vert k + \alpha \right \vert a \right)} \\[0.5em]
                    
                    && \qquad \qquad \qquad \qquad \qquad \qquad \qquad + 2 \pi \left \vert k - \alpha \right \vert a \frac{K_{i\nu}' \left( 2 \pi \left \vert k - \alpha \right \vert a \right)}{K_{i\nu} \left( 2 \pi \left \vert k - \alpha \right \vert a \right)}.
                \end{array}
            \end{equation}
        \end{definition}
        
        \begin{remark}
            Using the properties of modified Bessel functions (see \cite[App. C]{dutour:cuspDirichlet}), the function $f$ is entire, and its zeros characterize the solutions of \eqref{eq:AWBoundaryCondition4}.
        \end{remark}
        
        \begin{proposition}
        \label{prop:zerosfk}
            Assume we have $a > 1/\left( 4 \pi \left( 1 - \alpha \right) \right)$. Then, for any non-zero integer $k \in \mathbb{Z} \setminus \left \{ 0 \right \}$, the zeros of $f_k$ are real and simple.
        \end{proposition}
        
        \begin{proof}
            Consider a zero $\nu$ of $f_k$, which we assume is neither real nor purely imaginary. Using \cite[Prop. C.9]{dutour:cuspDirichlet}, we have $g_k \left( \nu \right) = 0$, since $\nu$ cannot be a zero of the modified Bessel functions in factor of $g_k$ in \eqref{eq:functionf}. This result stems from the formula
            \begin{equation}
            \label{eq:intBesselZero1}
                \begin{array}{lll}
                    \int_0^u \; \frac{1}{v} \left \vert K_{i \nu} \left( v \right) \right \vert^2 \; \mathrm{d}v & = & - \frac{u}{\overline{\nu}^2 - \nu^2} \left[ K_{i \nu} \left( u \right) K_{i \overline{\nu}}' \left( u \right) -  K_{i \overline{\nu}} \left( u \right) K_{i \nu}' \left( u \right) \right]
                \end{array}
            \end{equation}
            for any real number $u>0$. We then have
            \begin{equation}
            \label{eq:intBesselZero2}
                \begin{array}{lll}
                    \multicolumn{2}{l}{\frac{1}{\left \vert K_{i\nu} \left( 2 \pi \left \vert k + \alpha \right \vert a \right) \right \vert^2} \int_0^{2 \pi \left \vert k + \alpha \right \vert a} \frac{1}{v} \left \vert K_{i \nu} \left( v \right) \right \vert^2 \; \mathrm{d}v} \\[1em]
                    \multicolumn{2}{l}{\qquad \qquad \qquad + \frac{1}{\left \vert K_{i\nu} \left( 2 \pi \left \vert k - \alpha \right \vert a \right) \right \vert^2} \int_0^{2 \pi \left \vert k - \alpha \right \vert a} \frac{1}{v} \left \vert K_{i \nu} \left( v \right) \right \vert^2 \; \mathrm{d}v} \\[1em]
                    
                    = & - \frac{2 \pi a}{\overline{\nu}^2 - \nu^2} \left[ \left \vert k + \alpha \right \vert \frac{K_{i\overline{\nu}}' \left( 2 \pi \left \vert k + \alpha \right \vert a \right)}{K_{i\overline{\nu}} \left( 2 \pi \left \vert k + \alpha \right \vert a \right)} + \left \vert k - \alpha \right \vert \frac{K_{i\overline{\nu}}' \left( 2 \pi \left \vert k - \alpha \right \vert a \right)}{K_{i\overline{\nu}} \left( 2 \pi \left \vert k - \alpha \right \vert a \right)} \right] \\[1em]
                    & + \frac{2 \pi a}{\overline{\nu}^2 - \nu^2} \left[ \left \vert k + \alpha \right \vert \frac{K_{i\nu}' \left( 2 \pi \left \vert k + \alpha \right \vert a \right)}{K_{i\nu} \left( 2 \pi \left \vert k + \alpha \right \vert a \right)} + \left \vert k - \alpha \right \vert \frac{K_{i\nu}' \left( 2 \pi \left \vert k - \alpha \right \vert a \right)}{K_{i\nu} \left( 2 \pi \left \vert k - \alpha \right \vert a \right)} \right] & \text{using \eqref{eq:intBesselZero1}} \\[1em]
                    
                    = & 0.
                \end{array}
            \end{equation}
            For the last equality, let us note that the Schwarz reflection principle and properties of the modified Bessel functions imply that both $\nu$ and $\overline{\nu}$ are zeros of $g_k$. For details on this symmetry, the reader is referred to the proof of \cite[Prop. C.9]{dutour:cuspDirichlet}. In turn, since both terms on the left-hand side of \eqref{eq:intBesselZero2} are positive, we find a contradiction, as modified Bessel functions cannot vanish identically. Therefore, the zero $\nu$ of $g_k$ is either real or purely imaginary. Let us assume that $\nu$ is purely imaginary, and write $\nu = - it$ with $t$ real. According to \cite[7.8.1.(8.03)]{olver}, we have
            \begin{equation}
                \begin{array}{lll}
                    K_t \left( x \right) & = & \int_0^{+\infty} e^{-x \cosh u} \cosh \left( t u \right) \; \mathrm{d}u
                \end{array}
            \end{equation}
            for any real number $x > 0$, and we can, on the same interval, differentiate under the integral with respect to $x$. Applying these observations to $f_k$, we have
            \begin{equation}
                \begin{array}{lllll}
                    f_k \left( - i t \right) & = & \int_0^{+ \infty} \int_0^{+\infty} \left( e^{-2 \pi \left \vert k + \alpha \right \vert a \cosh u} e^{-2 \pi \left \vert k - \alpha \right \vert a \cosh v} \cosh \left( t u \right)  \right. \\[0.5em]
                    
                    && \qquad \qquad \qquad \times \cosh \left( t v \right) \left[ 1 + 4 \pi \alpha a - 2 \pi \left \vert k + \alpha \right \vert a \cosh u \right. \\[0.5em]
                    
                    && \qquad \qquad \qquad \qquad \left. \left. - 2 \pi \left \vert k - \alpha \right \vert a \cosh v \right] \right) \mathrm{d}u \mathrm{d}v.
                \end{array}
            \end{equation}
            We now note that we have
            \begin{equation}
            \label{eq:upperBoundIntegralRepModBessel}
                \begin{array}{lllllll}
                    \multicolumn{7}{l}{1 + 4 \pi \alpha a - 2 \pi \left \vert k + \alpha \right \vert a \cosh u - 2 \pi \left \vert k - \alpha \right \vert a \cosh v} \\[0.5em]
                    
                    \quad \quad & \leqslant & 1 - 2 \pi a \left[ \left \vert k + \alpha \right \vert + \left \vert k - \alpha \right \vert - 2 \alpha \right] & \leqslant & 1 - 4 \pi a \left( 1 - \alpha \right) & < & 0,
                \end{array}
            \end{equation}
            having assumed that we have $k \neq 0$ and $a > 1/\left( 4 \pi \left( 1 - \alpha \right) \right)$. Thus, the zero $\nu = - it$ of $f_k$ cannot be purely imaginary, and must therefore be real. Let us now prove that such a real zero $\nu$ can only be simple. Recall that we have
            \begin{equation}
                \begin{array}{lll}
                    f_k \left( \nu \right) & = & \left( 1 + 4 \pi \alpha a \right) K_{i\nu} \left( 2 \pi \left \vert k + \alpha \right \vert a \right) K_{i\nu} \left( 2 \pi \left \vert k - \alpha \right \vert a \right) \\[0.5em]
                    
                    && \qquad + 2 \pi \left \vert k + \alpha \right \vert a K_{i\nu}' \left( 2 \pi \left \vert k + \alpha \right \vert a \right) K_{i\nu} \left( 2 \pi \left \vert k - \alpha \right \vert a \right) \\[0.5em]
                    
                    && \qquad \qquad + 2 \pi \left \vert k - \alpha \right \vert a K_{i\nu}' \left( 2 \pi \left \vert k - \alpha \right \vert a \right) K_{i\nu} \left( 2 \pi \left \vert k + \alpha \right \vert a \right).
                \end{array}
            \end{equation}
            If $\nu$ satisfies $K_{i\nu} \left( 2 \pi \left \vert k + \alpha \right \vert a \right) = 0$, then we also have $K_{i\nu} \left( 2 \pi \left \vert k - \alpha \right \vert a \right) = 0$, since we cannot have $K_{i\nu}' \left( 2 \pi \left \vert k + \alpha \right \vert a \right) = 0$, as modified Bessel functions satisfy a second order linear differential equation and do not vanish identically. In this case, we can use \cite[Prop. C.9]{dutour:cuspDirichlet} to see that $\nu$ is a simple zero of both factors of $f_k$ in \eqref{eq:functionf}. It is also a simple pole of $g_k$, as it involves logarithmic derivatives whose residues do not cancel each other. We thus get a simple zero of $f_k$. We can now move on to the case where $\nu$ is a (real) zero of $g_k$. As indicated in the proof of \cite[Prop. C.9]{dutour:cuspDirichlet}, a difference quotient argument applied to \eqref{eq:intBesselZero1} yields
            \begin{equation}
            \label{eq:intModifiedBesselRealOrder}
                \begin{array}{lll}
                    \multicolumn{3}{l}{\int_0^u \; \frac{1}{v} \left \vert K_{i \nu} \left( v \right) \right \vert^2 \; \mathrm{d}v} \\[0.5em]
                    \qquad \qquad \quad & = & - i \frac{u}{2 \nu} \left[ K_{i \nu} \left( u \right) \frac{\partial}{\partial \beta}_{\vert \beta = i \nu} K_{\beta}' \left( u \right) - K_{i \nu}' \left( u \right) \frac{\partial}{\partial \beta}_{\vert \beta = i \nu} K_{\beta} \left( u \right) \right]
                \end{array}
            \end{equation}
            for any real number $u>0$. We can apply \eqref{eq:intModifiedBesselRealOrder} to get
            \begin{equation}
            \label{eq:orderDerivativeFunctiong}
                \begin{array}{lll}
                    g_k' \left( \nu \right) & = & - \frac{2 \nu}{\left \vert K_{i \nu} \left( 2 \pi \left \vert k + \alpha \right \vert a \right) \right \vert^2} \int_0^{2 \pi \left \vert k + \alpha \right \vert a} \; \frac{1}{v} \left \vert K_{i \nu} \left( v \right) \right \vert^2 \; \mathrm{d}v \\[0.5em]
                    
                    && \qquad \qquad \qquad - \frac{2 \nu}{\left \vert K_{i \nu} \left( 2 \pi \left \vert k - \alpha \right \vert a \right) \right \vert^2} \int_0^{2 \pi \left \vert k - \alpha \right \vert a} \; \frac{1}{v} \left \vert K_{i \nu} \left( v \right) \right \vert^2 \; \mathrm{d}v.
                \end{array}
            \end{equation}
            Note that the complex moduli around the modified Bessel functions can actually be removed here, as they take positive real values when the order is purely imaginary and the argument is real. Since these functions further do not vanish identically, we have $g_k' \left( \nu \right) \neq 0$, and $\nu$ is thus a simple (real) zero of $g_k$.
        \end{proof}
        
        \begin{remark}
            The reason we had to remove $k=0$ in this proposition and defer its treatment to a later proposition in this subsection was that \eqref{eq:upperBoundIntegralRepModBessel} leads to an upper bound of $1$ there, which cannot be used to get the absence of purely imaginary zeros of $f_k$. In fact, there are two such zeros in that case, located at $\nu = \pm i/2$. The argument will be modified to account for these differences.
        \end{remark}
        
        \begin{corollary}
        \label{cor:locZerosKNeq0}
            For any non-zero integer $k \in \mathbb{Z} \setminus \left \{ 0 \right \}$, the function
            \begin{equation}
                \begin{array}{lll}
                    s & \longmapsto & \left( 1 + 4 \pi \alpha a \right) K_{s-1/2} \left( 2 \pi \left \vert k + \alpha \right \vert a \right) K_{s-1/2} \left( 2 \pi \left \vert k - \alpha \right \vert a \right) \\[0.5em]
                    
                    && \quad + 2 \pi \left \vert k - \alpha \right \vert a K_{s-1/2} \left( 2 \pi \left \vert k + \alpha \right \vert a \right) K_{s-1/2} ' \left( 2 \pi \left \vert k - \alpha \right \vert a \right) \\[0.5em]
                    
                    && \quad \quad \quad + 2 \pi \left \vert k + \alpha \right \vert a K_{s-1/2} \left( 2 \pi \left \vert k - \alpha \right \vert a \right) K_{s-1/2} ' \left( 2 \pi \left \vert k + \alpha \right \vert a \right)
                \end{array}
            \end{equation}
            is entire. Assuming we have $a > 1/\left( 4 \pi \left( 1-\alpha \right) \right)$, its zeros, which are all simple, are given by a discrete set
            \begin{equation}
            \label{eq:zerosFK}
                \begin{array}{lll}
                    \left \{ \frac{1}{2} + i r_{k,j} \right \} & \subset & \left \{ \frac{1}{2} + ir, \; r \in \mathbb{R}^{\ast} \right \}.
                \end{array}
            \end{equation}
            The corresponding eigenvalues of the spectral problem \eqref{eq:spectralProblem1} are given by
            \begin{equation}
                \begin{array}{llll}
                    \lambda_{k,j} & = & \frac{1}{4} + r_{k,j}^2 & \text{for } k \neq 0.
                \end{array}
            \end{equation}
        \end{corollary}
        
        \begin{proof}
            This is a consequence of proposition \ref{prop:zerosfk}.
        \end{proof}
        
        \begin{remark}
            Note that \eqref{eq:spectralProblem1} has other eigenvalues if we have $\alpha \neq 0$, corresponding to $k=0$ in the Fourier decompositions \eqref{eq:fourierDecompositions}. They are described below.
        \end{remark}
        
        \begin{proposition}
        \label{prop:zerosFunctionF0}
            Assume we have $\alpha \neq 0$. The function $f_0$ is entire, its zeros are all simple and located on the real line, with the exception of two (simple zeros) located at $\nu = \pm i/2$.
        \end{proposition}
        
        \begin{proof}
            The fact that $f_0$ is entire follows from the properties of the modified Bessel functions of the second kind. Furthermore, it cannot vanish at any point $\nu$ which is neither real nor purely imaginary, using \eqref{eq:intBesselZero2} and \cite[Prop. C.9]{dutour:cuspDirichlet}. Let us now prove that the only purely imaginary zeros of $f$ are $\nu = \pm i/2$. We have
            \begin{equation}
                \begin{array}{lll}
                    f_0 \left( \nu \right) & = & \left( 1 + 4 \pi \alpha a \right) K_{i \nu} \left( 2 \pi \alpha a \right)^2 + 4 \pi \alpha a K_{i \nu} ' \left( 2 \pi \alpha a \right) K_{i \nu} \left( 2 \pi \alpha a \right).
                \end{array}
            \end{equation}
            Using \cite[Prop. C.9]{dutour:cuspDirichlet}, the function $\nu \mapsto K_{i \nu} \left( 2 \pi \alpha a \right)$ is entire, and all of its zeros are real simple. Furthermore, these zeros must be distinct from those of
            \begin{equation}
                \begin{array}{lll}
                    \nu & \longmapsto & \left( 1 + 4 \pi \alpha a \right) K_{i \nu} \left( 2 \pi \alpha a \right) + 4 \pi \alpha a K_{i \nu} ' \left( 2 \pi \alpha a \right)
                \end{array}
            \end{equation}
            as a zero of these two functions would satisfy $K_{i \nu} \left( 2 \pi \alpha a \right) = K_{i \nu} ' \left( 2 \pi \alpha a \right) = 0$, which would be absurd since modified Bessel functions of the second kind do not vanish everywhere, and satisfy a second order linear differential equation. Let us therefore prove that for any fixed real number $t > 0$, the function
            \begin{equation}
                \begin{array}{lllll}
                    h & : & \mathbb{R}_+^{\ast} & \longrightarrow & \mathbb{R} \\
                    && x & \longmapsto & \left( 1 + 2x \right) K_t \left( x \right) + 2x K_t' \left( x \right)
                \end{array}
            \end{equation}
            cannot vanish anywhere assuming we have $t \neq 1/2$. First, we note that we have
            \begin{equation}
            \label{eq:asymptHK0}
                \begin{array}{llll}
                    h \left( x \right) & = & \sqrt{\frac{\pi}{2x}} e^{-x} \left( \frac{3}{8} \left( 1 - 4t^2 \right) + O \left( \frac{1}{x} \right) \right) & \text{as } x \rightarrow + \infty.
                \end{array}
            \end{equation}
            This asymptotic expansion is a consequence of a similar one on $K_t \left( x \right)$ and $K_t' \left( x \right)$. For instance, they can be obtained from the integral representation \cite[Ex. 8.4]{olver}. Thus $h$ converges to $0$ at infinity and, for $x > 0$ large enough, it is strictly positive if we have $t>1/2$, and strictly negative if we have $t<1/2$. Let us assume that $h$ has at least one zero, and denote by $x_0$ the largest zero. Recall that we have
            \begin{equation}
            \label{eq:ODEModifiedBessel}
                \begin{array}{lll}
                    x^2 K_t'' \left( x \right) + x K_t' \left( x \right) - \left( x^2 + t^2 \right) K_t \left( x \right) & = & 0.
                \end{array}
            \end{equation}
            A direct computation using \eqref{eq:ODEModifiedBessel} now shows that we have
            \begin{equation}
            \label{eq:ODEFunctionH}
                \begin{array}{lllll}
                    x^2 \, \frac{\mathrm{d}^2h}{\mathrm{d}x^2} + x \, \frac{\mathrm{d} h}{\mathrm{d}x} & = & x \, \frac{\mathrm{d}}{\mathrm{d}x} \left( x \, \frac{\mathrm{d} h}{\mathrm{d}x} \right) & = & \left( 2x + x^2 + t^2 \right) h \left( x \right).
                \end{array}
            \end{equation}
            Using Rolle's theorem, there exists $y > x_0$ such that we have $h' \left( y \right) = 0$. Since $h$ is of constant sign after $x_0$, we may assume that $y$ is either a local minimum of $h$ if we have $t < 1/2$, or a local maximum of $h$ if we have $t > 1/2$. Thus $h'' \left( y \right)$ has the opposite sign to $h \left( y \right)$. Evaluating \eqref{eq:ODEFunctionH} at $y$, we get
            \begin{equation}
                \begin{array}{lll}
                    y^2 h'' \left( y \right) & = & \left( 2y + y^2 + t^2 \right) h \left( y \right)
                \end{array}.
            \end{equation}
            Since $y^2$ and $2y + y^2 + t^2$ are both strictly positive, we get a contradiction. Therefore, the function $h$ cannot vanish if we have $t \neq 1/2$. Furthermore, the function $h$ vanishes identically if we have $t=1/2$, a fact which can be proved by direct computation using the special values
            \begin{equation}
            \label{eq:specialValueHalfModBesselFunction}
                \begin{array}{lllllll}
                    K_{1/2} \left( x \right) & = & \sqrt{\frac{\pi}{2x}} \, e^{-x} & \text{ and } & \frac{\partial}{\partial t}_{\vert t = 1/2} \, K_t \left( x \right) & = & \sqrt{\frac{\pi}{2x}} \, \mathbb{E}_1 \left( 2x \right) e^x
                \end{array}
            \end{equation}
            given in \cite[Prop. C.7, C.8]{dutour:cuspDirichlet}. Using these, we can further prove that $\pm i/2$ is a simple zero of $f_0$. Finally, the real zeros of $f_0$ are also simple, by \cite[Prop C.9]{dutour:cuspDirichlet} and \eqref{eq:orderDerivativeFunctiong}. This concludes the proof of this proposition.
        \end{proof}
        
        \begin{corollary}
        \label{cor:locZerosK0}
            The function
            \begin{equation}
                \begin{array}{lll}
                    s & \longmapsto & \left( 1 + 4 \pi \alpha a \right) K_{s-1/2} \left( 2 \pi \alpha a \right)^2 \\[0.3em]
                    
                    && \qquad \qquad \qquad \qquad + 4 \pi \alpha a K_{s-1/2} ' \left( 2 \pi \alpha a \right) K_{s-1/2} \left( 2 \pi \alpha a \right)
                \end{array}
            \end{equation}
            is entire, and all of its zeros are simple, given by a discrete set
            \begin{equation}
            \label{eq:zerosFK0}
                \begin{array}{lll}
                    \left \{ 1 \right \} \; \; \sqcup \; \; \left \{ \frac{1}{2} + i r_{0,j} \right \} & \subset & \left \{ 1 \right \} \; \; \sqcup \; \; \left \{ \frac{1}{2} + ir, \; r \in \mathbb{R}^{\ast} \right \}.
                \end{array}
            \end{equation}
            The corresponding eigenvalues of the spectral problem \eqref{eq:spectralProblem1} are given by
            \begin{equation}
                \begin{array}{lllllll}
                    \lambda_{0,0} & = & 0 & \text{ and } & \lambda_{0,j} & = & \frac{1}{4} + r_{0,j}^2.
                \end{array}
            \end{equation}
        \end{corollary}
        
        \begin{proof}
            This is a direct consequence of proposition \ref{prop:zerosFunctionF0}.
        \end{proof}

\section{Integral representation of the spectral zeta function}

    In this paper, we wish to define and study the zeta-regularized determinant of the pseudo-Laplacian with Alvarez--Wentworth boundary conditions by setting
    \begin{equation}
    \label{eq:firstDefDeterminant}
        \begin{array}{lll}
            \log \det \left( \Delta_{L,0} + \mu \right) & = & - \zeta'_{L,\mu} \left( 0 \right)
        \end{array}
    \end{equation}
    for any real number $\mu \geqslant 0$, where $\zeta_{L,\mu}$ is the spectral zeta function from \eqref{eq:spectralZetaFunction}. However, equality \eqref{eq:firstDefDeterminant} does not make sense \textit{a priori}, as the spectral zeta function is only defined and holomorphic on the half-plane $\Re s > 1$. In order to prove that $\zeta_{L,\mu}$ has a holomorphic continuation to some open subset of $\mathbb{C}$ containing $0$, and study some asymptotic properties of $\zeta'_{L,\mu} \left( 0 \right)$, we must find an integral representation of the spectral zeta function. Such is the object of the present section.

    \subsection{Contours of integration}
    \label{SubSec:contoursOfIntegration}
    
        In order to write $\zeta_{L,\mu}$ as an integral, we will need to use the \textit{argument principle}, for which the reader is referred to \cite[Sec. 3.4]{stein-shakarchi:complexAnalysis}. For this technique to work, we need to find a suitable contour of integration for the function $f_k$. There will be a difference according to whether or not we have $k=0$, due to the slight change in the localization of zeros. This part of the argument is similar to \cite[Sec. 3.1]{dutour:cuspDirichlet}, and to \cite[Sec. 6.1]{freixas-vonPippich:RROrbifold}.

        \subsubsection{Large order asymptotics}
        \label{subsubsec:largeOrderAsympt}
        
            The argument principle can technically only be applied for closed, and in particular bounded, contours. However, in order to recover the full spectral zeta function, we need to capture infitely many eigenvalues, requiring the use of unbounded contours. The justification of why that is possible relies in part on large order asymptotics of the modified Bessel functions of the second kind, which we recall below.
            
            \begin{proposition}
            \label{prop:largeOrderAsymptoticModBessel}
                For any $\nu, z \in \mathbb{C}$, with $\Re \nu > 0$ and $\Re z > 0$, we have
                \begin{equation}
                \label{eq:largeOrderAsymptoticModBessel}
                    \begin{array}{lll}
                        K_{\nu} \left( \nu z \right) & = & \sqrt{\frac{\pi}{2 \nu}} e^{-\nu \xi \left( z \right)} \left( 1+z^2 \right)^{-1/4} \left[ \; \sum\limits_{k=0}^{n-1} \frac{\left( -1 \right)^k}{\nu^k} U_k \left( p \left( z \right) \right) + \eta_n \left( \nu, z \right) \; \right]
                    \end{array}
                \end{equation}
                for any integer $n \geqslant 0$, where the functions $\xi$ and $p$ are defined by
                \begin{equation}
                \label{eq:defFunctionXi}
                    \begin{array}{lllllll}
                        \xi \left( z \right) & = & \sqrt{1+z^2} + \log \frac{z}{1 + \sqrt{1+z^2}} & \text{and} & p \left( z \right) & = & \left( 1+z^2 \right)^{-1/2}
                    \end{array}
                \end{equation}
                for any $z \in \mathbb{C}$ with $\Re z > 0$. In these equalities, the principal branch of the logarithm is considered, to define the square root as well. Furthermore $\nu^n \eta_n$ is uniformly bounded in $\left \vert \nu \right \vert \geqslant A$ for any fixed $A > 0$, and in $z$ on the sector $\left \vert \arg z \right \vert \leqslant \pi/2 - \delta$ for any fixed $\delta > 0$, and we have
                \begin{equation}
                    \begin{array}{llll}
                        \left \vert \eta_n \left( \nu, z \right) \right \vert & \leqslant & \frac{2}{\nu^n} \exp \left( \frac{2}{\nu} \mathcal{V}_{0,p \left( z \right)} \left( U_1 \right) \right) \mathcal{V}_{0,p \left( z \right)} \left( U_n \right) & \text{if } \nu > 0,
                    \end{array}
                \end{equation}
                where $\mathcal{V}_{0, p \left( z \right)}$ denotes the total variation of a function between the points $0$ and $p \left( z \right)$, and the polynomials $U_k$ are defined inductively by
                \begin{equation}
                    \begin{array}{lll}
                        U_0 \left( t \right) & = & 1, \\[0.5em]
                        U_{k+1} \left( t \right) & = & \frac{1}{2} t^2 \left( 1-t^2 \right) U_k' \left( t \right) + \frac{1}{8} \int_0^t \left( 1-5x^2 \right) U_k \left( x \right) \mathrm{d}x.
                    \end{array}
                \end{equation}
            \end{proposition}
            
            \begin{proof}
                The expansion for complex orders is given in \cite[Sec. 3]{olver:asymptoticExpBesselLargeOrder}, while the estimate on the remainder $\eta_2$ whenever $\nu$ is real is \cite[Eq. (7.15)]{olver}.
            \end{proof}
            
            \begin{remark}
                For more information on the total variation of functions, the reader is referred to \cite[App. C.1]{dutour:cuspDirichlet}, which itself is based on \cite[Sec. 1.11]{olver}. In particular, for any real numbers $\nu > 0$ and $x > 0$, we have
                \begin{equation}
                    \begin{array}{lll}
                        \left \vert \eta_2 \left( \nu, \frac{x}{\nu} \right) \right \vert & \leqslant & \frac{B}{x^2 + t^2}
                    \end{array}
                \end{equation}
                for some constant $B > 0$.
            \end{remark}
            
            Let us now deduce from proposition \ref{prop:largeOrderAsymptoticModBessel} the asymptotics for large order of the logarithmic derivative with respect to the order of the modified Bessel functions of the second kind.
            
            \begin{corollary}
            \label{cor:boundLogDerModBessel}
                For any fixed $z \in \mathbb{C}$ with $\Re z > 0$, the function
                \begin{equation}
                \label{eq:largeOrderBoundednessLogDerModBessel}
                    \begin{array}{lll}
                        \nu & \longmapsto & \frac{1}{K_{\nu} \left( z \right)} \frac{\partial}{\partial \nu} K_{\nu} \left( z \right)
                    \end{array}
                \end{equation}
                is uniformly bounded, for $\left \vert \nu \right \vert \geqslant A$ large enough with $\left \vert \arg \nu \right \vert < \pi/2 - \delta$, where $\delta > 0$ is some constant.
            \end{corollary}
            
            \begin{proof}
                This is a combination of proposition \ref{prop:largeOrderAsymptoticModBessel} and Cauchy's formula.
            \end{proof}
            
            In order to study asymptotics of the functions $f_k$ from \eqref{eq:functionf}, we also need large order asymptotics for the derivative with respect to the argument $z$ of the modified Bessel functions of the second kind. The notations of proposition \ref{prop:largeOrderAsymptoticModBessel} are kept here.
            
            \begin{proposition}
            \label{prop:largeOrderAsymptDerivativeModBessel}
                For any $\nu, z \in \mathbb{C}$, with $\Re \nu > 0$ and $\Re z > 0$, we have
                \begin{equation}
                \label{eq:largeOrderAsymptDerivativeModBessel}
                    \begin{array}{lll}
                        K_{\nu}' \left( \nu z \right) & = & - \sqrt{\frac{\pi}{2 \nu}} e^{-\nu \xi \left( z \right)} \frac{1}{z} \left( 1+z^2 \right)^{1/4} \left[ \; \sum\limits_{k=0}^{n-1} \frac{\left( -1 \right)^k}{\nu^k} V_k \left( p \left( z \right) \right) \right. \\[1em]
                        
                        && \qquad \left. + \frac{\left( -1 \right)^n}{\nu^n} \left( V_n - U_n \right) \left( p \left( z \right) \right) - \kappa_n \left( \nu, z \right) + \frac{z^2 p \left( z \right)^3}{2 \nu} \eta_n \left( \nu, z \right) \; \right]
                    \end{array}
                \end{equation}
                for any integer $n \geqslant 0$. The additional remainder $\kappa_n$ satisfies the same properties and estimates as $\eta_n$ from proposition \ref{prop:largeOrderAsymptoticModBessel}. Furthermore, the polynomials $V_k$ are defined inductively by
                \begin{equation}
                    \begin{array}{lll}
                        V_0 \left( t \right) & = & 1, \\[0.5em]
                        
                        V_k \left( t \right) & = & U_k - \frac{1}{2} t \left( 1 - t^2 \right) U_{k-1} -t^2 \left( 1 - t^2 \right) U_{k-1} ' \left( t \right).
                    \end{array}
                \end{equation}
            \end{proposition}
            
            \begin{proof}
                This asymptotic expansion can be obtained from proposition \ref{prop:largeOrderAsymptoticModBessel} by applying Cauchy's formula. It is also given in \cite[Ex. 10.7.2]{olver}.
            \end{proof}
            
            \begin{remark}
                Note that \eqref{eq:largeOrderAsymptDerivativeModBessel} is written that way so as to have a more precise control on the remainder than \cite[Eq. (2.21)]{olver:asymptoticExpBesselLargeOrder}.
            \end{remark}
            
            \begin{proposition}
            \label{prop:largeOrderAsymptLogDerivativeModBessel}
                For any $\nu \in \mathbb{C}$ with $\Re \nu > 0$, and any $x > 0$, we have
                \begin{equation}
                \label{eq:largeOrderAsymptLogDerivativeModBessel}
                    \begin{array}{lll}
                        x \frac{K_{\nu}' \left( x \right)}{K_{\nu} \left( x \right)} & = & - \nu \sqrt{1 + \frac{x^2}{\nu^2}} \left( 1 - \frac{1}{\nu} \left( U_1 \left( p \left( \frac{x}{\nu} \right) \right) + V_1 \left( p \left( \frac{x}{\nu} \right) \right) \right) + \eta_{2,1} \left( \nu, x \right) \right)
                    \end{array}
                \end{equation}
                where, for each $x$, the term $\nu^2 \eta_{2,1} \left( \nu, x \right)$ is uniformly bounded in $\nu$ for $\left \vert \nu \right \vert \geqslant A$ large enough, with $A > 0$ being a constant. Furthermore, we have
                \begin{equation}
                    \begin{array}{llll}
                        \left \vert \eta_{2,1} \left( \nu, x \right) \right \vert & \leqslant & \frac{C}{x^2 + \nu^2} & \text{if } \nu > 0
                    \end{array}
                \end{equation}
                for some constant $C > 0$, and either $\nu \geqslant A$ large enough or $x \geqslant B$ large enough.
            \end{proposition}
            
            \begin{proof}
                This asymptotic expansion is found using propositions \ref{prop:largeOrderAsymptoticModBessel} and \ref{prop:largeOrderAsymptDerivativeModBessel}, and the estimates on the remainder are proved using a power series expansion for the denominator appearing after the factors common to $K_{\nu} \left( x \right)$ and $K_{\nu}' \left( x \right)$ are removed.
            \end{proof}
            
            \begin{definition}
                For any complex number $\nu$ and any real numbers $x,y>0$, we set
                \begin{equation}
                \label{eq:functiongtxy}
                    \begin{array}{lll}
                        g \left( \nu, x, y \right) & = & 1 + 4 \pi \alpha a + x \frac{K_{\nu}' \left( x \right)}{K_{\nu} \left( x \right)} + y \frac{K_{\nu}' \left( y \right)}{K_{\nu} \left( y \right)},
                    \end{array}
                \end{equation}
                with the understanding that this function of $\nu$ has poles. We further set
                \begin{equation}
                \label{eq:functionftxy}
                    \begin{array}{lll}
                        f \left( \nu, x, y \right) & = & K_{\nu} \left( x \right) K_{\nu} \left( y \right) g \left( t, x, y \right).
                    \end{array}
                \end{equation}
            \end{definition}
            
            \begin{remark}
                It is more convenient to work with \eqref{eq:functiongtxy} than $g_k$, as it lightens the notation. Nevertheless, for any integer $k \neq 0$, we have
                \begin{equation}
                    \begin{array}{lll}
                        g_k \left( t \right) & = & g \left( it, 2 \pi \left \vert k + \alpha \right \vert a, 2 \pi \left \vert k - \alpha \right \vert a \right), \\[0.5em]
                        
                        f_k \left( t \right) & = & f \left( it, 2 \pi \left \vert k + \alpha \right \vert a, 2 \pi \left \vert k - \alpha \right \vert a \right)
                    \end{array}
                \end{equation}
            \end{remark}
            
            \begin{proposition}
            \label{prop:largeOrderAsymptFunctionGtxy}
                For any $x,y > 0$, and $t \in \mathbb{C}$ with $\arg t < \pi/2$, we have
                \begin{equation}
                \label{eq:asymptLargeOrderFunctionGtxy}
                    \begin{array}{lll}
                        g \left( t, x, y \right) & = & - \left( \sqrt{x^2 + t^2} + \sqrt{y^2 + t^2} \right) \left[ 1 - \frac{4 \pi \alpha a}{\sqrt{x^2 + t^2} + \sqrt{y^2 + t^2}} \right. \\[1em]
                        
                        && \; \; \; \left. - \frac{1}{2} \cdot \frac{1}{\sqrt{x^2 + t^2} + \sqrt{y^2 + t^2}} \cdot \left( \frac{t^2}{x^2 + t^2} + \frac{t^2}{y^2 + t^2} \right) + \eta_{2,2} \left( t,x,y \right) \right],
                    \end{array}
                \end{equation}
                where $t^2 \eta_{2,2} \left( t, x, y \right)$ is uniformly bounded in $t$ for $\left \vert t \right \vert \geqslant A$ large enough. Furthermore, we have
                \begin{equation}
                \label{eq:estimateRemainderEta22}
                    \begin{array}{llll}
                        \left \vert \eta_{2,2} \left( t, x, y \right) \right \vert & \leqslant & C \left[ \frac{1}{x^2+t^2} + \frac{1}{y^2+t^2} \right] & \text{if } t > 0
                    \end{array}
                \end{equation}
                for some constant $C > 0$, with either $t \geqslant A$ or $x, y \geqslant B > 0$ large enough. In particular, the function $t \longmapsto 1/g \left( t, x, y \right)$ is uniformly bounded in $t$ for $\left \vert t \right \vert \geqslant A$ large enough with $\left \vert \arg t \right \vert < \pi/2 - \delta$ for some constant $\delta > 0$, and we also have, for such complex numbers $t$,
                \begin{equation}
                \label{eq:boundOrderFunctionGtxy}
                    \begin{array}{lll}
                        \left \vert g \left( t, x, y \right) \right \vert & \leqslant & D \left \vert t \right \vert
                    \end{array}
                \end{equation}
                for some constant $D > 0$. As a consequence, we have
                \begin{equation}
                \label{eq:boundDerOrderFunctionGtxy}
                    \begin{array}{lll}
                        \left \vert \frac{\partial}{\partial t} g \left( t, x, y \right) \right \vert & \leqslant & D \left( \frac{1}{\varepsilon} \left \vert t \right \vert + 1 \right)
                    \end{array}
                \end{equation}
                again for $\left \vert t \right \vert$ large enough, for any constant $\varepsilon > 0$ small enough.
            \end{proposition}
            
            \begin{proof}
                The first part of this result is a consequence of proposition \ref{prop:largeOrderAsymptLogDerivativeModBessel}, while \eqref{eq:boundDerOrderFunctionGtxy} is obtained by applying Cauchy's formula, using circles of radius $\varepsilon$ small enough so that \eqref{eq:boundOrderFunctionGtxy} still holds on them.
            \end{proof}
            
            \begin{corollary}
            \label{cor:boundLogDerFunctionFtxy}
                For any fixed real numbers $x,y > 0$ and any complex number $t$ with $\left \vert \arg t \right \vert > \delta$ for some constant $\delta > 0$ and $\left \vert t \right \vert \geqslant A$ large enough, we have
                \begin{equation}
                \label{eq:boundDerOrderFunctionFtxy}
                    \begin{array}{lll}
                        \left \vert \frac{1}{f \left( t, x, y \right)} \frac{\partial}{\partial t} f \left( t, x, y \right) \right \vert & \leqslant & E \left( \left \vert t \right \vert + 1 \right)
                    \end{array}
                \end{equation}
                where $E > 0$ is a constant.
            \end{corollary}
            
            \begin{proof}
                This is a consequence of corollary \ref{cor:boundLogDerModBessel} and proposition \ref{prop:largeOrderAsymptFunctionGtxy}.
            \end{proof}

        \subsubsection{Contour of integration for $k \neq 0$}
        
            Let us first describe the contour of integration we will use to detect the zeros of $f_k$ for non-zero integers $k$.
            
            \begin{definition}
                For any $\vartheta \in \left] 0, \pi/2 \right[$, the unbounded contour $\gamma_{\vartheta,\pm}$ is defined by
                \begin{equation}
                \label{eq:contourGammaPM}
                    \begin{array}{lll}
                        \gamma_{\vartheta,\pm} & = & \left \{ r e^{i \vartheta}, \; r \geqslant 0 \right \} \; \cup \; \left \{ r e^{- i \vartheta}, \; r \geqslant 0 \right \}
                    \end{array} .
                \end{equation}
                Its bounded counterpart $\gamma_{\vartheta, R, \pm}$ is defined for any $R > 0$ by
                \begin{equation}
                \label{eq:boundedContourGammaPM}
                    \begin{array}{lll}
                        \gamma_{\vartheta,R,\pm} & = & \left \{ r e^{i \vartheta}, \; R \geqslant r \geqslant 0 \right \} \; \cup \; \left \{ r e^{- i \vartheta}, \; R \geqslant r \geqslant 0 \right \} \; \cup \; \mathcal{C}_{\vartheta, R}
                    \end{array}
                \end{equation}
                where $\mathcal{C}_{\vartheta, R} = \left \{ R e^{- i \varphi}, \; \varphi \in \left[ - \vartheta, \vartheta \right] \right \}$ is the circular arc at radius $R$ and angles $\pm \vartheta$.
            \end{definition}
            
            \begin{remark}
                This contour, as well as its bounded counterpart, which is capped by a circular arc, are represented below. In order to use the argument principle with the infinite contour, it will be necessary to prove that the contribution due to the circular arc vanishes as the radius goes to infinity.
            
                \begin{figure}[H]
    
                    \begin{center}
        
                        \begin{subfigure}{0.45\textwidth}
                
                            \centering
                
                            \begin{tikzpicture}[scale=0.9]
        
                                \draw[thick] (-2,0) -- (0,0);
                                \draw[thick] (0,-2) -- (0,0);
                                \draw[thick,->] (0,0) coordinate (2) node {} -- (3.3,0) coordinate (1) node {};
                                \draw[thick, ->] (0,0) -- (0,2);
            
                                \draw[thick, red] (0,0) -- (3,3/2) coordinate (3) node {};
                                \draw[thick, red] (0,0) -- (3,-3/2) coordinate (4) node {};
            
                                \draw pic["$\color{orange} \vartheta$", draw=orange, ->, angle eccentricity=1.2, angle radius=1.3cm] {angle=1--2--3};
            
                                \node (P) at (2.5,0.7) {$\color{red} \gamma_{\vartheta, \pm}$};
            
                                \draw[thick, red, ->] (2.5,1.25) -- (2,1);
                                \draw[thick, red, ->] (2,-1) -- (2.5,-1.25);
                    
                            \end{tikzpicture}
                    
                            \caption{Unbounded Integration contour $\gamma_{\vartheta,\pm}$}
                            \label{fig:integrationContourPM}
        
                            \end{subfigure}
                            \hspace{0.9cm}
                            \begin{subfigure}{0.44\textwidth}
            
                            \centering
                            
                            \begin{tikzpicture}[scale=0.9]
        
                                \draw[thick] (-2,0) -- (0,0);
                                \draw[thick] (0,-1.8) -- (0,0);
                                \draw[thick,->] (0,0) coordinate (2) node {} -- (3.6,0) coordinate (1) node {};
                                \draw[thick, ->] (0,0) -- (0,1.8);
            
                                \draw[thick, blue] (0,0) -- (3,3/2) coordinate (3) node {};
                                \draw[thick, blue] (0,0) -- (3,-3/2) coordinate (4) node {};
            
                                \draw pic["$\color{orange} \vartheta$", draw=orange, ->, angle eccentricity=1.2, angle radius=1.3cm] {angle=1--2--3};
                                
                                \draw pic["", draw=blue, thick, -, angle eccentricity=1, angle radius=3cm] {angle=4--2--3};
            
                                \node (P) at (2.4,1.7) {$\color{blue} \gamma_{\vartheta, R, \pm}$};
                                
                                \node (Q) at (3.6, 0.4) {$\color{blue} R$};
            
                                \draw[thick, blue, ->] (2.5,1.25) -- (2,1);
                                \draw[thick, blue, ->] (2,-1) -- (2.5,-1.25);
                                \draw[thick, blue, ->] (3.23,-0.8) -- (3.24,-0.75);
                                \draw[thick, blue, ->] (3.24,0.75) -- (3.23,0.8);
                    
                            \end{tikzpicture}
                            
                            \caption{Bounded integration contour $\gamma_{\vartheta, R, \pm}$}
                            \label{fig:boundedIntegrationContourPM}
        
                        \end{subfigure}
            
                    \end{center}
        
                    \caption{Integration contours for $k \neq 0$}
                    \label{fig:integrationContoursKNeq0}
        
                \end{figure}
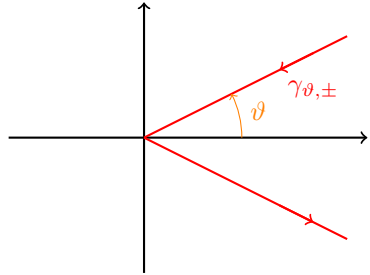
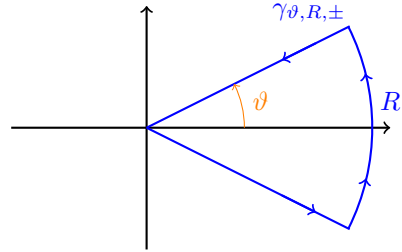
        
            \end{remark}

        \subsubsection{Contour of integration for $k = 0$}
        
            The integration contour for $k = 0$ is slightly more complicated, as we need to avoid $\pm i/2$ for an angle $\vartheta$ close to $\pi/2$. To simplify, we will not mention the bounded counterpart of the contour of integration.
            
            \begin{definition}
                For any $\vartheta \in \left] \pi/6, \pi/2 \right[$ and any $1 > \varepsilon > 0$ small enough, we define the contour of integration $\gamma_{\vartheta,\varepsilon}$ by
                \begin{equation}
                \label{eq:contourGamma0}
                    \begin{array}{lll}
                        \gamma_{\vartheta,0,\varepsilon} & = & \scriptstyle \left \{ r e^{i \vartheta}, \; r \geqslant \sqrt{\frac{1}{4} + \varepsilon \mu} \right \} \; \cup \; \left \{ \frac{1}{2} e^{i \vartheta} + \left( \sqrt{\frac{1}{4} + \varepsilon \mu} - \frac{1}{2} \right) e^{i \varphi}, \; \varphi \in \left] -\pi, 0 \right[ \right \} \\[1em]
                    
                        && \scriptstyle \quad \cup \; \left \{ r e^{i \vartheta}, \; r \leqslant 1 - \sqrt{\frac{1}{4} + \varepsilon \mu} \right \} \; \cup \; \left \{ r e^{-i \vartheta}, \; r \leqslant 1 - \sqrt{\frac{1}{4} + \varepsilon \mu} \right \} \\[1em]
                        
                        && \scriptstyle \quad \quad \cup \; \left \{ \frac{1}{2} e^{- i \vartheta} + \left( \sqrt{\frac{1}{4} + \varepsilon \mu} - \frac{1}{2} \right) e^{i \varphi}, \; \varphi \in \left] 0, \pi \right[ \right \} \; \cup \; \left \{ r e^{-i \vartheta}, \; r \geqslant \sqrt{\frac{1}{4} + \varepsilon \mu} \right \}.
                    \end{array}
                \end{equation}
            \end{definition}
            
            \begin{remark}
                In order for this contour to make sense, we need the different parts to not intersect each other. This is possible for $\varepsilon > 0$ small enough, if we impose a lower bound on the angle, which is why we restricted ourselves to $\vartheta > \pi/6$.
            \end{remark}
            
            \begin{remark}
                The contour $\gamma_{\vartheta, 0, \varepsilon}$ is represented below.
                
                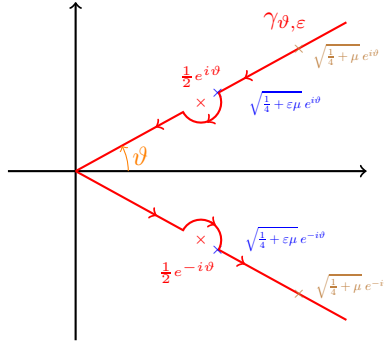
\begin{figure}[H]
    
                    \begin{center}
        
                        \centering
                
                            \begin{tikzpicture}[scale=0.9]
        
                                \draw[thick] (-1,0) -- (0,0);
                                \draw[thick] (0,-2.5) -- (0,0);
                                \draw[thick,->] (0,0) coordinate (2) node {} -- (4.3,0) coordinate (1) node {};
                                \draw[thick, ->] (0,0) -- (0,2.5);
            
                                \draw[thick, red] (0,0) -- (1.6,0.55*1.6) coordinate (3) node {};
                                \draw[thick, red] (2.1,0.55*2.1) -- (4,0.55*4) coordinate (5) node {};
                                \draw[thick, red] (0,0) -- (1.6,-0.55*1.6) coordinate (4) node {};
                                \draw[thick, red] (2.1,-0.55*2.1) -- (4,-0.55*4) coordinate (8) node {};
                                
                                \node (6) at (1.85,0.55*1.85) {$\color{red} \scriptscriptstyle \times$};
                                \node (Q) at (1.85,1.4) {$\color{red} \scriptscriptstyle \frac{1}{2} e^{i \vartheta}$};
                                
                                \node (7) at (1.85,-0.55*1.85) {$\color{red} \scriptscriptstyle \times$};
                                \node (R) at (1.65,-1.5) {$\color{red} \scriptscriptstyle \frac{1}{2} e^{-i \vartheta}$};
                                
                                \node (A) at (2.1,0.55*2.1) {$\color{blue} \scriptscriptstyle \times$};
                                \node (B) at (2.1,-0.55*2.1) {$\color{blue} \scriptscriptstyle \times$};
                                
                                \node (9) at (3.1,1) {\color{blue} \scalebox{0.5}{$\sqrt{\frac{1}{4} + \varepsilon \mu} \; e^{i \vartheta}$}};
                                \node (10) at (3.1,-1) {\color{blue} \scalebox{0.5}{$\sqrt{\frac{1}{4} + \varepsilon \mu} \; e^{-i \vartheta}$}};
                                
                                \node (C) at (3.3,0.55*3.3) {$\color{brown} \scriptscriptstyle \times$};
                                \node (D) at (3.3,-0.55*3.3) {$\color{brown} \scriptscriptstyle \times$};
                                
                                \node (11) at (4,1.7) {\color{brown} \scalebox{0.5}{$\sqrt{\frac{1}{4} + \mu} \; e^{i \vartheta}$}};
                                \node (12) at (4.1,-1.7) {\color{brown} \scalebox{0.5}{$\sqrt{\frac{1}{4} + \mu} \; e^{-i \vartheta}$}};
                                
                                \draw[thick, red] pic["", draw=red, -, angle eccentricity=1, angle radius=0.27cm] {angle=2--6--5};
                                \draw[thick, red, ->] (1.95,0.725) -- (1.94,0.718);
                                
                                \draw[thick, red] pic["", draw=red, -, angle eccentricity=1, angle radius=0.27cm] {angle=8--7--2};
                                \draw[thick, red, ->] (2,-0.755) -- (2.07,-0.795);
            
                                \draw pic["$\color{orange} \; \vartheta$", draw=orange, ->, angle eccentricity=1.2, angle radius=0.7cm] {angle=1--2--3};
            
                                \node (P) at (3.1,2.2) {$\color{red} \gamma_{\vartheta,\varepsilon}$};
                                
                                \draw[thick, red, ->] (2.5,0.55*2.5) -- (2.49,0.55*2.49);
                                \draw[thick, red, ->] (1.2,0.55*1.2) -- (1.18,0.55*1.18);
                                \draw[thick, red, ->] (1.18,-0.55*1.18) -- (1.2,-0.55*1.2);
                                \draw[thick, red, ->] (2.49,-0.55*2.49) -- (2.5,-0.55*2.5);
                    
                            \end{tikzpicture}
                    
                            \caption{Integration contour $\gamma_{\vartheta,0,\varepsilon}$}
                            \label{fig:integrationContourK0}
        
                    \end{center}
                \end{figure}
            \end{remark}

        \subsubsection{Integral representation}
        \label{subsubsec:integralRepresentation}
        
            We now wish to apply the argument principle with the contours of integration $\gamma_{\vartheta,\pm}$, defined in \eqref{eq:contourGammaPM}, and $\gamma_{\vartheta,0,\varepsilon}$, defined in \eqref{eq:contourGamma0}. However, as indicated in the introduction of paragraph \ref{subsubsec:largeOrderAsympt}, the argument principle can \textit{a priori} only be applied for closed, bounded contours. With the following lemma, we will justify why the unbounded contour $\gamma_{\vartheta,\pm}$ can be used. To lighten the presentation, we omit the corresponding lemma for $\gamma_{\vartheta,0,\varepsilon}$.
            
            \begin{lemma}
            \label{lem:circularContributionArgPrinciple}
                Consider fixed real numbers $x, y > 0$, and a sequence $\left( R_m \right)_m$ of strictly positive real numbers with $R_m \rightarrow + \infty$, and which are not zeros of $f \left( \cdot, x, y \right)$. For any complex number $s$ with $\Re s > 1$, we have
                \begin{equation}
                    \begin{array}{lll}
                        \int_{\mathcal{C}_{\vartheta, R_m}} \left( \frac{1}{4} + t^2 + \mu \right)^{-s} \frac{\partial}{\partial t} \log f \left( t, x, y \right) \mathrm{d}t & \underset{m \rightarrow + \infty}{\longrightarrow} & 0.
                    \end{array}
                \end{equation}
            \end{lemma}
            
            \begin{proof}
                Let us first note that the estimate provided by corollary \ref{cor:boundLogDerFunctionFtxy} cannot be applied on any $\mathcal{C}_{\vartheta, R_m}$, as this circular arc intersects the axis $\Im t = 0$. Furthermore, the hypothesis according to which $R_m$ is not a zero of $f \left( \cdot, x, y \right)$ is necessary to ensure integrability. Let us set
                \begin{equation}
                    \begin{array}{lll}
                        \gamma_{\varphi, R} & = & \left \{ r e^{i \varphi}, \; r \geqslant R \right \}
                    \end{array}
                \end{equation}
                with $\varphi$ being either $\vartheta$ or $- \vartheta$ and any $R > 0$. Using corollary \ref{cor:boundLogDerFunctionFtxy}, we have
                \begin{equation}
                    \begin{array}{lll}
                        \int_{\gamma_{\varphi, R_m}} \left( \frac{1}{4} + t^2 + \mu \right)^{-s} \frac{\partial}{\partial t} \log f \left( t, x, y \right) \mathrm{d}t & \underset{m \rightarrow + \infty}{\longrightarrow} & 0
                    \end{array}
                \end{equation}
                for $\varphi = \pm \vartheta$ and $\Re s > 1$. The absolute convergence of the spectral zeta function on the half-plane $\Re s > 1$ and the argument principle now prove that the sequence
                \begin{equation}
                    \begin{array}{c}
                        \left( \int_{\mathcal{C}_{\vartheta, R_m}} \left( \frac{1}{4} + t^2 + \mu \right)^{-s} \frac{\partial}{\partial t} \log f \left( t, x, y \right) \mathrm{d}t \right)_m
                    \end{array}
                \end{equation}
                is Cauchy, and thus converges, for any $s$ with $\Re s > 1$. Extracting part of the complex power with an upper bound in terms of $R$ finally proves that the limit must be $0$, thus concluding the proof.
            \end{proof}
            
            \begin{proposition-definition}
            \label{prop:spectralZetaIntegralRepresentation}
                On the half-plane $\Re s > 1$, we have
                \begin{equation}
                    \begin{array}{lll}
                        \zeta_{L,\mu} \left( s \right) & = & \left \{ \begin{array}{ll} \zeta_{L,\mu}^{\pm} \left( s \right) & \text{if } \alpha = 0 \\[0.5em] \zeta_{L,\mu}^{\pm} \left( s \right) + \zeta_{L,\mu}^{0} \left( s \right) & \text{if } \alpha \neq 0 \end{array} \right.
                    \end{array}
                \end{equation}
                where the partial spectral zeta function $\zeta_{L,\mu}^{\pm}$ is defined for any $\vartheta \in \left[ 0 , \pi/2 \right[$ by
                \begin{equation}
                \label{eq:partialZetaPM}
                    \begin{array}{lll}
                        \zeta_{L,\mu}^{\pm} \left( s \right) & = & \frac{1}{2 i \pi} \sum\limits_{k \neq 0} \int_{i \gamma_{\vartheta,\pm}} \left( \frac{1}{4} - t^2 + \mu \right)^{-s} \frac{\partial}{\partial t} \log f_k \left( i t \right) \mathrm{d}t,
                    \end{array}
                \end{equation}
                and, if we have $\alpha \neq 0$, the partial spectral zeta function $\zeta_{L,\mu}^{0}$ is defined by
                \begin{equation}
                \label{eq:partialZeta0}
                    \begin{array}{lll}
                        \zeta_{L,\mu}^{0} \left( s \right) & \hspace{-4pt} = \hspace{-4pt} & \left \{ \begin{array}{ll} \frac{1}{2 i \pi} \int_{i \gamma_{\vartheta, \varepsilon}} \left( \frac{1}{4} - t^2 \right)^{-s} \frac{\partial}{\partial t} \log f_0 \left( i t \right) \mathrm{d}t & \text{if } \mu = 0 \\[0.5em] \frac{1}{2 i \pi} \int_{i \gamma_{\vartheta, \varepsilon}} \left( \frac{1}{4} - t^2 + \mu \right)^{-s} \frac{\partial}{\partial t} \log f_0 \left( i t \right) \mathrm{d}t + \mu^{-s} & \text{if } \mu > 0 \end{array} \right.
                    \end{array}
                \end{equation}
                for any angle $\vartheta \in \left] \pi/6, \pi/2 \right[$.
            \end{proposition-definition}
            
            \begin{proof}
                This result is a consequence of the argument principle, which may be applied with the unbounded contours $\gamma_{\vartheta, \pm}$ and $\gamma_{\vartheta, 0, \varepsilon}$ using lemma \ref{lem:circularContributionArgPrinciple} and its analogue for $\gamma_{\vartheta, 0, \varepsilon}$, which is not explicitely stated to shorten the presentation. Note that a change of variable is performed, to rotate the contours by $\pi/2$.
            \end{proof}

    \subsection{Factorization by the Dirichlet contribution}
    \label{subsec:factorizationDirichletContribution}
    
        In \eqref{eq:partialZetaPM} and \eqref{eq:partialZeta0}, the definitions of the partial zeta functions $\zeta_{L,\mu}^{\pm}$ and $\zeta_{L,\mu}^{0}$ are based on the argument principle for the functions $f_k$. However, in \eqref{eq:functionf}, we saw that the functions $f_k$ can be naturally factored. The contributions thus extracted correspond to the Dirichlet boundary conditions, which have been studied in \cite{dutour:cuspDirichlet}.

        \subsubsection{The zeta function $\zeta_{L,\mu}^{\pm}$}
        
            Let us first extract a Dirichlet component from $\zeta_{L,\mu}^{\pm}$, which is the the partial zeta function corresponding to the sum over $k \neq 0$ in the Fourier decomposition \eqref{eq:fourierDecompositions}.
        
            \begin{proposition-definition}
                For any $s \in \mathbb{C}$ with $\Re s > 1$, we have
                \begin{equation}
                    \begin{array}{lll}
                        \zeta_{L,\mu}^{\pm} \left( s \right) & = & \zeta_{L,\mu}^{\pm, \, D} \left( s \right) + \zeta_{L,\mu}^{\pm, \, AW} \left( s \right)
                    \end{array}
                \end{equation}
                for any real number $\mu \geqslant 0$, where we have set
                \begin{equation}
                \label{eq:zetaDandAWKNeq0}
                    \begin{array}{lll}
                        \zeta_{L,\mu}^{\pm, \, D} \left( s \right) & = & \frac{1}{i \pi} \sum\limits_{k \neq 0} \int_{i \gamma_{\vartheta,\pm}} \left( \frac{1}{4} - t^2 + \mu \right)^{-s} \frac{\partial}{\partial t} \log K_t \left( 2 \pi \left \vert k + \alpha \right \vert a \right) \mathrm{d}t, \\[1em]
                        
                        \zeta_{L,\mu}^{\pm, \, AW} \left( s \right) & = & \frac{1}{2 i \pi} \sum\limits_{k \neq 0} \int_{i \gamma_{\vartheta,\pm}} \left( \frac{1}{4} - t^2 + \mu \right)^{-s} \frac{\partial}{\partial t} \log g_k \left( i t \right) \mathrm{d}t,
                    \end{array}
                \end{equation}
                the function $g_k$ having been defined in \eqref{eq:functiong}.
            \end{proposition-definition}
            
            \begin{proof}
                This result is a direct consequence of the factorization \eqref{eq:functionf}.
            \end{proof}
            
            \begin{remark}
                Up to a factor $2$, the function $\zeta_{L,\mu}^{\pm, \, D}$ has been studied in \cite{dutour:cuspDirichlet}, hence the exponent ``D'' standing for Dirichlet. It should be noted that the results in \cite{dutour:cuspDirichlet} include the contribution from $k=0$, which we keep separate here. The other zeta function in \eqref{eq:zetaDandAWKNeq0} bears an exponent ``AW'' to emphasize that it is specific to the Alvarez--Wentworth boundary conditions.
            \end{remark}
            
            \begin{theorem}
            \label{thm:asymptDerAt0DirichletContributionPM}
                The function $\zeta_{L,\mu}^{\pm, \, D}$ has a holomorphic continuation to a domain of the complex plane containing $0$, and its derivative there satisfies
                \begin{equation}
                    \begin{array}{lll}
                        \zeta_{L,\mu}^{\pm, \, D} \, ' \left( 0 \right) & \hspace{-3pt} = \hspace{-3pt} & \frac{1}{2 \pi a} \mu \log \mu - \frac{1}{2 \pi a} \mu - 3 \sqrt{\mu} \log \mu + 2 \left[ - 3 \log 2 + \alpha \log \frac{1 + \alpha}{1 - \alpha} \right. \\[1em]
                            
                        && \; + \frac{1}{4a} - 2 \int_0^{+ \infty} \frac{1}{e^{2 \pi t} - 1} \left( \arctan \left( \frac{t}{1 + \alpha} \right) + \arctan \left( \frac{t}{1 - \alpha} \right) \right) \mathrm{d}t \\[1em]
                            
                        && \; \; \left. + 1 + \frac{3}{2} \log \left( 4 \pi^2 \left( 1 - \alpha^2 \right) a^2 \right) \right] \sqrt{\mu} - \frac{1}{2} \log \mu + o \left( 1 \right)
                    \end{array}
                \end{equation}
                as $\mu$ goes to infinity. The same derivative further satisfies
                \begin{equation}
                    \begin{array}{lll}
                        \zeta_{L,0}^{\pm, \, D} \, ' \left( 0 \right) & = & - \frac{2 \pi}{3} a - 4 \pi \alpha^2 a - \log a + \log \frac{\sin \left( \pi \alpha \right)}{\pi \alpha} + o \left( 1 \right).
                    \end{array}
                \end{equation}
                with $\mu = 0$, as $a$ goes to infinity.
            \end{theorem}
                
            \begin{proof}
                This is a consequence of the results proved in \cite{dutour:cuspDirichlet}, though they are not explicitely stated there, since the contribution from $k=0$ is included there.
            \end{proof}
            
            \begin{remark}
                The study of the function $\zeta_{L,\mu}^{\pm, \, AW}$ will start in subsection \ref{SubSec:thetaPiOver2}.
            \end{remark}

        \subsubsection{The function $\zeta_{L,\mu}^{0}$}
        
            Let us perform the same type of factorization on $\zeta_{L,\mu}^{0}$, which is the partial zeta function corresponding to the constant term $k=0$ in the Fourier decomposition \eqref{eq:fourierDecompositions}. In this paragraph, we assume that we have $\alpha \neq 0$.
            
            \begin{proposition-definition}
            \label{propDef:decompZetaFunction0DAWLogMu}
                For any $s \in \mathbb{C}$ with $\Re s > 1$, we have
                \begin{equation}
                \label{eq:decompZeta0}
                    \begin{array}{lll}
                        \zeta_{L,\mu}^{0} \left( s \right) & = & \zeta_{L,\mu}^{0, \, D} \left( s \right) + \zeta_{L,\mu}^{0, \, AW} \left( s \right) + \mu^{-s}
                    \end{array}
                \end{equation}
                for any real number $\mu > 0$, where we have set
                \begin{equation}
                \label{eq:zetaDandAWK0}
                    \begin{array}{lll}
                        \zeta_{L,\mu}^{0, \, D} \left( s \right) & = & \frac{1}{i \pi} \int_{i \gamma_{\vartheta, 0}} \left( \frac{1}{4} - t^2 + \mu \right)^{-s} \frac{\partial}{\partial t} \log K_t \left( 2 \pi \alpha a \right) \mathrm{d}t, \\[1em]
                        
                        \zeta_{L,\mu}^{0, \, AW} \left( s \right) & = & \frac{1}{2 i \pi} \int_{i \gamma_{\vartheta, 0}} \left( \frac{1}{4} - t^2 + \mu \right)^{-s} \frac{\partial}{\partial t} \log g_0 \left( i t \right) \mathrm{d}t,
                    \end{array}
                \end{equation}
                the function $g_0$ having been defined in \eqref{eq:functiong}.
            \end{proposition-definition}
            
            \begin{proof}
                This result is a direct consequence of the factorization \eqref{eq:functionf}.
            \end{proof}
            
            \begin{remark}
                The case $\mu = 0$ is obtained by removing the term $\mu^{-s}$ in \eqref{eq:decompZeta0}.
            \end{remark}
            
            \begin{theorem}
            \label{thm:asymptDerAt0DirichletContribution0}
                The function $\zeta_{L,\mu}^{0, \, D}$ has a holomorphic continuation to a domain of the complex plane containing $0$, and its derivative there satisfies
                \begin{equation}
                    \begin{array}{lll}
                        \zeta_{L,\mu}^{0, \, D} \, ' \left( 0 \right) & \hspace{-1pt} = \hspace{-1pt} & - \sqrt{\mu} \log \mu + 2 \left( \log \left( \pi \alpha a \right) - 1 \right) \sqrt{\mu} + \frac{1}{2} \log \mu + o \left( 1 \right)
                    \end{array}
                \end{equation}
                as $\mu$ goes to infinity. The same derivative further satisfies
                \begin{equation}
                    \begin{array}{lll}
                        \zeta_{L,0}^{0, \, D} \, ' \left( 0 \right) & = & 2 \pi \alpha a + \frac{1}{2} \log a + \frac{1}{2} \log \left( 2 \pi \alpha \right) + o \left( 1 \right)
                    \end{array}
                \end{equation}
                with $\mu = 0$, as $a$ goes to infinity.
            \end{theorem}
            
            \begin{proof}
                This is a consequence of the results proved in \cite{dutour:cuspDirichlet}, though they are not explicitely stated there, as this is only the result pertaining to $k=0$.
            \end{proof}
            
            \begin{remark}
                The study of the function $\zeta_{L,\mu}^{0, \, AW}$ will start in subsection \ref{SubSec:thetaPiOver2}.
            \end{remark}

    \subsection{\texorpdfstring{Letting $\vartheta$ go to $\pi/2$}{Letting theta go to pi/2}}
    \label{SubSec:thetaPiOver2}
    
        So far, the angle in the integration contours has to satisfy $\vartheta < \pi/2$. In order to let the angle converge to $\pi/2$ and retain a meaningful result, some care must be taken. The arguments used here follow \cite[Sec. 3.2]{dutour:cuspDirichlet}, which themselves are based on \cite[Sec. 6.1]{freixas-vonPippich:RROrbifold}. Two ingredients are needed in this section: a regularization of the functions being integrated in \eqref{eq:zetaDandAWKNeq0} and \eqref{eq:zetaDandAWK0} for the AW contribution, and large argument asymptotics. The first one will help take the limit of each integral individually, while the second one will justify that the limit can be interchanged with the sum over $k$.
        
        \subsubsection{Large argument asymptotics}
        
            Let us begin by presenting the large argument asymptotics of the modified Bessel functions of the second kind, and their corollary for the functions $g_k$.
        
            \begin{proposition}
            \label{prop:largeArgumentAsymptoticsBessel}
                For any complex number $z$ with $\Re z > 0$ and any complex number $\nu$, we have
                \begin{equation}
                \label{eq:asymptoticsArgumentBessel}
                    \begin{array}{lll}
                        K_{\nu} \left( z \right) & = & \sqrt{\frac{\pi}{2z}} e^{-z} \left[ 1 + \frac{1}{8} \left( 4\nu^2 - 1 \right) \cdot \frac{1}{z} \right. \\[0.5em]
                    
                        && \qquad \qquad \qquad + \left. \frac{1}{128} \left( 4 \nu^2 - 1 \right) \left( 4\nu^2-9 \right) \cdot \frac{1}{z^2} + \gamma_3 \left( \nu, z \right) \right], \\[1em]
                    
                        K_{\nu}' \left( z \right) & = & - \sqrt{\frac{\pi}{2z}} e^{-z} \left[ 1 + \frac{1}{8} \left( 4\nu^2 + 3 \right) \cdot \frac{1}{z} \right. \\[0.5em]
                    
                        && \qquad \qquad \qquad + \left. \frac{1}{128} \left( 4\nu^2 - 1 \right) \left( 4\nu^2 + 15 \right) \cdot \frac{1}{z^2} + \gamma_{3,1} \left( \nu, z \right) \right]. 
                    \end{array}
                \end{equation}
                where $\sqrt{\cdot}$ is the principal branch of the complex square root. Furthermore, the remainders $\gamma_3$ and $\gamma_{3,1}$ satisfy, for $\left \vert \nu \right \vert \geqslant A$ and $A > 0$ some constant, the estimates
                \begin{equation}
                    \begin{array}{llll}
                        \left \vert \gamma_3 \left( \nu, z \right) \right \vert & \leqslant & C_1 \cdot \frac{1}{\left \vert z \right \vert^3} \left \vert \nu \right \vert^6 e^{\frac{1}{\left \vert z \right \vert} \left \vert \nu^2 - \frac{1}{4} \right \vert} & \text{with } C_1 > 0 \\[1em]
                    
                        \left \vert \gamma_{3,1} \left( \nu, z \right) \right \vert & \leqslant & C_2 \frac{1}{\left \vert z \right \vert^3} \left \vert \nu \right \vert^6 e^{\frac{2}{\left \vert z \right \vert} \left \vert \nu^2 - \frac{1}{4} \right \vert} & \text{with } C_2 > 0
                    \end{array}
                \end{equation}
                where $\varepsilon > 0$ is a small enough real number, and $C_1$, $C_2$ are constants.
            \end{proposition}
        
            \begin{proof}
                The first asymptotic expansion of \eqref{eq:asymptoticsArgumentBessel} is given in \cite[Ex. 13.2]{olver}, while the second one can be obtained by applying Cauchy's formula and the first estimate, using a disk of radius $\varepsilon$ around $z$.
            \end{proof}
        
            \begin{corollary}
            \label{cor:asymptoticsLargeArgumentLogDerivBessel}
                For any complex numbers $\nu$ and $z$, with $\Re z > 0$, we have
                \begin{equation}
                    \begin{array}{lll}
                        z \frac{K_{\nu}' \left( z \right)}{K_{\nu} \left( z \right)} & = & - z - \frac{1}{2} - \frac{1}{8} \left( 4 \nu^2 - 1 \right) \cdot \frac{1}{z} + \gamma_{3,2} \left( \nu, z \right),
                    \end{array}
                \end{equation}
                where the remainder $\gamma_{3,2}$ satisfies the estimate
                \begin{equation}
                    \begin{array}{llll}
                        \left \vert \gamma_{3,2} \left( \nu, z \right) \right \vert & \leqslant & C_3 \frac{1}{\left \vert z \right \vert^2} \left \vert \nu \right \vert^6 e^{\frac{2}{\left \vert z \right \vert} \left \vert \nu^2 - \frac{1}{4} \right \vert} & \text{with } C_3 > 0,
                    \end{array}
                \end{equation}
                for $\left \vert \nu \right \vert \geqslant A$, where $C_3$ and $A > 0$ are constants.
            \end{corollary}
        
            \begin{proof}
                This is a direct consequence of proposition \ref{prop:largeArgumentAsymptoticsBessel}.
            \end{proof}
    
            \begin{corollary}
            \label{cor:asymptoticsLargeArgumentGK}
                For any real number $\mu \geqslant 0$ and any integer $k \in \mathbb{Z}$, with the exception of $k = 0$ should $\alpha$ vanish, we have, for any $\nu$ with $\left \vert \arg \nu \right \vert > \eta$ and $\eta > 0$,
                \begin{equation}
                \label{eq:asymptArgumentGk}
                    \begin{array}{lll}
                        g_k \left( \nu \right) & = & 4 \pi \alpha a - 2 \pi \left( \left \vert k + \alpha \right \vert + \left \vert k - \alpha \right \vert \right) a \\[0.5em]
                    
                        && \qquad \qquad \qquad + \frac{1}{16 \pi a} \left( 4 \nu^2 + 1 \right) \left( \frac{1}{\left \vert k + \alpha \right \vert} + \frac{1}{\left \vert k - \alpha \right \vert} \right) + \gamma_{3,3} \left( \nu, k \right)
                    \end{array}
                \end{equation}
                where the remainder $\gamma_{3,3}$ satisfies
                \begin{equation}
                    \begin{array}{lll}
                        \left \vert \gamma_{3,3} \left( \nu, k \right) \right \vert & \leqslant & \left \{ \begin{array}{ll} C_4 \cdot \frac{1}{k^2 a^2} \left \vert \nu \right \vert^6 e^{\frac{B}{\left \vert k \right \vert} \left \vert \nu^2 - \frac{1}{4} \right \vert} & \text{if } k \neq 0 \\[0.5em] C_4 \cdot \frac{1}{a^2} & \text{if } k = 0 \end{array} \right.
                    \end{array}
                \end{equation}
                for $\left \vert \nu \right \vert \geqslant A$, where $A > 0$, as well as $C_4 > 0$ and $B > 0$ are constants.
            \end{corollary}
        
            \begin{proof}
                This is a direct consequence of corollary \ref{cor:asymptoticsLargeArgumentLogDerivBessel}, noting that the function $g_k$ was defined in \eqref{eq:functiong}.
            \end{proof}

        \subsubsection{Auxiliary functions}
        
            The second ingredient needed to let $\vartheta$ go to $\pi/2$ in the AW contributions \eqref{eq:zetaDandAWKNeq0} and \eqref{eq:zetaDandAWK0} is a regularization of the integrated function. To do that, let us define an auxiliary function, as well as choose one of its primitive.
            
            \begin{definition}
            \label{def:fmuk}
                For any integer $k \in \mathbb{Z}$, except $k=0$ if we have $\alpha = 0$, we define the function $h_{\mu,k}$ on $\mathbb{C}$ by
                \begin{equation}
                \label{eq:functionHMuK}
                    \begin{array}{lll}
                        h_{\mu,k} \left( t \right) & = & i \left( \frac{\partial}{\partial t} \log g_k \right) \left( it \right) - \frac{2it}{\sqrt{4 \mu + 1}} \frac{\partial}{\partial t}_{\vert t = i\sqrt{\frac{1}{4} + \mu}} \log g_k \left( t \right),
                    \end{array}
                \end{equation}
                the function $g_k$ having been defined in \eqref{eq:functiong}. We further define, for $t \in \mathbb{R}$,
                \begin{equation}
                \label{eq:choicePrimitiveHMuk}
                    \begin{array}{lll}
                        H_{\mu,k} \left( t \right) & = & \log \left \vert g_k \left( it \right) \right \vert - \log \left \vert g_k \left( i\sqrt{\frac{1}{4} + \mu} \, \right) \right \vert \\[1em]
                        
                        && \qquad \qquad \qquad \qquad \quad - i \frac{t^2 - \left( 1/4 + \mu \right)}{\sqrt{4 \mu + 1}} \frac{\partial}{\partial t}_{\vert t = i\sqrt{\frac{1}{4} + \mu}} \log \left \vert g_k \left( t \right) \right \vert,
                    \end{array}
                \end{equation}
                which is a primitive of $h_{\mu,k}$ on $\mathbb{R}$.
            \end{definition}
            
            \begin{proposition}
            \label{prop:HmukBound}
                For any integer $k \in \mathbb{Z}$, we can write
                \begin{equation}
                    \begin{array}{lll}
                        H_{\mu,k} \left( t \right) & = & \left( \frac{1}{4} + \mu - t^2 \right)^2 \widetilde{H}_{\mu,k} \left( t \right)
                    \end{array}
                \end{equation}
                where $\widetilde{H}_{\mu,k}$ satisfies
                \begin{equation}
                    \begin{array}{lll}
                        \left \vert \widetilde{H}_{\mu,k} \left( t \right) \right \vert & < & \left \{ \begin{array}{lll} \frac{C_{\mu}}{\left \vert k \right \vert^3} & \text{if } k \neq 0 & \text{for } \left \vert t \right \vert \leqslant \sqrt{\frac{1}{4}+\mu} \\[1em] C_{\mu, \eta} & \text{if } k = 0 & \text{for } \left \vert t \right \vert \leqslant \sqrt{\frac{1}{4}+\mu} \text{ and } \left \vert t - \frac{1}{2} \right \vert \geqslant \eta \end{array} \right.,
                    \end{array}
                \end{equation}
                where $C_{\mu, \eta} > 0$ is a constant which does not depend on $\arg t$, and $\eta > 0$ is some constant. Furthermore, we have
                \begin{equation}
                \label{eq:boundHmukKDelta}
                    \begin{array}{lll}
                        \left \vert \widetilde{H}_{\mu,k} \left( t \right) \right \vert & \hspace{-3pt} < \hspace{-3pt} & \left \{ \begin{array}{lll} \frac{C_{\mu}}{\left \vert k \right \vert^{3-6 \delta}} & \text{if } k \neq 0 & \hspace{-3pt} \text{for } \left \vert t \right \vert \in \left[ \sqrt{\frac{1}{4}+\mu}, 2 \left \vert k \right \vert^{\delta} \sqrt{\frac{1}{4}+\mu} \; \right], \\[1em] C_{\mu} & \text{if } k = 0 & \hspace{-3pt} \text{for } \left \vert t \right \vert \in \left[ \sqrt{\frac{1}{4}+\mu}, 2 \sqrt{\frac{1}{4}+\mu} \; \right], \end{array} \right.
                    \end{array}
                \end{equation}
                where $\delta > 0$ is a constant.
            \end{proposition}
            
            \begin{proof}
                The proof of this result follows the same argument as \cite[Prop. 3.14]{dutour:cuspDirichlet}, which itself was based on \cite[Cor. 6.4]{freixas-vonPippich:RROrbifold}. The only change is that one needs to use the asymptotics \eqref{eq:asymptArgumentGk}.
            \end{proof}
            
            \begin{remark}
                The parameter $\delta$ will serve to control the convergence of certain series, and was first introduced in \cite[Sec. 6.1]{freixas-vonPippich:RROrbifold}. Its precise value will be controlled throughout this paper by a number of inequalities. For now, it is enough to think of it as a ``small enough'' constant.
            \end{remark}

        \subsubsection{The function $\zeta_{L,\mu}^{\pm, \, AW}$}
        \label{subsubsec:functionZetaAWPMThetaPiOver2}

        Let us first deal with the zeta function $\zeta_{L,\mu}^{\pm, \, AW}$, obtained in \eqref{eq:zetaDandAWKNeq0} by integration along the contour $\gamma_{\vartheta,\pm}$.
            
            \begin{proposition}
            \label{prop:intRepZetaAWPMTheta}
                For any real number $\mu \geqslant 0$, we have
                \begin{equation}
                \label{eq:integralRepHMuK}
                    \begin{array}{lll}
                        \zeta_{L,\mu}^{\pm, \, AW} \left( s \right) & = & \frac{1}{2 i \pi} \sum\limits_{k \neq 0} \int_{i \gamma_{\vartheta,\pm}} \left( \frac{1}{4} - t^2 + \mu \right)^{-s} h_{\mu,k} \left( t \right) \mathrm{d}t
                    \end{array}
                \end{equation}
                on the half-plane $\Re s > 1$.
            \end{proposition}
            
            \begin{proof}
                This is a consequence of the fact that we have
                \begin{equation}
                \label{eq:vanishingIntegral}
                    \begin{array}{lll}
                        \int_{i \gamma_{\vartheta,\pm}} \left( \frac{1}{4} - t^2 + \mu \right)^{-s} t \mathrm{d}t & = & 0,
                    \end{array}
                \end{equation}
                which comes from a direct computation of the primitive involved.
            \end{proof}
            
            Adding the vanishing integral coming from the second term on the right-hand side of \eqref{eq:functionHMuK} to get the integral representation \eqref{eq:integralRepHMuK} has a regularizing effect, allowing us to get a result at angle $\vartheta = \pi/2$. As in \cite[Sec. 3.2]{dutour:cuspDirichlet}, the contour $\gamma_{\vartheta, \pm}$ is split into four parts, according to where $1/4+\mu-t^2$ is when $t$ goes through $i \gamma_{\vartheta, \pm}$.
            
            \begin{definition}
                The four paths of integration $\gamma_{\vartheta, \pm}^{\left( 1 \right)}$, \ldots, $\gamma_{\vartheta, \pm}^{\left( 4 \right)}$ are defined by:
                \begin{equation}    
                    \begin{array}{llllll}
                        \gamma_{\vartheta, \pm}^{\left( 1 \right)} & = & \left \{ r e^{i \vartheta}, \; r \, \geqslant \, \sqrt{\frac{1}{4} + \mu} \right \}, &  \gamma_{\vartheta, \pm}^{\left( 2 \right)} & = & \left \{ r e^{-i \vartheta}, \; r \, \geqslant \, \sqrt{\frac{1}{4} + \mu} \right \}, \\[1em]
                    
                        \gamma_{\vartheta, \pm}^{\left( 3 \right)} & = & \left \{ r e^{i \vartheta}, \; r \, < \, \sqrt{\frac{1}{4} + \mu} \right \}, &  \gamma_{\vartheta, \pm}^{\left( 4 \right)} & = & \left \{ r e^{- i \vartheta}, \; r \, < \, \sqrt{\frac{1}{4} + \mu} \right \},
                    \end{array}
                \end{equation}
                oriented according to figure \ref{fig:integrationContoursKNeq0}. These definitions include the limit case $\vartheta = \pi/2$.
            \end{definition}
            
            The convergence of the four integrals as $\theta$ goes to $\pi/2$ can be  $k \neq 0$. To that effect, for any angle $0 < \vartheta < \pi/2$, we set
            \begin{equation}
                \begin{array}{lll}
                    \mathcal{C}_{\vartheta} & = & \left \{ \sqrt{\frac{1}{4} + \mu} \; e^{i \varphi}, \; \varphi \in \left[ \vartheta, \frac{\pi}{2} \right[ \right \}.
                \end{array}
            \end{equation}
            This path of integration is represented below.
            
            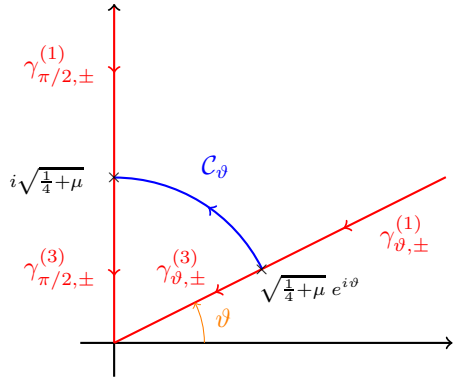
\begin{figure}[H]
    
                \begin{center}
        
                    \begin{tikzpicture}[scale=0.9]
        
                        \draw[thick] (-0.5,0) -- (0,0);
                        \draw[thick] (0,-0.5) -- (0,0);
                        \draw[thick,->] (0,0) coordinate (2) node {} -- (5,0) coordinate (1) node {};
                        \draw[thick, red] (0,0) -- (0,5) coordinate (4) node{};
                        \draw[thick, ->] (0,5) -- (0,5.01) coordinate (5) node{};
            
                        \draw[thick, red] (0,0) -- (4.9,2.45) coordinate (3) node {};
            
                        \draw pic["$\;\color{orange} \vartheta$", draw=orange, ->, angle eccentricity=1.2, angle radius=1.2cm] {angle=1--2--3};
                        
                        \draw pic["", draw=blue, thick, -, angle eccentricity=1, angle radius=2.2cm] {angle=3--2--4};
                        \node (R) at (2.9,0.8) {$\scriptstyle \sqrt{\frac{1}{4} + \mu} \; e^{i \vartheta}$};
                        \node (S) at (2.18,1.09) {$\scriptstyle \times$};
                        \node (T) at (-1,2.4) {$\scriptstyle i \sqrt{\frac{1}{4} + \mu}$};
                        \node (U) at (0,2.44) {$\scriptstyle \times$};
                        
                        \draw[thick, red, ->] (3.4,1.7) -- (3.39,1.695) coordinate (6) node {};
                        \node (V) at (4.3,1.6) {$\color{red} \gamma_{\vartheta,\pm}^{\left( 1 \right)}$};
                        
                        \draw[thick, red, ->] (1.5,0.75) -- (1.49,0.745) coordinate (7) node {};
                        \node (W) at (1,1.1) {$\color{red} \gamma_{\vartheta,\pm}^{\left( 3 \right)}$};
                        
                        \draw[thick, red, ->] (0,4) -- (0,3.99) coordinate (8) node {};
                        \node (X) at (-0.8,4.1) {$\color{red} \gamma_{\pi/2,\pm}^{\left( 1 \right)}$};
                        
                        \draw[thick, red, ->] (0,1) -- (0,0.99) coordinate (9) node {};
                        \node (Y) at (-0.8,1.1) {$\color{red} \gamma_{\pi/2,\pm}^{\left( 3 \right)}$};
                        
                        \node (Z) at (1.5,2.6) {$\color{blue} \mathcal{C}_{\vartheta}$};
                        \draw[thick, blue, ->] (1.4,2) -- (1.39,2.007) coordinate (10) node {};
                    \end{tikzpicture}
                    
                \end{center}
        
                \caption{Circular arc between $\vartheta$ and $\pi/2$}
                \label{fig:integrationContoursTwoAngles}
        
            \end{figure}
            
            \noindent Similarly, we define the circular arc $\mathcal{C}_{-\vartheta}$ between the angles $- \vartheta$ and $-\pi/2$. In the following, regarding contour integration, we will use the symbols $+$ and $-$ respectively to mean a disjoint union of contours, with preserved or reversed orientation.
            
            \begin{lemma}
            \label{lem:integralContours13ThetaPiOver2}
                On the strip $1 < \Re s < 2$, we have
                \begin{equation}
                \label{eq:resFormulaThetaPiOver2Contour1}
                    \begin{array}{lll}
                        e^{-i s \pi} \int_{i\gamma_{\pi/2,\pm}^{\left( 1 \right)} - i\mathcal{C}_{\vartheta} - i\gamma_{\vartheta,\pm}^{\left( 1 \right)}} \left( t^2 - \left( \frac{1}{4} + \mu \right) \right)^{-s} h_{\mu,k} \left( t \right) \mathrm{d}t & = & 0
                    \end{array}
                \end{equation}
                for any integer $k \neq 0$. On the same strip and for the same integers, we also have
                \begin{equation}
                \label{eq:resFormulaThetaPiOver2Contour3}
                    \begin{array}{lll}
                        \int_{i\gamma_{\pi/2,\pm}^{\left( 3 \right)} + i\mathcal{C}_{\vartheta} - i\gamma_{\vartheta,\pm}^{\left( 3 \right)}} \left( \frac{1}{4} - t^2 + \mu \right)^{-s} h_{\mu,k} \left( t \right) \mathrm{d}t & = & 0.
                    \end{array}
                \end{equation}
            \end{lemma}
            
            \begin{proof}
                This is a direct application of the residue formula. Note that its application for the unbounded contour considered in \eqref{eq:resFormulaThetaPiOver2Contour1} follows the same argument as in paragraph \ref{subsubsec:integralRepresentation}.
            \end{proof}
            
            Let us state the corresponding lemma for the contour $\mathcal{C}_{- \vartheta}$.
            
            \begin{lemma}
                On the strip $1 < \Re s < 2$, we have
                \begin{equation}
                \label{eq:resFormulaThetaPiOver2Contour2}
                    \begin{array}{lll}
                        e^{i s \pi} \int_{i\gamma_{\pi/2,\pm}^{\left( 2 \right)} - i\mathcal{C}_{-\vartheta} - i\gamma_{\vartheta,\pm}^{\left( 2 \right)}} \left( t^2 - \left( \frac{1}{4} + \mu \right) \right)^{-s} h_{\mu,k} \left( t \right) \mathrm{d}t & = & 0
                    \end{array}
                \end{equation}
                for any integer $k \neq 0$. On the same strip and for the same integers, we also have
                \begin{equation}
                \label{eq:resFormulaThetaPiOver2Contour4}
                    \begin{array}{lll}
                        \int_{i\gamma_{\pi/2,\pm}^{\left( 4 \right)} + i\mathcal{C}_{-\vartheta} - i\gamma_{\vartheta,\pm}^{\left( 4 \right)}} \left( \frac{1}{4} - t^2 + \mu \right)^{-s} h_{\mu,k} \left( t \right) \mathrm{d}t & = & 0.
                    \end{array}
                \end{equation}
            \end{lemma}
            
            \begin{proof}
                This is similar to lemma \ref{lem:integralContours13ThetaPiOver2}.
            \end{proof}
            
            \begin{remark}
                It should be noted that the exponentials in factor of the integrals in \eqref{eq:resFormulaThetaPiOver2Contour1} and \eqref{eq:resFormulaThetaPiOver2Contour2} are explicitely written so as to facilitate the combination of these results. They stem from the equalities
                \begin{equation}
                \label{eq:limitThetaPi2Gamma12}
                    \begin{array}{lll}
                        \left( \frac{1}{4} + \mu - t^2 \right)^{-s} & = & \left \{ \begin{array}{lcl} e^{-is\pi} \left( t^2 - \left( \frac{1}{4} + \mu \right) \right)^{-s} & \text{on} & i \gamma_{\vartheta, \pm}^{\left( 1 \right)} \cup i \mathcal{C}_{\vartheta} \\[0.8em] e^{is\pi} \left( t^2 - \left( \frac{1}{4} + \mu \right) \right)^{-s} & \text{on} & i \gamma_{\vartheta, \pm}^{\left( 2 \right)} \cup i \mathcal{C}_{-\vartheta} \end{array} \right.
                    \end{array}
                \end{equation}
                which allow us to give meaning to the complex powers even at $\vartheta = \pi/2$.
            \end{remark}
            
            We will now add the vanishing contributions \eqref{eq:resFormulaThetaPiOver2Contour3} and \eqref{eq:resFormulaThetaPiOver2Contour4} to the integral representation \eqref{eq:integralRepHMuK}.
            
            \begin{lemma}
            \label{lem:intRepZetaPMAWPathsAndCircularArcs}
                For any real number $\mu \geqslant 0$, we have
                \begin{equation}
                \label{eq:intRepZetaPMAWPathsAndCircularArcs}
                    \begin{array}{lll}
                        \zeta_{L,\mu}^{\pm, \, AW} \left( s \right) & = & \frac{1}{2 i \pi} \sum\limits_{k \neq 0} \left[ \int_{i \gamma_{\vartheta,\pm}^{\left( 1 \right)} + i \gamma_{\vartheta,\pm}^{\left( 2 \right)}} \left( \frac{1}{4} - t^2 + \mu \right)^{-s} h_{\mu,k} \left( t \right) \mathrm{d}t \right. \\[1em]
                        
                        && \qquad \qquad \qquad \left. +\int_{i \mathcal{C}_{\vartheta} + i \mathcal{C}_{-\vartheta}} \left( \frac{1}{4} - t^2 + \mu \right)^{-s} h_{\mu,k} \left( t \right) \mathrm{d}t \right]
                    \end{array}
                \end{equation}
                on the strip $1 < \Re s < 2$.
            \end{lemma}
            
            \begin{proof}
                First, we need to add the vanishing contributions \eqref{eq:resFormulaThetaPiOver2Contour3} and \eqref{eq:resFormulaThetaPiOver2Contour4} to the integral representation \eqref{eq:integralRepHMuK}. This yields
                \begin{equation}
                    \begin{array}{lll}
                        \zeta_{L,\mu}^{\pm, \, AW} \left( s \right) & = & \frac{1}{2 i \pi} \sum\limits_{k \neq 0} \left[ \int_{i \gamma_{\vartheta,\pm}^{\left( 1 \right)} + i \gamma_{\vartheta,\pm}^{\left( 2 \right)}} \left( \frac{1}{4} - t^2 + \mu \right)^{-s} h_{\mu,k} \left( t \right) \mathrm{d}t \right. \\[1em]
                        
                        && \qquad + \int_{i \mathcal{C}_{\vartheta} + i \mathcal{C}_{-\vartheta}} \left( \frac{1}{4} - t^2 + \mu \right)^{-s} h_{\mu,k} \left( t \right) \mathrm{d}t \\[0.5em]
                        
                        && \qquad \qquad \left. + \int_{i \gamma_{\pi/2,\pm}^{\left( 3 \right)} + i \gamma_{\pi/2,\pm}^{\left( 4 \right)}} \left( \frac{1}{4} - t^2 + \mu \right)^{-s} h_{\mu,k} \left( t \right) \mathrm{d}t \right],
                    \end{array}
                \end{equation}
                noting that the contributions of the contours $i \gamma_{\vartheta,\pm}^{\left( 3 \right)}$ and $i \gamma_{\vartheta,\pm}^{\left( 4 \right)}$ are canceled, as they appear twice with reversed orientations. Finally, we have
                \begin{equation}
                    \begin{array}{lllll}
                        \multicolumn{5}{l}{\int_{i \gamma_{\pi/2,\pm}^{\left( 3 \right)} + i \gamma_{\pi/2,\pm}^{\left( 4 \right)}} \left( \frac{1}{4} - t^2 + \mu \right)^{-s} h_{\mu,k} \left( t \right) \mathrm{d}t} \\[1em]
                        
                        \qquad \qquad \qquad \qquad & = & \int_{- \sqrt{1/4+\mu}}^{\sqrt{1/4+\mu}} \left( \frac{1}{4} - t^2 + \mu \right)^{-s} h_{\mu,k} \left( t \right) \mathrm{d}t & = & 0,
                    \end{array}
                \end{equation}
                because of the oddness of the function $h_{\mu,k}$. This concludes the proof.
            \end{proof}
            
            We can now obtain the integral representation of $\zeta_{L,\mu}^{\pm, \, AW}$ at angle $\pi/2$.
            
            \begin{theorem}
            \label{thm:intRepZetaFunctionPMAW}
                In the strip $1 < \Re s < 2$, we have
                \begin{equation}
                \label{eq:intRepZetaFunctionPMAW}
                    \begin{array}{lll}
                        \zeta_{L,\mu}^{\pm, \, AW} \left( s \right) & = & \frac{1}{\pi} \sin \left( \pi s \right) \sum\limits_{k \neq 0} \int_{\sqrt{\frac{1}{4} + \mu}}^{+ \infty} \left( t^2 - \left( \frac{1}{4} + \mu \right) \right)^{-s} h_{\mu,k} \left( t \right) \mathrm{d}t.
                    \end{array}
                \end{equation}
            \end{theorem}
            
            \begin{proof}
                This is a consequence of the integral representation \eqref{eq:intRepZetaPMAWPathsAndCircularArcs} from lemma \ref{lem:intRepZetaPMAWPathsAndCircularArcs}, to which we add the vanishing contributions \eqref{eq:resFormulaThetaPiOver2Contour1} and \eqref{eq:resFormulaThetaPiOver2Contour2}.
            \end{proof}

        \subsubsection{The function $\zeta_{L,\mu}^{0, \, AW}$}
        
            Let us now deal with the zeta function $\zeta_{L,\mu}^{0, \, AW}$, obtained in \eqref{eq:zetaDandAWK0} by integration along the contour $\gamma_{\vartheta,0}$. Throughout this paragraph, we assume that we have $\alpha \neq 0$.
            
            \begin{proposition}
            \label{prop:intRepZetaAW0Theta}
                For any real number $\mu > 0$, we have
                \begin{equation}
                \label{eq:integralRepHMu0}
                    \begin{array}{lll}
                        \zeta_{L,\mu}^{0, \, AW} \left( s \right) & = & \frac{1}{2 i \pi} \int_{i \gamma_{\vartheta,0}} \left( \frac{1}{4} - t^2 + \mu \right)^{-s} h_{\mu,0} \left( t \right) \mathrm{d}t
                    \end{array}
                \end{equation}
                on the half-plane $\Re s > 1$.
            \end{proposition}
            
            \begin{proof}
                This is a consequence of the fact that we have
                \begin{equation}
                    \begin{array}{lll}
                        \int_{i \gamma_{\vartheta,0}} \left( \frac{1}{4} - t^2 + \mu \right)^{-s} t \mathrm{d}t & = & 0,
                    \end{array}
                \end{equation}
                which comes from a direct computation of the primitive involved.
            \end{proof}
            
            As in paragraph \ref{subsubsec:functionZetaAWPMThetaPiOver2}, the integral representation \eqref{eq:integralRepHMu0} is better than \eqref{eq:zetaDandAWK0}, insofar as it will make sense at angle $\vartheta = \pi/2$. To that effect, we will split the countour $\gamma_{\vartheta,0}$ into eight parts. Let us consider a parameter $\varepsilon > 0$ small enough.
            
            \begin{definition}
                The eight paths of integration $\gamma_{\vartheta, 0}^{\left( 1 \right)}$, \ldots, $\gamma_{\vartheta, 0}^{\left( 8 \right)}$ are defined by:
                \begin{equation}    
                    \begin{array}{llllll}
                        \gamma_{\vartheta, 0}^{\left( 1 \right)} & = & \left \{ r e^{i \vartheta}, \; r \, \geqslant \, \sqrt{\frac{1}{4} + \mu} \right \}, &  \gamma_{\vartheta, 0}^{\left( 2 \right)} & = & \left \{ r e^{-i \vartheta}, \; r \, \geqslant \, \sqrt{\frac{1}{4} + \mu} \right \}, \\[1em]
                    
                        \gamma_{\vartheta, 0}^{\left( 3 \right)} & = & \multicolumn{4}{l}{\left \{ r e^{i \vartheta}, \; \sqrt{\frac{1}{4} + \varepsilon \mu} \, \leqslant \, r \, < \, \sqrt{\frac{1}{4} + \mu} \right \},} \\[1em]
                        
                        \gamma_{\vartheta, 0}^{\left( 4 \right)} & = & \multicolumn{4}{l}{\left \{ r e^{-i \vartheta}, \; \sqrt{\frac{1}{4} + \varepsilon \mu} \, \leqslant \, r \, < \, \sqrt{\frac{1}{4} + \mu} \right \},} \\[1em]
                        
                        \gamma_{\vartheta, 0}^{\left( 5 \right)} & = & \multicolumn{4}{l}{\left \{ \frac{1}{2} e^{i \vartheta} + \left( \sqrt{\frac{1}{4} + \mu} - \frac{1}{2} \right) e^{i \varphi}, \; \varphi \in \left[ -\pi, 0 \right[ \right \},} \\[1em]
                        
                        \gamma_{\vartheta, 0}^{\left( 6 \right)} & = & \multicolumn{4}{l}{\left \{ \frac{1}{2} e^{-i \vartheta} + \left( \sqrt{\frac{1}{4} + \mu} - \frac{1}{2} \right) e^{i \varphi}, \; \varphi \in \left[ 0, \pi \right[ \right \},} \\[1em]
                        
                        \gamma_{\vartheta, 0}^{\left( 7 \right)} & = & \multicolumn{4}{l}{\left \{ r e^{i \vartheta}, \; r \, < \, 1 - \sqrt{\frac{1}{4} + \varepsilon \mu} \right \},} \\[1em]
                        
                        \gamma_{\vartheta, 0}^{\left( 8 \right)} & = & \multicolumn{4}{l}{\left \{ r e^{-i \vartheta}, \; r \, < \, 1 - \sqrt{\frac{1}{4} + \varepsilon \mu} \right \},}
                    \end{array}
                \end{equation}
                oriented according to figure \ref{fig:integrationContourK0}.
            \end{definition}
            
            As we did in paragraph \ref{subsubsec:functionZetaAWPMThetaPiOver2}, we obtain the following result.
            
            \begin{lemma}
            \label{lem:intRepZetaAW0WithSemiCircles}
                For any real number $\mu > 0$, we have
                \begin{equation}
                \label{eq:integralRepHMu0PiOver2}
                    \begin{array}{lll}
                        \zeta_{L,\mu}^{0, \, AW} \left( s \right) & = & \frac{1}{\pi} \sin \left( \pi s \right) \int_{\sqrt{\frac{1}{4} + \mu}}^{+ \infty} \left( t^2 - \left( \frac{1}{4} + \mu \right) \right)^{-s} h_{\mu,0} \left( t \right) \mathrm{d}t \\[1em]
                        
                        && \qquad + \frac{1}{2 i \pi} \int_{i \gamma_{\pi/2,0}^{\left( 5 \right)}} \hspace{-1pt} \left( \frac{1}{4} - t^2 + \mu \right)^{-s} h_{\mu,0} \left( t \right) \mathrm{d}t \\[1em]
                        
                        && \qquad \qquad + \frac{1}{2 i \pi} \int_{i \gamma_{\pi/2,0}^{\left( 6 \right)}} \hspace{-1pt} \left( \frac{1}{4} - t^2 + \mu \right)^{-s} h_{\mu,0} \left( t \right) \mathrm{d}t
                    \end{array}
                \end{equation}
                on the strip $1 < \Re s < 2$.
            \end{lemma}
            
            \begin{proof}
                The argument is similar to the one used for theorem \ref{thm:intRepZetaFunctionPMAW}.
            \end{proof}
            
            The last two integrals on the right-hand side of \eqref{eq:integralRepHMu0PiOver2} cannot be easily computed. However, they depend on the parameter $\varepsilon > 0$, which controls the radius of the semi-circles centered at $\pm i/2$. Since no other term in \eqref{eq:integralRepHMu0PiOver2} depends on $\varepsilon$, we can consider the limit as $\varepsilon$ goes to $0^+$.
            
            \begin{lemma}
            \label{lem:limitEpsilon0IntSemiCircles}
                On the strip $1 < \Re s < 2$, we have the following limit
                \begin{equation}
                \label{eq:limitEpsilon0IntSemiCircleUp}
                    \begin{array}{lll}
                        \frac{1}{2 i \pi} \int_{i \gamma_{\pi/2,0}^{\left( 5 \right)}} \hspace{-1pt} \left( \frac{1}{4} - t^2 + \mu \right)^{-s} h_{\mu,0} \left( t \right) \mathrm{d}t & \underset{\varepsilon \rightarrow 0^+}{\longrightarrow} & - \frac{1}{2} \mu^{-s}.
                    \end{array}
                \end{equation}
                Similarly, we have
                \begin{equation}
                \label{eq:limitEpsilon0IntSemiCircleDown}
                    \begin{array}{lll}
                        \frac{1}{2 i \pi} \int_{i \gamma_{\pi/2,0}^{\left( 6 \right)}} \hspace{-1pt} \left( \frac{1}{4} - t^2 + \mu \right)^{-s} h_{\mu,0} \left( t \right) \mathrm{d}t & \underset{\varepsilon \rightarrow 0^+}{\longrightarrow} & - \frac{1}{2} \mu^{-s}.
                    \end{array}
                \end{equation}
            \end{lemma}
            
            \begin{proof}
                This is a consequence of the fact that the function
                \begin{equation}
                    \begin{array}{lll}
                        t & \longmapsto & \left( \frac{1}{4} - t^2 + \mu \right)^{-s} h_{\mu,0} \left( t \right)
                    \end{array}
                \end{equation}
                has a simple pole at $t=\pm 1/2$, for any $s$ such that we have $1 < \Re s < 2$, and any real number $\mu > 0$. The limits \eqref{eq:limitEpsilon0IntSemiCircleUp} and \eqref{eq:limitEpsilon0IntSemiCircleDown} may then computed using Laurent series expansions.
            \end{proof}
            
            We can now give a simpler integral representation of $\zeta_{L,\mu}^{0, \, AW}$.
            
            \begin{theorem}
            \label{thm:intRepZetaFunction0AW}
                For any real number $\mu > 0$, we have
                \begin{equation}
                \label{eq:integralRepZetaAW0}
                    \begin{array}{lll}
                        \zeta_{L,\mu}^{0, \, AW} \left( s \right) & = & - \mu^{-s} + \frac{1}{\pi} \sin \left( \pi s \right) \int_{\sqrt{\frac{1}{4} + \mu}}^{+ \infty} \left( t^2 - \left( \frac{1}{4} + \mu \right) \right)^{-s} h_{\mu,0} \left( t \right) \mathrm{d}t
                    \end{array}
                \end{equation}
                on the strip $1 < \Re s < 2$.
            \end{theorem}
            
            \begin{proof}
                This is a consequence of lemmas \ref{lem:intRepZetaAW0WithSemiCircles} and \ref{lem:limitEpsilon0IntSemiCircles}.
            \end{proof}

\section{The Alvarez--Wentworth contributions}

    In the previous section, we began the study of the spectral zeta function $\zeta_{L,\mu}$ with two distinct results:
    \begin{itemize}
        \item In subsection \ref{subsec:factorizationDirichletContribution}, we extracted Dirichlet-like contributions, in the form of the partial spectral zeta functions $\zeta_{L,\mu}^{\pm, \, D}$ and $\zeta_{L,\mu}^{0, \, D}$. These were proved, in \cite{dutour:cuspDirichlet}, to have holomorphic continuations to a domain containing $0$, and asymptotic expansions of their derivative at $s=0$ were obtained. This was recalled in theorems \ref{thm:asymptDerAt0DirichletContributionPM} and \ref{thm:asymptDerAt0DirichletContribution0}.
        
        \item The remaining partial spectral zeta functions $\zeta_{L,\mu}^{\pm, \, AW}$ and $\zeta_{L,\mu}^{0, \, AW}$, which constitute the core of this paper, require more care. An integral representation of these functions was found in theorems \ref{thm:intRepZetaFunctionPMAW} and \ref{thm:intRepZetaFunction0AW}.
    \end{itemize}
    
    This section is devoted to continuing the study of the Alvarez--Wentworth contributions $\zeta_{L,\mu}^{\pm, \, AW}$ and $\zeta_{L,\mu}^{0, \, AW}$. More precisely, we will use the integral representations of these partial spectral zeta functions to prove that they can be holomorphically continued near $0$, and obtain four asymptotic expansions. The first three concern:
    
    \begin{enumerate}
    \itemsep0.5em
        
        \item $\zeta_{L,\mu}^{\pm, \, AW} {}' \left( 0 \right)$ as $\mu$ goes to infinity, for any real number $a > 1/\left( 4 \pi \left( 1- \alpha \right) \right)$, this hypothesis being necessary in proposition \ref{prop:zerosfk};
        
        \item $\zeta_{L,0}^{\pm, \, AW} {}' \left( 0 \right)$ as $a$ goes to infinity;
        
        \item $\zeta_{L,\mu}^{0, \, AW} {}' \left( 0 \right)$ as $\mu$ goes to infinity, for any real number $a > 0$.
    \end{enumerate}
    
    The last asymptotic expansion we need is more complicated, since the spectral zeta function $\zeta_{L,\mu}^{0, \, AW}$ is not even defined when we have $\mu = 0$. However, we will set
    \begin{equation}
        \begin{array}[t]{lll}
            \zeta_{L,0}^{0, \, AW} {}' \left( 0 \right) & = & \lim\limits_{\mu \rightarrow 0^+} \zeta_{L,\mu}^{0, \, AW} {}' \left( 0 \right),
        \end{array}
    \end{equation}
    as an abuse of notation, after proving that this limit finite. This is closely related to the notion of \textit{modified} determinant, where an operator is restricted to the orthogonal complement of its kernel before considering its spectral zeta function, thereby excluding the eigenvalue $0$. The remaining asymptotic expansion concerns:
    
    \begin{enumerate}
        \setcounter{enumi}{3}
        \item $\zeta_{L,0}^{0, \, AW} {}' \left( 0 \right)$ as $a$ goes to infinity.
    \end{enumerate}
    
    To achieve these purposes, we will need to perform technical and lengthy computations, akin to \cite[Sec. 3]{dutour:cuspDirichlet}, which themselves used (in part) techniques developped by Freixas i Montplet and von Pippich in \cite{freixas-vonPippich:RROrbifold}. Finally, let us mention that the reader will find an overview of the methods used here in the introduction.

    \subsection{Splitting the interval of integration}
    \label{subsec:splittingIntervalIntegration}
        
        The two partial spectral zeta functions $\zeta_{L,0}^{0, \, AW}$ and $\zeta_{L,0}^{\pm, \, AW}$ have been written, in theorems \ref{thm:intRepZetaFunctionPMAW} and \ref{thm:intRepZetaFunction0AW}, as an integral and a sum of integrals respectively. In order to study these, we need to split the interval of integration, similarly to \cite[Sec. 3.3]{dutour:cuspDirichlet} and \cite[Sec. 6.1]{freixas-vonPippich:RROrbifold}.
    
        \subsubsection{The function $\zeta_{L,0}^{0, \, AW}$}
        \label{subsubsec:splittingK0}
        
            Let us begin with the simplest of the two contributions, which does not involve a series, but a single integral. Let us recall the statement of theorem \ref{thm:intRepZetaFunction0AW}. In the strip $1 < \Re s < 2$, we have, for any real number $\mu > 0$,
            \begin{equation}
                \begin{array}{lll}
                    \zeta_{L,\mu}^{0, \, AW} \left( s \right) & = & - \mu^{-s} + \frac{1}{\pi} \sin \left( \pi s \right) \int_{\sqrt{\frac{1}{4} + \mu}}^{+ \infty} \left( t^2 - \left( \frac{1}{4} + \mu \right) \right)^{-s} h_{\mu,0} \left( t \right) \mathrm{d}t.
                \end{array}
            \end{equation}
            It should be noted that, in this paragraph, we assume that we have $\alpha \neq 0$.
            
            \begin{definition}
                For any real number $\mu > 0$, we define the integral $I_{\mu,0}$ by
                \begin{equation}
                    \begin{array}{lll}
                        I_{\mu,0} \left( s \right) & = & \frac{1}{\pi} \sin \left( \pi s \right) \int_{\sqrt{\frac{1}{4} + \mu}}^{+ \infty} \left( t^2 - \left( \frac{1}{4} + \mu \right) \right)^{-s} h_{\mu,0} \left( t \right) \mathrm{d}t
                    \end{array}
                \end{equation}
                in the strip $1 < \Re s < 2$.
            \end{definition}
            
            As discussed in the introduction of this subsection, we need to split the interval of integration. Let us write
            \begin{equation}
            \label{eq:splittingIntervalIntegrationK0}
                \begin{array}{lll}
                    \left] \sqrt{\frac{1}{4} + \mu}, \; + \infty \right[ & = & \left] \sqrt{\frac{1}{4} + \mu}, \; 2 \sqrt{\frac{1}{4} + \mu} \; \right] \cup \left] 2 \sqrt{\frac{1}{4} + \mu}, \; + \infty \right[,
                \end{array}
            \end{equation}
            and make the following definition, according to this splitting.
            
            \begin{definition}
                For any real number $\mu > 0$, we define
                \begin{equation}
                \label{eq:intLMuK0}
                    \begin{array}{lll}
                        L_{\mu,0} \left( s \right) & = & \frac{1}{\pi} \sin \left( \pi s \right) \int_{\sqrt{\frac{1}{4} + \mu}}^{2 \sqrt{\frac{1}{4} + \mu}} \left( t^2 - \left( \frac{1}{4} + \mu \right) \right)^{-s} h_{\mu,0} \left( t \right) \mathrm{d}t
                    \end{array}
                \end{equation}
                in the strip $1 < \Re s < 2$, where we also define
                \begin{equation}
                \label{eq:intMMuK0}
                    \begin{array}{lll}
                        M_{\mu,0} \left( s \right) & = & \frac{1}{\pi} \sin \left( \pi s \right) \int_{2 \sqrt{\frac{1}{4} + \mu}}^{+ \infty} \left( t^2 - \left( \frac{1}{4} + \mu \right) \right)^{-s} h_{\mu,0} \left( t \right) \mathrm{d}t.
                    \end{array}
                \end{equation}
                These integrals are related by $I_{\mu,0} \left( s \right) = L_{\mu,0} \left( s \right) + M_{\mu,0} \left( s \right)$.
            \end{definition}
        
        \subsubsection{The function $\zeta_{L,0}^{\pm, \, AW}$}
        \label{subsubsec:splittingKNeq0}
        
            Let us now move on to the remaining partial spectral zeta function. As stated in theorem \ref{thm:intRepZetaFunctionPMAW}, we have, in the strip $1 < \Re s < 2$,
            \begin{equation}
                \begin{array}{lll}
                    \zeta_{L,\mu}^{\pm, \, AW} \left( s \right) & = & \frac{1}{\pi} \sin \left( \pi s \right) \sum\limits_{k \neq 0} \int_{\sqrt{\frac{1}{4} + \mu}}^{+ \infty} \left( t^2 - \left( \frac{1}{4} + \mu \right) \right)^{-s} h_{\mu,k} \left( t \right) \mathrm{d}t,
                    \end{array}
            \end{equation}
            for any real number $\mu \geqslant 0$. 
            
            \begin{definition}
                For any real number $\mu \geqslant 0$ and any integer $k \neq 0$, we define the integral $I_{\mu,k}$ in the strip $1 < \Re s < 2$ by
                \begin{equation}
                    \begin{array}{lll}
                        I_{\mu,k} \left( s \right) & = & \frac{1}{\pi} \sin \left( \pi s \right) \int_{\sqrt{\frac{1}{4} + \mu}}^{+ \infty} \left( t^2 - \left( \frac{1}{4} + \mu \right) \right)^{-s} h_{\mu,k} \left( t \right) \mathrm{d}t.
                    \end{array}
                \end{equation}
            \end{definition}

            In order for the splitting of the interval of integration to help with the convergence of the series, it must depend on the integer $k$. This choice was foreshadowed in \eqref{eq:boundHmukKDelta}. For any integer $k \neq 0$, we write
            \begin{equation}
            \label{eq:splittingIntervalIntegrationKNeq0}
                \begin{array}{lll}
                    \left] \sqrt{\frac{1}{4}+\mu}, + \infty \; \right[ & = & \left] \sqrt{\frac{1}{4}+\mu}, \; 2 \left \vert k \right \vert^{\delta} \sqrt{\frac{1}{4}+\mu} \; \right] \cup \left] 2 \left \vert k \right \vert^{\delta} \sqrt{\frac{1}{4}+\mu}, + \infty \; \right].
                \end{array}
            \end{equation}
            It should be noted that the same splitting was performed in \cite[Sec. 3.3]{dutour:cuspDirichlet}, and that it is based on the one effected in \cite[Sec. 6.1]{freixas-vonPippich:RROrbifold}.
            
            \begin{definition}
                For any real number $\mu \geqslant 0$ and any integer $k \neq 0$, we define
                \begin{equation}
                \label{eq:intLMuKNeq0}
                    \begin{array}{lll}
                        L_{\mu,k} \left( s \right) & = & \frac{1}{\pi} \sin \left( \pi s \right) \int_{\sqrt{\frac{1}{4} + \mu}}^{2 \left \vert k \right \vert^{\delta} \sqrt{\frac{1}{4} + \mu}} \left( t^2 - \left( \frac{1}{4} + \mu \right) \right)^{-s} h_{\mu,k} \left( t \right) \mathrm{d}t
                    \end{array}
                \end{equation}
                in the strip $1 < \Re s < 2$, where we also define
                \begin{equation}
                \label{eq:intMMuKNeq0}
                    \begin{array}{lll}
                        M_{\mu,k} \left( s \right) & = & \frac{1}{\pi} \sin \left( \pi s \right) \int_{2 \left \vert k \right \vert^{\delta} \sqrt{\frac{1}{4} + \mu}}^{+ \infty} \left( t^2 - \left( \frac{1}{4} + \mu \right) \right)^{-s} h_{\mu,k} \left( t \right) \mathrm{d}t.
                    \end{array}
                \end{equation}
                These integrals are related by $I_{\mu,k} \left( s \right) = L_{\mu,k} \left( s \right) + M_{\mu,k} \left( s \right)$.
            \end{definition}

    \subsection{\texorpdfstring{General study of the integrals $L_{\mu,k}$}{General study of the integrals Lmuk}}
    \label{subsec:generalStudyLMuK}
        
        The first part of the integrals $I_{\mu,k}$ to be studied is given in the splitting above by $L_{\mu,k}$. Note that the notations for these integrals has been uniformized, even though the cases $k=0$ and $k \neq 0$ are involved in different partial spectral zeta functions, as seen in paragraphs \ref{subsubsec:splittingK0} and \ref{subsubsec:splittingKNeq0}.
    
        \subsubsection{The case $k=0$}
        \label{subsubsec:generalStudyLMuK0}
        
            Let us begin by considering, for any real number $\mu > 0$, the integral $L_{\mu,0}$, which was defined in \eqref{eq:intLMuK0}. In this paragraph, we assume that we have $\alpha \neq 0$.
            
            \begin{proposition}
                For any real number $\mu > 0$, we have
                \begin{equation}
                \label{eq:LMuK0IntByParts}
                    \begin{array}{lll}
                        L_{\mu,0} \left( s \right) & = & \frac{1}{\pi} \sin \left( \pi s \right) \cdot \frac{1}{3^s} \left( \frac{1}{4} + \mu \right)^{-s} H_{\mu,0} \left( 2 \sqrt{\frac{1}{4} + \mu} \right) \\[0.5em]
                        
                        && \qquad + \frac{2}{\pi} s \sin \left( \pi s \right) \int_{\sqrt{\frac{1}{4} + \mu}}^{2 \sqrt{\frac{1}{4} + \mu}} t \left( t^2 - \left( \frac{1}{4} + \mu \right) \right)^{-s-1} H_{\mu,0} \left( t \right) \mathrm{d}t
                    \end{array}
                \end{equation}
                in the strip $1 < \Re s < 2$.
            \end{proposition}
            
            \begin{proof}
                This proposition results from an integration by parts performed on \eqref{eq:intLMuK0}, using proposition \ref{prop:HmukBound} to justify the vanishing of part of the integrated component.
            \end{proof}
            
            Let us now study separately both parts of the right-hand side of \eqref{eq:LMuK0IntByParts}.
            
            \begin{definition}
                For any real number $\mu > 0$, we define
                \begin{equation}
                \label{eq:defAMuK0}
                    \begin{array}{lll}
                        A_{\mu,0} \left( s \right) & = & \frac{1}{\pi} \sin \left( \pi s \right) \cdot \frac{1}{3^s} \left( \frac{1}{4} + \mu \right)^{-s} H_{\mu,0} \left( 2 \sqrt{\frac{1}{4} + \mu} \right)
                    \end{array}
                \end{equation}
                in the strip $1 < \Re s < 2$, where we also define
                \begin{equation}
                \label{eq:defBMuK0}
                    \begin{array}{lll}
                        B_{\mu,0} \left( s \right) & = & \frac{2}{\pi} s \sin \left( \pi s \right) \int_{\sqrt{\frac{1}{4} + \mu}}^{2 \sqrt{\frac{1}{4} + \mu}} t \left( t^2 - \left( \frac{1}{4} + \mu \right) \right)^{-s-1} H_{\mu,0} \left( t \right) \mathrm{d}t.
                    \end{array}
                \end{equation}
            \end{definition}
            
            The integral $B_{\mu,0}$ seems more complicated, as it involves an integral whose lower bound is problematic for the complex power. However, as in \cite[Prop. 3.19]{dutour:cuspDirichlet}, the derivative at $s=0$ of this term will vanish.
            
            \begin{proposition}
            \label{prop:studyBMuK0}
                The function
                \begin{equation}
                    \begin{array}{lll}
                        s & \longmapsto & B_{\mu,0} \left( s \right)
                    \end{array}
                \end{equation}
                is holomorphic on the half-plane $\Re s < 2$, and its derivative at $s=0$ vanishes.
            \end{proposition}
            
            \begin{proof}
                This comes from the study of the integral in \eqref{eq:defBMuK0} with proposition \ref{prop:HmukBound}. The factor $s \sin \left( \pi s \right)$ then forces the vanishing of the derivative at $s=0$.
            \end{proof}

        \subsubsection{The case $k \neq 0$}
        \label{subsubsec:generalStudyLMukStudyBMuk}
        
            We will now consider the sum of the integrals $L_{\mu,k}$ over the integers $k \neq 0$, for any real number $\mu \geqslant 0$. Let us insist upon the fact that $\alpha$ is not assumed to be non-vanishing here.
            
            \begin{proposition}
                For any real number $\mu \geqslant 0$ and any integer $k \neq 0$, we have
                \begin{equation}
                \label{eq:LMuKNeq0IntByParts}
                    \begin{array}{lll}
                        L_{\mu,k} \left( s \right) & = & \frac{1}{\pi} \sin \left( \pi s \right) \cdot \left( 4 \mu + 1 \right)^{-s} \left( \left \vert k \right \vert^{2 \delta} - \frac{1}{4} \right)^{-s} H_{\mu,k} \left( 2 \left \vert k \right \vert^{\delta} \sqrt{\frac{1}{4} + \mu} \right) \\[0.8em]
                        
                        && \quad + \frac{2}{\pi} s \sin \left( \pi s \right) \int_{\sqrt{\frac{1}{4} + \mu}}^{2 \left \vert k \right \vert^{\delta} \sqrt{\frac{1}{4} + \mu}} t \left( t^2 - \left( \frac{1}{4} + \mu \right) \right)^{-s-1} H_{\mu,k} \left( t \right) \mathrm{d}t
                    \end{array}
                \end{equation}
                in the strip $1 < \Re s < 2$.                
            \end{proposition}
            
            \begin{proof}
                This proposition results from an integration by parts performed on \eqref{eq:intLMuKNeq0}, using proposition \ref{prop:HmukBound} to justify the vanishing of part of the integrated component.
            \end{proof}
            
            \begin{definition}
                For any real number $\mu \geqslant 0$ and any integer $k \neq 0$, we define
                \begin{equation}
                \label{eq:defAMuKNeq0}
                    \begin{array}{lll}
                        A_{\mu,k} \left( s \right) & = & \frac{1}{\pi} \sin \left( \pi s \right) \cdot \left( 4 \mu + 1 \right)^{-s} \left( \left \vert k \right \vert^{2\delta} - \frac{1}{4} \right)^{-s} H_{\mu,k} \left( 2 \left \vert k \right \vert^{\delta} \sqrt{\frac{1}{4} + \mu} \right)
                    \end{array}
                \end{equation}
                in the strip $1 < \Re s < 2$, where we also define
                \begin{equation}
                \label{eq:defBMuKNeq0}
                    \begin{array}{lll}
                        B_{\mu,k} \left( s \right) & = & \frac{2}{\pi} s \sin \left( \pi s \right) \int_{\sqrt{\frac{1}{4} + \mu}}^{2 \left \vert k \right \vert^{\delta} \sqrt{\frac{1}{4} + \mu}} t \left( t^2 - \left( \frac{1}{4} + \mu \right) \right)^{-s-1} H_{\mu,k} \left( t \right) \mathrm{d}t.
                    \end{array}
                \end{equation}
            \end{definition}
            
            Let us first study the series over $k \neq 0$ whose general term is given by $B_{\mu,k}$.
            
            \begin{proposition}
            \label{prop:studyBMuKNeq0}
                The function
                \begin{equation}
                \label{eq:functionSumKNeq0BMuK}
                    \begin{array}{lll}
                        s & \longmapsto & \sum\limits_{k \neq 0} B_{\mu,k} \left( s \right)
                    \end{array}
                \end{equation}
                is holomorphic on the strip $6 - 1/\delta < \Re s < 2$, which contains $0$ if we have $\delta < 1/6$, in which case its derivative at $s=0$ vanishes.
            \end{proposition}
            
            \begin{proof}
                Using the bound \eqref{eq:boundHmukKDelta} from proposition \ref{prop:HmukBound}, the function
                \begin{equation}
                    \begin{array}{lll}
                        s & \longmapsto & \sum\limits_{k \neq 0} \int_{\sqrt{\frac{1}{4} + \mu}}^{2 \left \vert k \right \vert^{\delta} \sqrt{\frac{1}{4} + \mu}} t \left( t^2 - \left( \frac{1}{4} + \mu \right) \right)^{-s-1} H_{\mu,k} \left( t \right) \mathrm{d}t
                    \end{array}
                \end{equation}
                is holomorphic on the strip $6 - 1/\delta < \Re s < 2$, and, assuming we have $\delta < 1/6$, the derivative at $s=0$ of the function \eqref{eq:functionSumKNeq0BMuK} vanishes, because of the factor $s \sin \left( \pi s \right)$.
            \end{proof}

    \subsection{\texorpdfstring{Study of the terms $A_{\mu,k}$}{Study of the terms Amuk}}
    \label{subsec:studyAMuK}
    
        In the previous subsection, we used the splittings of the interval of integration \eqref{eq:splittingIntervalIntegrationK0} and \eqref{eq:splittingIntervalIntegrationKNeq0} to get decompositions
        \begin{equation}
        \label{eq:decompLMukIntoAMuKBMuK}
            \begin{array}{lll}
                L_{\mu,k} \left( s \right) & = & A_{\mu,k} \left( s \right) + B_{\mu,k} \left( s \right)
            \end{array}
        \end{equation}
        on the strip $1 < \Re s < 2$. With propositions \ref{prop:studyBMuK0} and \ref{prop:studyBMuKNeq0}, we completed the study of the terms $B_{\mu,k}$, and we now move on to the terms $A_{\mu,k}$, which require more care. One of the main ingredients will be a consequence of proposition .

        \subsubsection{The case $k = 0$}
        \label{subsubsec:studyAMuK0}
        
            Let us begin with the study of the single term $A_{\mu,0}$, assuming we have $\alpha \neq 0$. We will first prove that the function $A_{\mu,0}$ is holomorphic near $0$, and find an asymptotic expansion as $a$ goes to infinity of the quantity
            \begin{equation}
            \label{eq:defFinitePartAPrimeMu0AsMuGoesTo0}
                \begin{array}{lll}
                    \Fp_{\mu \rightarrow 0^+} A_{\mu,0} ' \left( 0 \right) & = & \lim\limits_{\mu \rightarrow 0^+} \left[ \log \mu + A_{\mu,0}' \left( 0 \right) \right]
                \end{array}
            \end{equation}
            which will be proved to be finite. The notation $\Fp$ stands for ``finite part'', and corresponds to the constant term in an asymptotic expansion. 
            
            \begin{proposition}
            \label{prop:studyAMu0AsAGoesToInfinity}
                The function
                \begin{equation}
                    \begin{array}{lll}
                        s & \longmapsto & A_{\mu,0} \left( s \right)
                    \end{array}
                \end{equation}
                is entire, and we have
                \begin{equation}
                \label{eq:AMu0DerAt0}
                    \begin{array}{lll}
                        A_{\mu,0} ' \left( 0 \right) & =  & - \log \mu + \log \left \vert g_0 \left( i \right) \right \vert - \log \left \vert i g_0 ' \left( \frac{i}{2} \right) \right \vert \\[0.5em]
                        
                        && \qquad \qquad \qquad \qquad \qquad - \frac{3}{4} i \frac{\partial}{\partial t}_{\left \vert t = \frac{i}{2} \right.} \log \left \vert g_0 \left( t \right) \right \vert + o \left( 1 \right)
                    \end{array}
                \end{equation}
                as $\mu$ goes to $0^+$. The quantity \eqref{eq:defFinitePartAPrimeMu0AsMuGoesTo0} is therefore finite, and satisfies
                \begin{equation}
                \label{eq:asymptAMu0DerAt0AsAGoesToInfinity}
                    \begin{array}{lll}
                        \Fp_{\mu \rightarrow 0^+} A_{\mu,0} ' \left( 0 \right) & = & \log \left \vert g_0 \left( i \right) \right \vert - \frac{3}{4} i \frac{\partial}{\partial t}_{\left \vert t = \frac{i}{2} \right.} \log \left \vert g_0 \left( t \right) \right \vert \\[0.5em]
                        
                        && \qquad \qquad \qquad \qquad \qquad - \log a - \log \pi \alpha + o \left( 1 \right)
                    \end{array}
                \end{equation}
                as $a$ goes to infinity.
            \end{proposition}
            
            \begin{proof}
                In order to get \eqref{eq:AMu0DerAt0}, it is enough to note that we have $g_0 \left( i/2 \right) = 0$, which is a consequence of proposition \ref{prop:zerosFunctionF0}. The term $\log \mu$ serves to make a difference quotient appear, and $i g_0 ' \left( i/2 \right)$ is then a real number. Let us move on to the asymptotics \eqref{eq:asymptAMu0DerAt0AsAGoesToInfinity}. Using \cite[Prop. C.8]{dutour:cuspDirichlet}, where special values of the modified Bessel functions recalled, we have
                \begin{equation}
                    \begin{array}{lllll}
                        i g_0' \left( \frac{i}{2} \right) & = & - 4 \pi \alpha a \left( 2 \mathbb{E}_1 \left( 4 \pi \alpha a \right) e^{4 \pi \alpha a} - \frac{1}{4 \pi \alpha a} \right) & = & - \frac{1}{\pi \alpha a} + o \left( \frac{1}{a} \right),
                    \end{array}
                \end{equation}
                as $a$ goes to infinity, where $\mathbb{E}_1$ denotes the exponential integral, the asymptotics used here also being present in \cite[Prop. C.8]{dutour:cuspDirichlet}.
            \end{proof}
            
            \begin{remark}
                It is possible to fully compute the $a$-asymptotic expansion \eqref{eq:asymptAMu0DerAt0AsAGoesToInfinity}, but unnecessary, as the non-explicit terms will be canceled later by other contributions.
            \end{remark}
            
            Let us now study the derivative $A_{\mu,0}' \left( 0 \right)$ as $\mu$ goes to infinity. To that effect, we need to state a consequence of proposition \ref{prop:largeOrderAsymptFunctionGtxy}.
            
            \begin{lemma}
            \label{lem:asymptExpLargeOrderFunctionLogGK0}
                For any real parameter $t > 0$ large enough, we have
                \begin{equation}
                    \begin{array}{lll}
                        \log \left \vert g_0 \left( i t \right) \right \vert & = & \log 2 + \frac{1}{2} \log \left( 4 \pi^2 \alpha^2 a^2 + t^2 \right) - \frac{2 \pi \alpha a}{\sqrt{4 \pi^2 \alpha^2 a^2 + t^2}} \\[0.5em]
                        
                        && \qquad \qquad - \frac{1}{2} \cdot \frac{1}{\sqrt{4 \pi^2 \alpha^2 a^2 + t^2}} \cdot \frac{t^2}{4 \pi^2 \alpha^2 a^2 + t^2} + \eta_{2,3} \left( t \right)
                    \end{array}
                \end{equation}
                where the remainder $\eta_{2,3}$ satisfies
                \begin{equation}
                    \begin{array}{lll}
                        \left \vert \eta_{2,3} \left( t \right) \right \vert & \leqslant & \frac{C}{4 \pi^2 \alpha^2 a^2 + t^2},
                    \end{array}
                \end{equation}
                where $C > 0$ is a constant.
            \end{lemma}
            
            \begin{proof}
                In \eqref{eq:asymptLargeOrderFunctionGtxy}, one notes that taking $t$ large enough guarantees the positivity of the last factor. Taking the logarithm of this asymptotic expansion then gives this lemma, using the power series expansion of the logarithm, which converges for $t$ large enough, to control the remainder $\eta_{2,3}$.
            \end{proof}
            
            \begin{proposition}
            \label{prop:studyAMu0AsMuGoesToInfinity}
                We have
                \begin{equation}
                \label{eq:asymptLargeMuDerAMu0}
                    \begin{array}{lll}
                        A_{\mu,0}' \left( 0 \right) & = & \log \left \vert g_0 \left( 2 i \sqrt{\frac{1}{4} + \mu} \right) \right \vert - \frac{3}{2} i \sqrt{\frac{1}{4} + \mu} \; \frac{\partial}{\partial t}_{\left \vert t = i \sqrt{\frac{1}{4} + \mu} \right.} \log \left \vert g_0 \left( t \right) \right \vert\\[1em]
                        
                        && \qquad \qquad \qquad \qquad \qquad \qquad \qquad - \frac{1}{2} \log \mu + \log 2 + o \left( 1 \right)
                    \end{array}
                \end{equation}
                as $\mu$ goes to infinity.
            \end{proposition}
            
            \begin{proof}
                This is a direct consequence of lemma \ref{lem:asymptExpLargeOrderFunctionLogGK0}.
            \end{proof}
            
            \begin{remark}
                The right-hand side of \eqref{eq:asymptLargeMuDerAMu0} contains terms whose $\mu$-asymptotic behavior are not explicit. These contributions will be canceled in the final results by others appearing later.
            \end{remark}

        \subsubsection{The case $k \neq 0$}
        
            We now consider the terms $A_{\mu,k}$ for every $k \neq 0$, and their sum over such integers, for any real number $\mu \geqslant 0$. Let us begin by proving that this series has a holomorphic continuation near $0$, and by finding the asymptotic expansion of its derivative at $s=0$ as $a$ goes to infinity.
            
            \begin{remark}
            \label{rmk:abuseOfNotationSeriesContinuation}
                In this subsection, we will conflate the notations for a function defined by a series and its continuation. This abuse of notation is designed to keep track more easily of which term is being considered.
            \end{remark}
            
            \begin{proposition}
            \label{prop:globalStudyAMuKNeq0}
                For any real number $\mu \geqslant 0$, the function
                \begin{equation}
                    \begin{array}[t]{lll}
                        s & \longmapsto & \sum\limits_{k \neq 0} \; A_{\mu,k} \left( s \right)
                    \end{array}
                \end{equation}
                is holomorphic on the half-plane $\Re s > 3-1/\delta$, which contains $0$ if we have $\delta < 1/3$. The derivative at $s=0$ of this funtion further satisfies
                \begin{equation}
                    \begin{array}[t]{lll}
                        \frac{\partial}{\partial s}_{\left \vert s = 0 \right.} \sum\limits_{k \neq 0} \; A_{\mu,k} \left( s \right) & = & o \left( 1 \right)
                    \end{array}
                \end{equation}
                as $a$ goes to infinity.
            \end{proposition}
            
            \begin{proof}
                This is a direct consequence of estimate \eqref{eq:boundHmukKDelta} from proposition \ref{prop:HmukBound}.
            \end{proof}
            
            Having obtained the existence of the continuation and the $a$-asymptotic behavior of the derivative at $s=0$, we now turn our attention to the $\mu$-asymptotic expansion. This requires more care, and, as in \cite[Sec. 3.4.3]{dutour:cuspDirichlet}, the use of Taylor expansions and of the \textit{Ramanujan summation}. The reader is referred to \cite{candelpergher:ramanujan} for a detailed treatment of the latter, and to \cite[App. B]{dutour:cuspDirichlet} for an overview. We will require the following lemma, stated in a generality which will be useful later.
            
            \begin{lemma}
            \label{lem:logExpansionGkLargeOrder}
                For any real parameter $t > 0$, and any integer $k \neq 0$, we have
                \begin{equation}
                \label{eq:logExpansionGkLargeOrder}
                    \begin{array}{ll}
                        \multicolumn{2}{l}{\log \left \vert g_k \left( it \right) \right \vert} \\[1em]
                        
                        = \hspace{-2pt} & \scriptstyle \log \left( \sqrt{4 \pi^2 \left( k + \alpha \right)^2 a^2 + t^2} + \sqrt{4 \pi^2 \left( k - \alpha \right)^2 a^2 + t^2} \right) - \frac{4 \pi \alpha a}{\sqrt{4 \pi^2 \left( k + \alpha \right)^2 a^2 + t^2} + \sqrt{4 \pi^2 \left( k - \alpha \right)^2 a^2 + t^2}} \\[1em]
                        
                        & \quad \scriptstyle - \frac{1}{2} \cdot \frac{t^2}{\sqrt{4 \pi^2 \left( k + \alpha \right)^2 a^2 + t^2} + \sqrt{4 \pi^2 \left( k - \alpha \right)^2 a^2 + t^2}} \left( \frac{1}{4 \pi^2 \left( k + \alpha \right)^2 a^2 + t^2} + \frac{1}{4 \pi^2 \left( k - \alpha \right)^2 a^2 + t^2} \right) \\[1em]
                        
                        & \quad \quad \scriptstyle + \eta_{2,3} \left( t, a, k \right),
                    \end{array}
                \end{equation}
                with either $t \geqslant A$ or $a \geqslant B$ being large enough. The remainder $\eta_{2,3}$ satisfies
                \begin{equation}
                \label{eq:estimateRemainderAsymptLargeOrderGKNeq0}
                    \begin{array}{lll}
                        \left \vert \eta_{2,3} \left( t, a, k \right) \right \vert & \leqslant & \frac{C}{4 \pi^2 k^2 \left( 1- \alpha \right)^2 a^2 + t^2},
                    \end{array}
                \end{equation}
                for some constant $C > 0$ and $t \geqslant A$ large enough.
            \end{lemma}
            
            \begin{proof}
                After taking the logarithm of the absolute value of the right-hand side of expansion \eqref{eq:asymptLargeOrderFunctionGtxy}, we note that we have to study the term
                \begin{equation}
                \label{eq:logExpansionGkLargeOrderRemainingTerm}
                    \begin{array}{l}
                        \scriptstyle \log \left[ \; \left \vert 1 - \frac{4 \pi \alpha a}{\sqrt{4 \pi^2 \left( k + \alpha \right)^2 a^2 + t^2} + \sqrt{4 \pi^2 \left( k - \alpha \right)^2 a^2 + t^2}} \right. \right. \\[0.5em]
                        
                        \scriptstyle \qquad - \frac{1}{2} \cdot \frac{t^2}{\sqrt{4 \pi^2 \left( k + \alpha \right)^2 a^2 + t^2} + \sqrt{4 \pi^2 \left( k - \alpha \right)^2 a^2 + t^2}} \cdot \left( \frac{1}{4 \pi^2 \left( k + \alpha \right)^2 a^2 + t^2} + \frac{1}{4 \pi^2 \left( k - \alpha \right)^2 a^2 + t^2} \right) \\[0.5em]
                        
                        \scriptstyle \qquad \qquad + \eta_{2,2} \left( t, 2 \pi \left \vert k + \alpha \right \vert a, 2 \pi \left \vert k - \alpha \right \vert a \right) \Big \vert \; \Big].
                    \end{array}
                \end{equation}
                For any real number $t > 0$, and any integer $k \neq 0$, we have
                \begin{equation}
                    \begin{array}[t]{lllllll}
                        0 & < & \frac{4 \pi \alpha a}{\sqrt{4 \pi^2 \left( k + \alpha \right)^2 a^2 + t^2} + \sqrt{4 \pi^2 \left( k - \alpha \right)^2 a^2 + t^2}} & < & \frac{2 \alpha}{\left \vert k + \alpha \right \vert + \left \vert k - \alpha \right \vert} & \leqslant & \alpha,
                    \end{array}
                \end{equation}
                due to the fact that we have $\left \vert k \right \vert \geqslant 1$ and $0 \leqslant \alpha < 1$. Thus, for either $t$ large enough or $a$ large enough, we can remove the absolute value in \eqref{eq:logExpansionGkLargeOrderRemainingTerm} and use the power series expansion of the logarithm. We get the required expansion \eqref{eq:logExpansionGkLargeOrder}, with
                \begin{equation}
                \label{eq:remainderLargeOrderAsymptLogGKNeq0}
                    \begin{array}{lll}
                        \multicolumn{3}{l}{\scriptstyle \eta_{2,3} \left( t, a, k \right)} \\
                        
                        & = & \scriptstyle - \eta_{2,2} \left( t, 2 \pi \left \vert k + \alpha \right \vert a, 2 \pi \left \vert k - \alpha \right \vert a \right) - \sum\limits_{n=2}^{+ \infty} \frac{1}{n} \left[ \frac{4 \pi \alpha a}{\sqrt{4 \pi^2 \left( k + \alpha \right)^2 a^2 + t^2} + \sqrt{4 \pi^2 \left( k - \alpha \right)^2 a^2 + t^2}} \right. \\[0.5em]
                        
                        && \quad \scriptstyle + \frac{1}{2} \cdot \frac{t^2}{\sqrt{4 \pi^2 \left( k + \alpha \right)^2 a^2 + t^2} + \sqrt{4 \pi^2 \left( k - \alpha \right)^2 a^2 + t^2}} \cdot \left( \frac{1}{4 \pi^2 \left( k + \alpha \right)^2 a^2 + t^2} + \frac{1}{4 \pi^2 \left( k - \alpha \right)^2 a^2 + t^2} \right) \\[0.5em]
                        
                        && \qquad \qquad \qquad \qquad \qquad \qquad \qquad \qquad \quad \scriptstyle  - \eta_{2,2} \left( t, 2 \pi \left \vert k + \alpha \right \vert a, 2 \pi \left \vert k - \alpha \right \vert a \right) \Big]^n.
                    \end{array}
                \end{equation}
                The estimate \eqref{eq:estimateRemainderAsymptLargeOrderGKNeq0} on the remainder $\eta_{2,3}$ is then a consequence of formula \eqref{eq:remainderLargeOrderAsymptLogGKNeq0}.
            \end{proof}
            
            \begin{remark}
                Note that, although the expression \eqref{eq:remainderLargeOrderAsymptLogGKNeq0} holds for $a$ large enough, it does not provide an estimate of the type \eqref{eq:estimateRemainderAsymptLargeOrderGKNeq0}, because the first term in the sum over $n$ of \eqref{eq:remainderLargeOrderAsymptLogGKNeq0} does not vanish as $a$ goes to infinity.
            \end{remark}
            
            We will now break down the study of the terms $A_{\mu,k}$, and of their sum over all integers $k \neq 0$, as $\mu$ goes to infinity into different parts, according to the asymptotic expansion \eqref{eq:logExpansionGkLargeOrder} from lemma \eqref{lem:logExpansionGkLargeOrder}. Let us first recall that we have
            \begin{equation}
\begin{array}{lll}
                    A_{\mu,k} \left( s \right) & = & \frac{\sin \left( \pi s \right)}{\pi} \cdot \left( 4 \mu + 1 \right)^{-s} \left( \left \vert k \right \vert^{2\delta} - \frac{1}{4} \right)^{-s} H_{\mu,k} \left( 2 \left \vert k \right \vert^{\delta} \sqrt{\frac{1}{4} + \mu} \right)
                \end{array}
            \end{equation}
            on the strip $1 < \Re s < 2$, for any integer $k \neq 0$ and any real number $\mu \geqslant 0$, with
            \begin{equation}
            \label{eq:recallDefHMuK}
                \begin{array}{lll}
                    H_{\mu,k} \left( t \right) & = & \log \left \vert g_k \left( i t \right) \right \vert - \log \left \vert g_k \left( i\sqrt{\frac{1}{4} + \mu} \, \right) \right \vert \\[1em]
                        
                    && \qquad \qquad \qquad \qquad \quad - i \frac{t^2 - \left( 1/4 + \mu \right)}{\sqrt{4 \mu + 1}} \frac{\partial}{\partial t}_{\vert t = i\sqrt{\frac{1}{4} + \mu}} \log \left \vert g_k \left( t \right) \right \vert,
                \end{array}
            \end{equation}
            for any real parameter $t > 0$. The first two terms on the right-hand side of \eqref{eq:recallDefHMuK} will be studied separately, using \eqref{eq:logExpansionGkLargeOrder}.
            
            \subsubsection*{\textbf{First part}}
            \label{AMuKNeq0:part1}
            
                We begin with first term where the remainder $\eta_{2,3}$ appears.
                
                \begin{proposition}
                \label{prop:studyAMuKNeq0FirstPart}
                    The function
                    \begin{equation}
                        \begin{array}{lll}
                            s & \longmapsto & \frac{\sin \left( \pi s \right)}{\pi} \left( 4 \mu + 1 \right)^{-s} \sum\limits_{k \neq 0} \left( \left \vert k \right \vert^{2\delta} - \frac{1}{4} \right)^{-s} \eta_{2,3} \left( 2 \left \vert k \right \vert^{\delta} \sqrt{\frac{1}{4} + \mu}, a, k \right)
                        \end{array}
                    \end{equation}
                    is holomorphic around $0$, vanishes at $0$, and its derivative at $s=0$ satisfies
                    \begin{equation}
                        \begin{array}{lll}
                            \sum\limits_{k \neq 0} \eta_{2,3} \left( 2 \left \vert k \right \vert^{\delta} \sqrt{\frac{1}{4} + \mu}, a, k \right) & = & o \left( 1 \right)
                        \end{array}
                    \end{equation}
                    as $\mu$ goes to infinity.
                \end{proposition}
                
                \begin{proof}
                    This is a direct consequence of estimate \eqref{eq:estimateRemainderAsymptLargeOrderGKNeq0} from lemma \ref{lem:logExpansionGkLargeOrder}.
                \end{proof}

            \subsubsection*{\textbf{Second part}}
            \label{AMuKNeq0:part2}
            
                Let us now take care of the other term involving the remainder $\eta_{2,3}$.
                
                \begin{proposition}
                \label{prop:studyAMuKNeq0SecondPart}
                    The function
                    \begin{equation}
                        \begin{array}{lll}
                            s & \longmapsto & - \frac{\sin \left( \pi s \right)}{\pi} \left( 4 \mu + 1 \right)^{-s} \sum\limits_{k \neq 0} \left( \left \vert k \right \vert^{2\delta} - \frac{1}{4} \right)^{-s} \eta_{2,3} \left( \sqrt{\frac{1}{4} + \mu}, a, k \right)
                        \end{array}
                    \end{equation}
                    is holomorphic around $0$, vanishes at $0$, and its derivative at $s=0$ satisfies
                    \begin{equation}
                        \begin{array}{lll}
                            \sum\limits_{k \neq 0} \eta_{2,3} \left( \sqrt{\frac{1}{4} + \mu}, a, k \right) & = & o \left( 1 \right)
                        \end{array}
                    \end{equation}
                    as $\mu$ goes to infinity.
                \end{proposition}
                
                \begin{proof}
                    This is a direct consequence of estimate \eqref{eq:estimateRemainderAsymptLargeOrderGKNeq0} from lemma \ref{lem:logExpansionGkLargeOrder}.
                \end{proof}

            \subsubsection*{\textbf{Third part}}
            \label{AMuKNeq0:part3}
            
                We now move on to the first non-remainder term of \eqref{eq:logExpansionGkLargeOrder}.
                
                \begin{proposition}
                \label{prop:studyAMuKNeq0ThirdPart}
                    Assume that $1/\left( 2 \delta \right)$ is not an integer. The function
                    \begin{equation}
                    \label{eq:functionAMukThirdPart}
                        \begin{array}{lll}
                            s & \longmapsto & \scriptstyle \frac{\sin \left( \pi s \right)}{\pi} \left( 4 \mu + 1 \right)^{-s} \sum\limits_{k \neq 0} \left( \left( \left \vert k \right \vert^{2\delta} - \frac{1}{4} \right)^{-s} \log \left( \sqrt{4 \pi^2 \left( k + \alpha \right)^2 a^2 + \left( 4 \mu + 1 \right) \left \vert k \right \vert^{2 \delta}} \right. \right. \\[0.5em]
                            
                            && \qquad \qquad \qquad \qquad \qquad \qquad \qquad \scriptstyle  \left. \left. + \sqrt{4 \pi^2 \left( k - \alpha \right)^2 a^2 + \left( 4 \mu + 1 \right) \left \vert k \right \vert^{2 \delta}} \right) \right),
                        \end{array}
                    \end{equation}
                    which is holomorphic on the half-plane $\Re s > 1/\left( 2 \delta \right)$, has a holomorphic continuation to an open region of $\mathbb{C}$ containing $0$, which vanishes at $s=0$, and whose derivative at $s=0$ equals
                    \begin{equation}
                    \label{eq:studyAMukThirdPartDerAt0}
                        \begin{array}{l}
                            \scriptstyle\frac{\partial}{\partial s}_{\vert s = 0} \left[ \frac{\sin \left( \pi s \right)}{\pi} \left( 4 \mu + 1 \right)^{-s} \sum\limits_{k \neq 0} \left( \left \vert k \right \vert^{- 2 \delta s} \log \left( \sqrt{4 \pi^2 \left( k + \alpha \right)^2 a^2 + \left( 4 \mu + 1 \right) \left \vert k \right \vert^{2 \delta}} \right. \right. \right. \\[0.5em]
                            
                            \qquad \qquad \qquad \qquad \qquad \qquad \qquad \qquad \qquad \scriptstyle \left. \left. \left. + \sqrt{4 \pi^2 \left( k - \alpha \right)^2 a^2 + \left( 4 \mu + 1 \right) \left \vert k \right \vert^{2 \delta}} \right) \right) \right],
                        \end{array}
                    \end{equation}
                    the abuse of notation being explained in remark \ref{rmk:abuseOfNotationSeriesContinuation}.
                \end{proposition}
                
                \begin{proof}
                    We first note that, using the symmetry $k \leftrightarrow - k$, we need only study
                    \begin{equation}
                    \label{eq:studyAmuKNeq0ThirdPart01}
                        \begin{array}{lll}
                            s & \longmapsto & \scriptstyle 2 \frac{\sin \left( \pi s \right)}{\pi} \left( 4 \mu + 1 \right)^{-s} \sum\limits_{k = 1}^{+ \infty} \left( \left( k ^{2\delta} - \frac{1}{4} \right)^{-s} \log \left( \sqrt{4 \pi^2 \left( k + \alpha \right)^2 a^2 + \left( 4 \mu + 1 \right) k^{2 \delta}} \right. \right. \\[0.5em]
                            
                            && \qquad \qquad \qquad \qquad \qquad \qquad \qquad \scriptstyle  \left. \left. + \sqrt{4 \pi^2 \left( k - \alpha \right)^2 a^2 + \left( 4 \mu + 1 \right)  k^{2 \delta}} \right) \right).
                        \end{array}
                    \end{equation}
                    Using the binomial formula, which is recalled in \cite[Prop. C.26]{dutour:cuspDirichlet}, we have
                    \begin{equation}
                    \label{eq:studyAmuKNeq0ThirdPart02}
                        \begin{array}{ll}
                            \multicolumn{2}{l}{\scriptstyle \sum\limits_{k = 1}^{+ \infty} \left( k ^{2\delta} - \frac{1}{4} \right)^{-s} \log \left( \sqrt{4 \pi^2 \left( k + \alpha \right)^2 a^2 + \left( 4 \mu + 1 \right) k^{2 \delta}} + \sqrt{4 \pi^2 \left( k - \alpha \right)^2 a^2 + \left( 4 \mu + 1 \right)  k^{2 \delta}} \right)} \\[0.5em]
                            
                            = \hspace{-6pt} & \scriptstyle \sum\limits_{j=0}^{+ \infty} \frac{\left( s \right)_j}{j! 4^j} \sum\limits_{k = 1}^{+ \infty} k^{- 2 \delta \left( s + j \right)} \log \left( \sqrt{4 \pi^2 \left( k + \alpha \right)^2 a^2 + \left( 4 \mu + 1 \right) k^{2 \delta}} + \sqrt{4 \pi^2 \left( k - \alpha \right)^2 a^2 + \left( 4 \mu + 1 \right)  k^{2 \delta}} \right)
                        \end{array}
                    \end{equation}
                    on the half-plane $\Re s > 1/\left( 2 \delta \right)$. We now note that we have
                    \begin{equation}
                    \label{eq:studyAmuKNeq0ThirdPart03}
                        \begin{array}{lll}
                            \multicolumn{3}{l}{\scriptstyle \log \left( \sqrt{4 \pi^2 \left( k + \alpha \right)^2 a^2 + \left( 4 \mu + 1 \right) k^{2 \delta}} + \sqrt{4 \pi^2 \left( k - \alpha \right)^2 a^2 + \left( 4 \mu + 1 \right)  k^{2 \delta}} \right)} \\[0.5em]
                            
                            \qquad & = & \scriptstyle \frac{1}{2} \log \left( 4 \pi^2 \left( k + \alpha \right)^2 a^2 + \left( 4 \mu + 1 \right) k^{2 \delta} \right) + \log \left( 1 + \sqrt{1 - \frac{16 \pi^2 \alpha k a^2}{4 \pi^2 \left( k + \alpha \right)^2 a^2 + \left( 4 \mu + 1 \right) k^{2 \delta}}} \right).
                        \end{array}
                    \end{equation}
                    Note that the parameter inside the last logarithm of \eqref{eq:studyAmuKNeq0ThirdPart03} is in $\left] 0, 1 \right[$. Let us start with the term induced by the first part of the right-hand side of \eqref{eq:studyAmuKNeq0ThirdPart03}, given by
                    \begin{equation}
                    \label{eq:studyAmuKNeq0ThirdPart04}
                        \begin{array}{c}
                            \scriptstyle \frac{\sin \left( \pi s \right)}{\pi} \left( 4 \mu + 1 \right)^{-s} \sum\limits_{j=0}^{+ \infty} \frac{\left( s \right)_j}{j! 4^j} \sum\limits_{k = 1}^{+ \infty} k^{- 2 \delta \left( s + j \right)} \log \left( 4 \pi^2 \left( k + \alpha \right)^2 a^2 + \left( 4 \mu + 1 \right) k^{2 \delta} \right).
                        \end{array}
                    \end{equation}
                    For any integer $k \geqslant 1$, we have
                    \begin{equation}
                    \label{eq:studyAmuKNeq0ThirdPart05}
                        \begin{array}{lllll}
                            \multicolumn{3}{l}{\scriptstyle \log \left( 4 \pi^2 \left( k + \alpha \right)^2 a^2 + \left( 4 \mu + 1 \right) k^{2 \delta} \right)} & = & \scriptstyle 2 \log \left( 2 \pi a \right) + 2 \log k + \log \left( \left( 1 + \frac{\alpha}{k} \right)^2 + \frac{4 \mu + 1}{4 \pi^2 a^2} \cdot k^{2 \delta - 2} \right) \\[0.5em]
                            
                            \qquad \qquad & = & \multicolumn{3}{l}{\scriptstyle 2 \log \left( 2 \pi a \right) + 2 \log k + \frac{2 \alpha}{k} + \left( \log \left( 1 + \frac{2 \alpha}{k} + \frac{\alpha^2}{k^2} + \frac{4 \mu + 1}{4 \pi^2 a^2} \cdot k^{2 \delta - 2} \right) - \frac{2 \alpha}{k} \right).}
                        \end{array}
                    \end{equation}
                    Using the fundamental theorem of calculus, we now get
                    \begin{equation}
                    \label{eq:studyAmuKNeq0ThirdPart06}
                        \begin{array}{lllll}
                            \multicolumn{3}{l}{\scriptstyle \log \left( 1 + \frac{2 \alpha}{k} + \frac{\alpha^2}{k^2} + \frac{4 \mu + 1}{4 \pi^2 a^2} \cdot k^{2 \delta - 2} \right) - \frac{2 \alpha}{k}} & = & \scriptstyle \int_0^{1/k} \left( \frac{2 \alpha + 2 \alpha^2 t + \frac{4 \mu + 1}{2 \pi^2 a^2} \left( 1 - \delta \right) t^{1-2 \delta}}{1 + 2 \alpha t + \alpha^2 t^2 + \frac{4 \mu + 1}{4 \pi^2 a^2} t^{2 - 2 \delta}} - 2 \alpha \right) \mathrm{d}t \\[1em]
                            
                            \qquad & = & \multicolumn{3}{l}{\scriptstyle \int_0^{1/k} \frac{1}{1 + 2 \alpha t + \alpha^2 t^2 + \frac{4 \mu + 1}{4 \pi^2 a^2} t^{2 - 2 \delta}} \left( \cancel{2 \alpha} + 2 \alpha^2 t + \frac{4 \mu + 1}{2 \pi^2 a^2} \left( 1 - \delta \right) t^{1 - 2 \delta} - \cancel{2 \alpha} \right.} \\[0.5em] 
                            
                            && \multicolumn{3}{l}{\qquad \qquad \qquad \qquad \qquad \qquad \qquad \scriptstyle \left. - 4 \alpha^2 t - 2 \alpha^3 t^2 - \frac{4 \mu + 1}{2 \pi^2 a^2} \alpha t^{2 - 2 \delta} \right) \mathrm{d}t,}
                        \end{array}
                    \end{equation}
                    which gives
                    \begin{equation}
                    \label{eq:studyAmuKNeq0ThirdPart07}
                        \begin{array}{lll}
                            \scriptstyle \left \vert \log \left( 1 + \frac{2 \alpha}{k} + \frac{\alpha^2}{k^2} + \frac{4 \mu + 1}{4 \pi^2 a^2} \cdot k^{2 \delta - 2} \right) - \frac{2 \alpha}{k} \right \vert & \hspace{-2pt} \leqslant \hspace{-2pt} & \scriptstyle \frac{\alpha^2}{k^2} + \frac{4 \mu + 1}{4 \pi^2 a^2} \cdot \frac{1}{k^{2 - 2 \delta}} + \frac{2 \alpha^3}{3 k^3} + \frac{4 \mu + 1}{2 \pi^2 a^2} \cdot \frac{\alpha}{3 - 2 \delta} \cdot \frac{1}{k^{3 - 2 \delta}}.
                        \end{array}
                    \end{equation}
                    Therefore, the function
                    \begin{equation}
                    \label{eq:studyAmuKNeq0ThirdPart08}
                        \begin{array}{lll}
                            s & \longmapsto & \scriptstyle \frac{\sin \left( \pi s \right)}{\pi} \left( 4 \mu + 1 \right)^{-s} \sum\limits_{j=1}^{+ \infty} \frac{\left( s \right)_j}{j! 4^j} \sum\limits_{k = 1}^{+ \infty} k^{- 2 \delta \left( s + j \right)} \left( \log \left( 1 + \frac{2 \alpha}{k} + \frac{\alpha^2}{k^2} + \frac{4 \mu + 1}{4 \pi^2 a^2} \cdot k^{2 \delta - 2} \right) - \frac{2 \alpha}{k} \right)
                        \end{array}
                    \end{equation}
                    is holomorphic around $0$, both it and its derivative vanish at $s=0$ because of $\left( s \right)_j$ for $j \geqslant 1$. The function corresponding to $j=0$, given by
                    \begin{equation}
                    \label{eq:studyAmuKNeq0ThirdPart09}
                        \begin{array}{lll}
                            s & \longmapsto & \scriptstyle \frac{\sin \left( \pi s \right)}{\pi} \left( 4 \mu + 1 \right)^{-s} \sum\limits_{k = 1}^{+ \infty} k^{- 2 \delta s} \left( \log \left( 1 + \frac{2 \alpha}{k} + \frac{\alpha^2}{k^2} + \frac{4 \mu + 1}{4 \pi^2 a^2} \cdot k^{2 \delta - 2} \right) - \frac{2 \alpha}{k} \right)
                        \end{array}
                    \end{equation}
                    is also holomorphic around $0$ and vanishes at $s=0$. Its derivative at $s=0$ is part of \eqref{eq:studyAMukThirdPartDerAt0}. Going back to decomposition \eqref{eq:studyAmuKNeq0ThirdPart05}, we now note that the function
                    \begin{equation}
                    \label{eq:studyAmuKNeq0ThirdPart10}
                        \begin{array}{lll}
                            s & \longmapsto & \scriptstyle \frac{\sin \left( \pi s \right)}{\pi} \left( 4 \mu + 1 \right)^{-s} \sum\limits_{j=1}^{+ \infty} \frac{\left( s \right)_j}{j! 4^j} \sum\limits_{k = 1}^{+ \infty} k^{- 2 \delta \left( s + j \right)} \left( 2 \log \left( 2 \pi a \right) + 2 \log k + \frac{2 \alpha}{k} \right)
                        \end{array}
                    \end{equation}
                    has a holomorphic continuation near $0$, since the series above converges absolutely for $j$ large enough, and the finitely many remaining terms can be written in terms of the Riemann zeta function. The function \eqref{eq:studyAmuKNeq0ThirdPart10} and its derivative at $s=0$ vanish because of $\left( s \right)_j$ for $j \geqslant 1$. The remaining term, corresponding to $j=0$, is given by
                    \begin{equation}
                    \label{eq:studyAmuKNeq0ThirdPart11}
                        \begin{array}{lll}
                            \multicolumn{3}{l}{\scriptstyle \frac{\sin \left( \pi s \right)}{\pi} \left( 4 \mu + 1 \right)^{-s} \sum\limits_{k = 1}^{+ \infty} k^{- 2 \delta s} \left( 2 \log \left( 2 \pi a \right) + 2 \log k + \frac{2 \alpha}{k} \right)} \\[0.5em]
                            
                            \qquad \qquad \qquad & = & \scriptstyle \frac{\sin \left( \pi s \right)}{\pi} \left( 4 \mu + 1 \right)^{-s} \left( 2 \log \left( 2 \pi a \right) \zeta \left( 2 \delta s \right) - 2 \zeta' \left( 2 \delta s \right) + 2 \alpha \zeta \left( 2 \delta s + 1 \right) \right).
                        \end{array}
                    \end{equation}
                    It also induces a holomorphic function near $0$, as the factor $\sin \left( \pi s \right)$ cancels the simple pole of $\zeta \left( 2 \delta s + 1 \right)$. Its derivative at $s=0$ is part of \eqref{eq:studyAMukThirdPartDerAt0}, and its value at~$0$ will be canceled later. Looking back at \eqref{eq:studyAmuKNeq0ThirdPart03}, we must now deal with
                    \begin{equation}
                    \label{eq:studyAmuKNeq0ThirdPart12}
                        \begin{array}{c}
                            \scriptstyle 2 \frac{\sin \left( \pi s \right)}{\pi} \left( 4 \mu + 1 \right)^{-s} \sum\limits_{j=0}^{+ \infty} \frac{\left( s \right)_j}{j! 4^j} \sum\limits_{k = 1}^{+ \infty} k^{- 2 \delta \left( s + j \right)} \log \left( 1 + \sqrt{1 - \frac{16 \pi^2 \alpha k a^2}{4 \pi^2 \left( k + \alpha \right)^2 a^2 + \left( 4 \mu + 1 \right) k^{2 \delta}}} \; \right).
                        \end{array}
                    \end{equation}
                    Applying Taylor's formula at order $2$, we get
                    \begin{equation}
                    \label{eq:studyAmuKNeq0ThirdPart13}
                        \begin{array}{lll}
                            \multicolumn{3}{l}{\scriptscriptstyle \log \left( 1 + \sqrt{1 - \frac{16 \pi^2 \alpha k a^2}{4 \pi^2 \left( k + \alpha \right)^2 a^2 + \left( 4 \mu + 1 \right) k^{2 \delta}}} \; \right)} \\[0.5em]
                            
                            \quad & = & \scriptstyle \log 2 - \frac{1}{4} \cdot \frac{16 \pi^2 \alpha k a^2}{4 \pi^2 \left( k + \alpha \right)^2 a^2 + \left( 4 \mu + 1 \right) k^{2 \delta}} - \frac{1}{4} \int_0^{\frac{16 \pi^2 \alpha k a^2}{4 \pi^2 \left( k + \alpha \right)^2 a^2 + \left( 4 \mu + 1 \right) k^{2 \delta}}} \left[ \left( 1 - x \right)^{-3/2} \right. \\[0.5em]
                            
                            && \qquad \qquad \qquad \qquad \qquad \quad \scriptstyle  \left. \cdot \frac{1 + 2 \sqrt{1-x}}{\left( 1 + \sqrt{1-x} \right)^2} \left( \frac{16 \pi^2 \alpha k a^2}{4 \pi^2 \left( k + \alpha \right)^2 a^2 + \left( 4 \mu + 1 \right) k^{2 \delta}} - x \right) \right] \mathrm{d}x.
                        \end{array}
                    \end{equation}
                    We now note that we have, for any real number $t \in \left] 0, 1 \right[$,
                    \begin{equation}
                    \label{eq:studyAmuKNeq0ThirdPart14}
                        \begin{array}{lll}
                            \scriptstyle \left \vert \int_0^t \left( 1 - x \right)^{-3/2} \cdot \frac{1 + 2 \sqrt{1-x}}{\left( 1 + \sqrt{1-x} \right)^2} \left( t - x \right) \mathrm{d}x \right \vert & \leqslant & \scriptstyle t^2 \left( 1 - t \right)^{-3/2}
                        \end{array}
                    \end{equation}
                    which gives, for any integer $k \geqslant 1$,
                    \begin{equation}
                    \label{eq:studyAmuKNeq0ThirdPart15}
                        \begin{array}{lll}
                            \multicolumn{3}{l}{\scriptstyle \left \vert \int_0^{\frac{16 \pi^2 \alpha k a^2}{4 \pi^2 \left( k + \alpha \right)^2 a^2 + \left( 4 \mu + 1 \right) k^{2 \delta}}} \left( 1 - x \right)^{-3/2} \cdot \frac{1 + 2 \sqrt{1-x}}{\left( 1 + \sqrt{1-x} \right)^2} \left( \frac{16 \pi^2 \alpha k a^2}{4 \pi^2 \left( k + \alpha \right)^2 a^2 + \left( 4 \mu + 1 \right) k^{2 \delta}} - x \right) \mathrm{d}x \right \vert} \\[1em]
                            
                            \qquad \qquad \qquad \quad & \leqslant & \scriptstyle \frac{256 \pi^4 \alpha^2 k^2 a^4}{\sqrt{\left( 4 \mu + 1 \right) k^{2 \delta} + 4 \pi^2 \left( k + \alpha \right)^2 a^2}} \cdot \left( \left( 4 \mu + 1 \right) k^{2 \delta} + 4 \pi^2 \left( k - \alpha \right)^2 a^2 \right)^{-3/2}.
                        \end{array}
                    \end{equation}
                    Therefore, the function
                    \begin{equation}
                    \label{eq:studyAmuKNeq0ThirdPart16}
                        \begin{array}{lll}
                            s & \hspace{-2pt} \mapsto \hspace{-2pt} & \scriptstyle - \frac{1}{2} \frac{\sin \left( \pi s \right)}{\pi} \left( 4 \mu + 1 \right)^{-s} \sum\limits_{j=1}^{+ \infty} \frac{\left( s \right)_j}{j! 4^j} \sum\limits_{k = 1}^{+ \infty} k^{- 2 \delta \left( s + j \right)} \int_0^{\frac{16 \pi^2 \alpha k a^2}{4 \pi^2 \left( k + \alpha \right)^2 a^2 + \left( 4 \mu + 1 \right) k^{2 \delta}}}  \left[ \left( 1 - x \right)^{-3/2} \right. \\[0.5em]
                            
                            && \qquad \qquad \qquad \qquad \qquad \qquad \scriptstyle \left. \cdot \frac{1 + 2 \sqrt{1-x}}{\left( 1 + \sqrt{1-x} \right)^2} \left( \frac{16 \pi^2 \alpha k a^2}{4 \pi^2 \left( k + \alpha \right)^2 a^2 + \left( 4 \mu + 1 \right) k^{2 \delta}} - x \right) \mathrm{d}x \right]
                        \end{array}
                    \end{equation}
                    is holomorphic around $0$. Both it and its derivative vanish at $s=0$, because of the Pochhammer symbol $\left( s \right)_j$ for $j \geqslant 1$. The function corresponding to $j=0$, given by
                    \begin{equation}
                    \label{eq:studyAmuKNeq0ThirdPart17}
                        \begin{array}{lll}
                            s & \longmapsto & \scriptstyle - \frac{1}{2} \frac{\sin \left( \pi s \right)}{\pi} \left( 4 \mu + 1 \right)^{-s} \sum\limits_{k = 1}^{+ \infty} k^{- 2 \delta s} \int_0^{\frac{16 \pi^2 \alpha k a^2}{4 \pi^2 \left( k + \alpha \right)^2 a^2 + \left( 4 \mu + 1 \right) k^{2 \delta}}}  \left[ \left( 1 - x \right)^{-3/2} \right. \\[0.5em]
                            
                            && \qquad \qquad \qquad \qquad \qquad \scriptstyle \left. \cdot \frac{1 + 2 \sqrt{1-x}}{\left( 1 + \sqrt{1-x} \right)^2} \left( \frac{16 \pi^2 \alpha k a^2}{4 \pi^2 \left( k + \alpha \right)^2 a^2 + \left( 4 \mu + 1 \right) k^{2 \delta}} - x \right) \mathrm{d}x \right]
                        \end{array}
                    \end{equation}
                    is also holomorphic around $0$, vanishes at $s=0$, and its derivative at $s=0$ is part of \eqref{eq:studyAMukThirdPartDerAt0}. Going back to \eqref{eq:studyAmuKNeq0ThirdPart13}, we must now study the term
                    \begin{equation}
                    \label{eq:studyAmuKNeq0ThirdPart18}
                        \begin{array}{c}
                            \scriptstyle - \frac{1}{2} \frac{\sin \left( \pi s \right)}{\pi} \left( 4 \mu + 1 \right)^{-s} \sum\limits_{j=0}^{+ \infty} \frac{\left( s \right)_j}{j! 4^j} \sum\limits_{k = 1}^{+ \infty} k^{- 2 \delta \left( s + j \right)} \frac{16 \pi^2 \alpha k a^2}{4 \pi^2 \left( k + \alpha \right)^2 a^2 + \left( 4 \mu + 1 \right) k^{2 \delta}}.
                        \end{array}
                    \end{equation}
                    Considering the sum over $j \geqslant 1$, the function
                    \begin{equation}
                    \label{eq:studyAmuKNeq0ThirdPart19}
                        \begin{array}{lll}
                            s & \longmapsto & \scriptstyle - \frac{1}{2} \frac{\sin \left( \pi s \right)}{\pi} \left( 4 \mu + 1 \right)^{-s} \sum\limits_{j=1}^{+ \infty} \frac{\left( s \right)_j}{j! 4^j} \sum\limits_{k = 1}^{+ \infty} k^{- 2 \delta \left( s + j \right)} \frac{16 \pi^2 \alpha k a^2}{4 \pi^2 \left( k + \alpha \right)^2 a^2 + \left( 4 \mu + 1 \right) k^{2 \delta}}
                        \end{array}
                    \end{equation}
                    is holomorphic near $0$. Both it and its derivative vanish at $s=0$ because of $\left( s \right)_j$. The term corresponding to $j=0$ is given by
                    \begin{equation}
                    \label{eq:studyAmuKNeq0ThirdPart20}
                        \begin{array}{c}
                            \scriptstyle - \frac{1}{2} \frac{\sin \left( \pi s \right)}{\pi} \left( 4 \mu + 1 \right)^{-s} \sum\limits_{k = 1}^{+ \infty} k^{- 2 \delta s} \frac{16 \pi^2 \alpha k a^2}{4 \pi^2 \left( k + \alpha \right)^2 a^2 + \left( 4 \mu + 1 \right) k^{2 \delta}}.
                        \end{array}
                    \end{equation}
                    For any integer $k \geqslant 1$, we now have
                    \begin{equation}
                    \label{eq:studyAmuKNeq0ThirdPart21}
                        \begin{array}{lll}
                            \scriptstyle \frac{16 \pi^2 \alpha k a^2}{4 \pi^2 \left( k + \alpha \right)^2 a^2 + \left( 4 \mu + 1 \right) k^{2 \delta}} & = & \scriptstyle \frac{4 \alpha k}{\left( k + \alpha \right)^2} - \frac{\alpha \left( 4 \mu + 1 \right) k^{1 + 2 \delta}}{\pi^2 \left( k + \alpha \right)^4 a^2} \cdot \frac{4 \pi^2 \left( k + \alpha \right)^2 a^2}{4 \pi^2 \left( k + \alpha \right)^2 a^2 + \left( 4 \mu + 1 \right) k^{2 \delta}}.
                        \end{array}
                    \end{equation}
                    Induced by the second part of the right-hand side of \eqref{eq:studyAmuKNeq0ThirdPart21}, the function
                    \begin{equation}
                    \label{eq:studyAmuKNeq0ThirdPart22}
                        \begin{array}{lll}
                            s & \mapsto & \scriptstyle 2 \alpha \frac{\sin \left( \pi s \right)}{\pi} \left( 4 \mu + 1 \right)^{-s+1} \sum\limits_{k = 1}^{+ \infty} k^{1 - 2 \delta \left( s - 1 \right)} \frac{1}{\left( k + \alpha \right)^2} \cdot \frac{1}{4 \pi^2 \left( k + \alpha \right)^2 a^2 + \left( 4 \mu + 1 \right) k^{2 \delta}}
                        \end{array}
                    \end{equation}
                    is holomorphic around $0$, vanishes at $0$, and its derivative at $s=0$ is part of \eqref{eq:studyAMukThirdPartDerAt0}. The function induced by the first part of the right-hand side of \eqref{eq:studyAmuKNeq0ThirdPart21} is given by
                    \begin{equation}
                    \label{eq:studyAmuKNeq0ThirdPart24}
                        \begin{array}{lll}
                            s & \longmapsto & \scriptstyle - 2 \alpha \frac{\sin \left( \pi s \right)}{\pi} \left( 4 \mu + 1 \right)^{-s} \sum\limits_{k = 1}^{+ \infty} k^{1 - 2 \delta s} \frac{1}{\left( k + \alpha \right)^2}.
                        \end{array}
                    \end{equation}
                    We can break this function apart by writing, for any integer $k \geqslant 1$,
                    \begin{equation}
                    \label{eq:studyAmuKNeq0ThirdPart25}
                        \begin{array}{lll}
                            \frac{1}{\left( k + \alpha \right)^2} & = & \frac{1}{k^2} - \frac{\alpha}{k} \cdot \frac{1}{\left( k + \alpha \right)^2} \left( 2 + \frac{\alpha}{k} \right).
                        \end{array}
                    \end{equation}
                    Induced by the second part of the right-hand side of \eqref{eq:studyAmuKNeq0ThirdPart25}, the function
                    \begin{equation}
                    \label{eq:studyAmuKNeq0ThirdPart26}
                        \begin{array}{lll}
                            s & \longmapsto & \scriptstyle 2 \alpha^2 \frac{\sin \left( \pi s \right)}{\pi} \left( 4 \mu + 1 \right)^{-s} \sum\limits_{k = 1}^{+ \infty} k^{- 2 \delta s} \frac{1}{\left( k + \alpha \right)^2} \left( 2 + \frac{\alpha}{k} \right)
                        \end{array}
                    \end{equation}
                    is holomorphic around $0$, vanishes at $0$, and its derivative at $s=0$ is part of \eqref{eq:studyAMukThirdPartDerAt0}. The remaining term from \eqref{eq:studyAmuKNeq0ThirdPart25} induces
                    \begin{equation}
                    \label{eq:studyAmuKNeq0ThirdPart27}
                        \begin{array}{lll}
                            \scriptstyle - 2 \alpha \frac{\sin \left( \pi s \right)}{\pi} \left( 4 \mu + 1 \right)^{-s} \sum\limits_{k = 1}^{+ \infty} k^{-1 - 2 \delta s} & = & \scriptstyle - 2 \alpha \frac{\sin \left( \pi s \right)}{\pi} \left( 4 \mu + 1 \right)^{-s} \zeta \left( 1 + 2 \delta s \right)
                        \end{array}
                    \end{equation}
                    and the associated function of $s$ has a holomorphic continuation near $0$, whose derivative at $s=0$ is part of \eqref{eq:studyAMukThirdPartDerAt0}. Note that the sum of \eqref{eq:studyAmuKNeq0ThirdPart11} and \eqref{eq:studyAmuKNeq0ThirdPart27} vanishes at $s=0$. We now study the function induced by the last part of \eqref{eq:studyAmuKNeq0ThirdPart13}, which is given by
                    \begin{equation}
                    \label{eq:studyAmuKNeq0ThirdPart28}
                        \begin{array}{lll}
                            s & \longmapsto & \scriptstyle 2 \log 2 \frac{\sin \left( \pi s \right)}{\pi} \left( 4 \mu + 1 \right)^{-s} \sum\limits_{j=0}^{+ \infty} \frac{\left( s \right)_j}{j! 4^j} \zeta \left( 2 \delta \left( s + j \right) \right).
                        \end{array}
                    \end{equation}
                    Under the assumption that $1/\left( 2 \delta \right)$ is not an integer, so that we never have $2 \delta j = 1$, the function
                    \begin{equation}
                    \label{eq:studyAmuKNeq0ThirdPart29}
                        \begin{array}{lll}
                            s & \longmapsto & \scriptstyle 2 \log 2 \frac{\sin \left( \pi s \right)}{\pi} \left( 4 \mu + 1 \right)^{-s} \sum\limits_{j=1}^{+ \infty} \frac{\left( s \right)_j}{j! 4^j} \zeta \left( 2 \delta \left( s + j \right) \right)
                        \end{array}
                    \end{equation}
                    has a holomorphic continuation near $0$. Both it and its derivative vanish at $s=0$. Similarly, the function corresponding to $j=0$, given by
                    \begin{equation}
                    \label{eq:studyAmuKNeq0ThirdPart30}
                        \begin{array}{lll}
                            s & \longmapsto & \scriptstyle 2 \log 2 \frac{\sin \left( \pi s \right)}{\pi} \left( 4 \mu + 1 \right)^{-s} \zeta \left( 2 \delta s \right)
                        \end{array}
                    \end{equation}
                    has a holomorphic continuation near $0$, whose derivative at $s=0$ is part of \eqref{eq:studyAMukThirdPartDerAt0}. The proof of the proposition is thus complete.
                \end{proof}

            \subsubsection*{\textbf{Fourth part}}
            \label{AMuKNeq0:part4}
            
                We now move on to the second term of the right-hand side of \eqref{eq:logExpansionGkLargeOrder}. Once again, we use the notation $\Fp$, which stands for ``finite part'', here in the context of Laurent expansions at a point.
                
                \begin{proposition}
                \label{prop:studyAMuKNeq0FourthPart}
                    Assume that $1/\left( 2 \delta \right)$ is not an integer. The function
                    \begin{equation}
                        \begin{array}{lll}
                            s & \longmapsto & \scriptstyle - \frac{\sin \left( \pi s \right)}{\pi} \left( 4 \mu + 1 \right)^{-s} \sum\limits_{k \neq 0} \left[ \left( \left \vert k \right \vert^{2\delta} - \frac{1}{4} \right)^{-s} \right. \\
                            
                            && \qquad \qquad \qquad \left. \scriptstyle \cdot \frac{4 \pi \alpha a}{\sqrt{4 \pi^2 \left( k + \alpha \right)^2 a^2 + \left( 4 \mu + 1 \right) \left \vert k \right \vert^{2 \delta}} + \sqrt{4 \pi^2 \left( k - \alpha \right)^2 a^2 + \left( 4 \mu + 1 \right) \left \vert k \right \vert^{2 \delta}}} \right],
                        \end{array}
                    \end{equation}
                    which is holomorphic on the half-plane $\Re s > 0$, has a holomorphic continuation to an open region of $\mathbb{C}$ containing $0$, whose value at $s=0$ equals
                    \begin{equation}
                    \label{eq:studyAMukFourthPartValueAt0}
                        \begin{array}{l}
                            \Fp_{s=0} \left[ - 2 \alpha \frac{\sin \left( \pi s \right)}{\pi} \left( 4 \mu + 1 \right)^{-s} \zeta \left( 2 \delta s + 1 \right) \right]
                        \end{array}
                    \end{equation}
                    and whose derivative at $s=0$ equals
                    \begin{equation}
                    \label{eq:studyAMukFourthPartDerAt0}
                        \begin{array}{l}
                            \scriptstyle \frac{\partial}{\partial s}_{\vert s = 0} \left[ - \frac{\sin \left( \pi s \right)}{\pi} \left( 4 \mu + 1 \right)^{-s} \sum\limits_{k \neq 0} \left \vert k \right \vert^{-2 \delta s} \right. \\
                            
                            \qquad \qquad \qquad \qquad \qquad \left. \scriptstyle \cdot \frac{4 \pi \alpha a}{\sqrt{4 \pi^2 \left( k + \alpha \right)^2 a^2 + \left( 4 \mu + 1 \right) \left \vert k \right \vert^{2 \delta}} + \sqrt{4 \pi^2 \left( k - \alpha \right)^2 a^2 + \left( 4 \mu + 1 \right) \left \vert k \right \vert^{2 \delta}}} \right]
                        \end{array}
                    \end{equation}
                \end{proposition}
                
                \begin{proof}
                    Using the binomial formula, which is recalled in \cite[Prop. C.26]{dutour:cuspDirichlet}, we have
                    \begin{equation}
                    \label{eq:studyAmuKNeq0FourthPart01}
                        \begin{array}{lll}
                            \multicolumn{3}{l}{\scriptstyle \sum\limits_{k \neq 0} \left( \left \vert k \right \vert^{2\delta} - \frac{1}{4} \right)^{-s} \cdot \frac{4 \pi \alpha a}{\sqrt{4 \pi^2 \left( k + \alpha \right)^2 a^2 + \left( 4 \mu + 1 \right) \left \vert k \right \vert^{2 \delta}} + \sqrt{4 \pi^2 \left( k - \alpha \right)^2 a^2 + \left( 4 \mu + 1 \right) \left \vert k \right \vert^{2 \delta}}}} \\
                            
                            & = & \scriptstyle \sum\limits_{j=0}^{+ \infty} \frac{\left( s \right)_j}{j! 4^j} \sum\limits_{k \neq 0} \left \vert k \right \vert^{- 2 \delta \left( s + j \right)} \cdot \frac{4 \pi \alpha a}{\sqrt{4 \pi^2 \left( k + \alpha \right)^2 a^2 + \left( 4 \mu + 1 \right) \left \vert k \right \vert^{2 \delta}} + \sqrt{4 \pi^2 \left( k - \alpha \right)^2 a^2 + \left( 4 \mu + 1 \right) \left \vert k \right \vert^{2 \delta}}}
                        \end{array}
                    \end{equation}
                    Considering the sum over $j \geqslant 1$, we note that the function
                    \begin{equation}
                    \label{eq:studyAmuKNeq0FourthPart02}
                        \begin{array}{lll}
                            s & \longmapsto & \scriptstyle - \frac{\sin \left( \pi s \right)}{\pi} \left( 4 \mu + 1 \right)^{-s} \sum\limits_{j=1}^{+ \infty} \frac{\left( s \right)_j}{j! 4^j} \sum\limits_{k \neq 0} \left[ \left \vert k \right \vert^{- 2 \delta \left( s + j \right)} \right. \\
                            
                            && \qquad \qquad \qquad \quad \left. \scriptstyle \cdot \frac{4 \pi \alpha a}{\sqrt{4 \pi^2 \left( k + \alpha \right)^2 a^2 + \left( 4 \mu + 1 \right) \left \vert k \right \vert^{2 \delta}} + \sqrt{4 \pi^2 \left( k - \alpha \right)^2 a^2 + \left( 4 \mu + 1 \right) \left \vert k \right \vert^{2 \delta}}} \right]
                        \end{array}
                    \end{equation}
                    is holomorphic around $0$, vanishes at $0$, and that its derivative at $s=0$ also vanishes, because of the Pochhammer symbol $\left( s \right)_j$ for $j \geqslant 1$. Thus, we only have to study the term corresponding to $j=0$, which is given by
                    \begin{equation}
                    \label{eq:studyAmuKNeq0FourthPart03}
                        \begin{array}{l}
                            \scriptstyle - \frac{\sin \left( \pi s \right)}{\pi} \left( 4 \mu + 1 \right)^{-s} \sum\limits_{k \neq 0} \left \vert k \right \vert^{- 2 \delta s} \cdot \frac{4 \pi \alpha a}{\sqrt{4 \pi^2 \left( k + \alpha \right)^2 a^2 + \left( 4 \mu + 1 \right) \left \vert k \right \vert^{2 \delta}} + \sqrt{4 \pi^2 \left( k - \alpha \right)^2 a^2 + \left( 4 \mu + 1 \right) \left \vert k \right \vert^{2 \delta}}}.
                        \end{array}
                    \end{equation}
                    Using the symmetry $k \leftrightarrow -k$, we need only look at the function
                    \begin{equation}
                    \label{eq:studyAmuKNeq0FourthPart04}
                        \begin{array}{lll}
                            s & \longmapsto & \scriptstyle - 8 \pi \alpha a \frac{\sin \left( \pi s \right)}{\pi} \left( 4 \mu + 1 \right)^{-s} \sum\limits_{k = 1}^{+ \infty} \left[ k^{- 2 \delta s} \right. \\
                            
                            && \qquad \qquad \qquad \left. \scriptstyle \cdot \frac{4 \pi \alpha a}{\sqrt{4 \pi^2 \left( k + \alpha \right)^2 a^2 + \left( 4 \mu + 1 \right) k^{2 \delta}} + \sqrt{4 \pi^2 \left( k - \alpha \right)^2 a^2 + \left( 4 \mu + 1 \right) k^{2 \delta}}} \right].
                        \end{array}
                    \end{equation}
                    Furthermore, we can assume that we have $\alpha \neq 0$, or the function we study vanishes identically. For any integer $k \geqslant 1$, we have
                    \begin{equation}
                    \label{eq:studyAmuKNeq0FourthPart05TrickConjugateRoot}
                        \begin{array}{lll}
                            \multicolumn{3}{l}{\scriptstyle \frac{1}{\sqrt{4 \pi^2 \left( k + \alpha \right)^2 a^2 + \left( 4 \mu + 1 \right) k^{2 \delta}} + \sqrt{4 \pi^2 \left( k - \alpha \right)^2 a^2 + \left( 4 \mu + 1 \right)  k^{2 \delta}}}} \\[1em]
                            
                            \qquad \qquad & = & \scriptstyle \frac{1}{16 \pi^2 \alpha k a^2} \left( \sqrt{4 \pi^2 \left( k + \alpha \right)^2 a^2 + \left( 4 \mu + 1 \right) k^{2 \delta}} - \sqrt{4 \pi^2 \left( k - \alpha \right)^2 a^2 + \left( 4 \mu + 1 \right) k^{2 \delta}} \right)
                        \end{array}
                    \end{equation}
                    The term induced by the first part of the right-hand side of \eqref{eq:studyAmuKNeq0FourthPart02} is given by
                    \begin{equation}
                    \label{eq:studyAmuKNeq0FourthPart06}
                        \begin{array}{lll}
                            s & \longmapsto & \scriptstyle - \frac{1}{2 \pi a} \frac{\sin \left( \pi s \right)}{\pi} \left( 4 \mu + 1 \right)^{-s} \sum\limits_{k = 1}^{+ \infty} k^{- 1 - 2 \delta s} \sqrt{4 \pi^2 \left( k + \alpha \right)^2 a^2 + \left( 4 \mu + 1 \right) k^{2 \delta}}.
                        \end{array}
                    \end{equation}
                    For any integer $k \geqslant 1$, we have
                    \begin{equation}
                    \label{eq:studyAmuKNeq0FourthPart07}
                        \begin{array}{lll}
                            \scriptstyle \sqrt{4 \pi^2 \left( k + \alpha \right)^2 a^2 + \left( 4 \mu + 1 \right) k^{2 \delta}} & = & \scriptstyle 2 \pi \left( k + \alpha \right) a \sqrt{1 + \frac{\left( 4 \mu + 1 \right) k^{2 \delta}}{4 \pi^2 \left( k + \alpha \right)^2 a^2}} \\[0.5em]
                            
                            & = & \scriptstyle 2 \pi \left( k + \alpha \right) a + \frac{\left( 4 \mu + 1 \right) k^{2 \delta}}{4 \pi^2 \left( k + \alpha \right)^2 a^2 + \sqrt{4 \pi^2 \left( k + \alpha \right)^2 a^2 + \left( 4 \mu + 1 \right) k^{2 \delta}}}.
                        \end{array}
                    \end{equation}
                    The function induced by the second part of the right-hand side of \eqref{eq:studyAmuKNeq0FourthPart07}, given by
                    \begin{equation}
                    \label{eq:studyAmuKNeq0FourthPart08}
                        \begin{array}{lll}
                            s & \hspace{-3pt} \mapsto \hspace{-3pt} & \scriptstyle - \frac{1}{2 \pi a} \frac{\sin \left( \pi s \right)}{\pi} \left( 4 \mu + 1 \right)^{-s+1} \sum\limits_{k = 1}^{+ \infty} k^{- 1 - 2 \delta \left( s-1 \right)} \frac{1}{4 \pi^2 \left( k + \alpha \right)^2 a^2 + \sqrt{4 \pi^2 \left( k + \alpha \right)^2 a^2 + \left( 4 \mu + 1 \right) k^{2 \delta}}},
                        \end{array}
                    \end{equation}
                    is holomorphic near $0$ and vanishes at $0$. Its derivative at $s=0$ is part of \eqref{eq:studyAMukFourthPartDerAt0}. Let us move on to the term induced by the first part of the right-hand side of \eqref{eq:studyAmuKNeq0FourthPart07}, which is given by
                    \begin{equation}
                    \label{eq:studyAmuKNeq0FourthPart09}
                        \begin{array}{lll}
                            \scriptstyle - \frac{\sin \left( \pi s \right)}{\pi} \left( 4 \mu + 1 \right)^{-s} \sum\limits_{k = 1}^{+ \infty} k^{- 1 - 2 \delta s} \cdot \left( k + \alpha \right) & = & \scriptstyle - \frac{\sin \left( \pi s \right)}{\pi} \left( 4 \mu + 1 \right)^{-s} \left( \zeta \left( 2 \delta s \right) + \alpha \zeta \left( 1 + 2 \delta s \right) \right).
                        \end{array}
                    \end{equation}
                    The term \eqref{eq:studyAmuKNeq0FourthPart09} induces a holomorphic function around $0$, whose value at $s=0$ is part of \eqref{eq:studyAMukFourthPartValueAt0}, and whose derivative at $s=0$ is part of \eqref{eq:studyAMukFourthPartDerAt0}. The term induced by the second part of the right-hand side of \eqref{eq:studyAmuKNeq0FourthPart02} is given by
                    \begin{equation}
                    \label{eq:studyAmuKNeq0FourthPart10}
                        \begin{array}{lll}
                            s & \longmapsto & \scriptstyle \frac{1}{2 \pi a} \frac{\sin \left( \pi s \right)}{\pi} \left( 4 \mu + 1 \right)^{-s} \sum\limits_{k = 1}^{+ \infty} k^{- 1 - 2 \delta s} \sqrt{4 \pi^2 \left( k - \alpha \right)^2 a^2 + \left( 4 \mu + 1 \right) k^{2 \delta}}.
                        \end{array}
                    \end{equation}
                    Similarly to the study performed above, this function has a holomorphic continuation near $0$, whose value at $s=0$ is given by
                    \begin{equation}
                    \label{eq:studyAmuKNeq0FourthPart11}
                        \begin{array}{l}
                            \scriptstyle \Fp_{s=0} \left[ \frac{\sin \left( \pi s \right)}{\pi} \left( 4 \mu + 1 \right)^{-s} \left( \zeta \left( 2 \delta s \right) - \alpha \zeta \left( 1 + 2 \delta s \right) \right) \right]
                        \end{array}
                    \end{equation}
                    and whose derivative at $s=0$ is part of \eqref{eq:studyAMukFourthPartDerAt0}. This concludes the proof of the proposition, as we note that summing the value at $0$ of \eqref{eq:studyAmuKNeq0FourthPart09} with \eqref{eq:studyAmuKNeq0FourthPart11} gives the contribution \eqref{eq:studyAMukFourthPartValueAt0}.
                \end{proof}

            \subsubsection*{\textbf{Fifth part}}
            \label{AMuKNeq0:part5}
            
                We now turn our attention to the last non-remainder term appearing in the right-hand side of \eqref{eq:logExpansionGkLargeOrder}, in lemma \ref{lem:logExpansionGkLargeOrder}.
                
                \begin{proposition}
                \label{prop:studyAMuKNeq0FifthPart}
                    The function
                    \begin{equation}
                        \begin{array}{lll}
                            s & \mapsto & \scriptstyle \frac{1}{2} \frac{\sin \left( \pi s \right)}{\pi} \left( 4 \mu + 1 \right)^{-s+1} \sum\limits_{k \neq 0} \left[ \left( \left \vert k \right \vert^{2\delta} - \frac{1}{4} \right)^{-s} \left \vert k \right \vert^{2\delta} \right. \\[0.5em]
                            
                            && \scriptstyle \qquad \qquad \qquad \cdot \frac{1}{\sqrt{4 \pi^2 \left( k + \alpha \right)^2 a^2 + \left( 4 \mu + 1 \right) \left \vert k \right \vert^{2 \delta}} + \sqrt{4 \pi^2 \left( k - \alpha \right)^2 a^2 + \left( 4 \mu + 1 \right) \left \vert k \right \vert^{2 \delta}}} \\[0.5em]
                            
                            && \scriptstyle \qquad \qquad \qquad \quad \left. \cdot \left( \frac{1}{4 \pi^2 \left( k + \alpha \right)^2 a^2 + \left( 4 \mu + 1 \right) \left \vert k \right \vert^{2 \delta}} + \frac{1}{4 \pi^2 \left( k - \alpha \right)^2 a^2 + \left( 4 \mu + 1 \right) \left \vert k \right \vert^{2 \delta}} \right) \right]
                        \end{array}
                    \end{equation}
                    is holomorphic around $0$ and vanishes at $0$.
                \end{proposition}
                
                \begin{proof}
                    This is a direct study of the series involved.
                \end{proof}

            \subsubsection*{\textbf{Sixth part}}
            \label{AMuKNeq0:part6}
                
                We must now study the first in a series of more complicated terms, for which we will need to use the Ramanujan summation process. It is induced by the logarithm on the right-hand side of \eqref{eq:logExpansionGkLargeOrder}, in lemma \ref{lem:logExpansionGkLargeOrder}.
                
                \begin{proposition}
                \label{prop:studyAMuKNeq0SixthPart}
                    Assume that $1 / \left( 2 \delta \right)$ is not an integer. The function
                    \begin{equation}
                        \label{eq:functionAMukSixthPart}
                        \begin{array}{lll}
                            s & \longmapsto & \scriptstyle - \frac{\sin \left( \pi s \right)}{\pi} \left( 4 \mu + 1 \right)^{-s} \sum\limits_{k \neq 0} \left( \left( \left \vert k \right \vert^{2\delta} - \frac{1}{4} \right)^{-s} \log \left( \sqrt{4 \pi^2 \left( k + \alpha \right)^2 a^2 + \frac{1}{4} + \mu} \right. \right. \\[0.5em]
                            
                            && \qquad \qquad \qquad \qquad \qquad \qquad \qquad \qquad \scriptstyle  \left. \left. + \sqrt{4 \pi^2 \left( k - \alpha \right)^2 a^2 + \frac{1}{4} + \mu} \right) \right),
                        \end{array}
                    \end{equation}
                    which is holomorphic on the half-plane $\Re s > 1/ \left( 2 \delta \right)$, has a holomorphic continuation near $0$, which vanishes at $0$, and whose derivative at $s=0$ satisfies
                    \begin{equation}
                    \label{eq:studyAMukSixthPartDerAt0}
                        \begin{array}{c}
                            \frac{1}{2} \log \mu - \frac{1}{2 a} \sqrt{\mu} + \log 2 + o \left( 1 \right)
                        \end{array}
                    \end{equation}
                    as $\mu$ goes to infinity.
                \end{proposition}
                
                \begin{proof}
                    Let us begin by noting that, if $\alpha$ vanishes, the result is \cite[Prop. 3.49]{dutour:cuspDirichlet}, up to a constant and a direct computation. Therefore, we assume that we have $\alpha \neq 0$. Now, we note that the symmetry $k \leftrightarrow -k$ means that we only need to study
                    \begin{equation}
                    \label{eq:studyAmuKNeq0SixthPart01}
                        \begin{array}{l}
                            \scriptstyle - 2 \frac{\sin \left( \pi s \right)}{\pi} \left( 4 \mu + 1 \right)^{-s} \sum\limits_{k = 1}^{+ \infty} \left( k^{2\delta} - \frac{1}{4} \right)^{-s} \log \left( \sqrt{4 \pi^2 \left( k + \alpha \right)^2 a^2 + \frac{1}{4} + \mu} + \sqrt{4 \pi^2 \left( k - \alpha \right)^2 a^2 + \frac{1}{4} + \mu} \right).
                        \end{array}
                    \end{equation}
                    We can then use the binomial formula, which yields
                    \begin{equation}
                    \label{eq:studyAmuKNeq0SixthPart02}
                        \begin{array}{lll}
                            \multicolumn{3}{l}{\scriptstyle - 2 \frac{\sin \left( \pi s \right)}{\pi} \left( 4 \mu + 1 \right)^{-s} \sum\limits_{k = 1}^{+ \infty} \left( k^{2\delta} - \frac{1}{4} \right)^{-s} \log \left( \sqrt{4 \pi^2 \left( k + \alpha \right)^2 a^2 + \frac{1}{4} + \mu} + \sqrt{4 \pi^2 \left( k - \alpha \right)^2 a^2 + \frac{1}{4} + \mu} \right)} \\[0.5em]
                            
                            & = & \scriptstyle - 2 \frac{\sin \left( \pi s \right)}{\pi} \left( 4 \mu + 1 \right)^{-s} \sum\limits_{j=0}^{+ \infty} \frac{\left( s \right)_j}{j! 4^j} \sum\limits_{k = 1}^{+ \infty} k^{-2 \delta \left( s + j \right)} \log \left( \sqrt{4 \pi^2 \left( k + \alpha \right)^2 a^2 + \frac{1}{4} + \mu} \right. \\
                            
                            && \scriptstyle \qquad \qquad \qquad \qquad \qquad \qquad \qquad \qquad \qquad \quad \left. + \sqrt{4 \pi^2 \left( k - \alpha \right)^2 a^2 + \frac{1}{4} + \mu} \right)
                        \end{array}
                    \end{equation}
                    In order to use the Ramanujan summation, let us define the function
                    \begin{equation}
                    \label{eq:studyAmuKNeq0SixthPart03}
                        \begin{array}{lllll}
                            f_{s,j} & : & z & \longmapsto & \scriptstyle z^{-2 \delta \left( s + j \right)} \log \left( \sqrt{4 \pi^2 \left( z + \alpha \right)^2 a^2 + \frac{1}{4} + \mu} + \sqrt{4 \pi^2 \left( z - \alpha \right)^2 a^2 + \frac{1}{4} + \mu} \right)
                        \end{array}
                    \end{equation}
                    which is holomorphic on the half-plane $\Re z > \alpha$, where $s$ is a complex number such that we have $\Re s > 1/\left( 2 \delta \right)$. Note that we consider the principal branch of the complex logarithm and square root. Each function $f_{s,j}$ is of moderate growth, as defined in \cite[Def. B.3]{dutour:cuspDirichlet} and in \cite[Def. 1.3.1]{candelpergher:ramanujan}. Let us now verify the two hypotheses required for the Ramanujan summation to behave well, recalled in \cite[Thm. B.8]{dutour:cuspDirichlet} and in \cite[Sec. 1.4.3]{candelpergher:ramanujan}. The first one is that we have
                    \begin{equation}
                    \label{eq:studyAmuKNeq0SixthPart04}
\begin{array}{lll}
                            \lim\limits_{k \rightarrow + \infty} f_{s,j} \left( k \right) & = & 0
                        \end{array}
                    \end{equation}
                    which holds, since we have in particular $\Re s > 0$. Let us now prove that we have
                    \begin{equation}
                    \label{eq:studyAmuKNeq0SixthPart05}
                        \begin{array}{lll}
                            \lim\limits_{k \rightarrow + \infty} \int_0^{+ \infty} \frac{f_{s,j} \left( k + i t \right) - f_{s,j} \left( k - it \right)}{e^{2 \pi t} - 1} \mathrm{d}t & = & 0.
                        \end{array}
                    \end{equation}
                    Writing the difference $f_{s,j} \left( k + i t \right) - f_{s,j} \left( k - it \right)$ as an integral using the fundamental theorem of calculus, and finding estimates of the derivative $f_{s,j}'$, we can apply the dominated convergence theorem to prove \eqref{eq:studyAmuKNeq0SixthPart05}. We thus have
                    \begin{equation}
                    \label{eq:studyAmuKNeq0SixthPart06}
                        \begin{array}{lll}
                            \sum\limits_{k=1}^{+ \infty} f_{s,j} \left( k \right) & = & \sum\limits_{k \geqslant 1}^{\left( \mathcal{R} \right)} f_{s,j} \left( k \right) + \int_1^{+ \infty} f_{s,j} \left( x \right) \mathrm{d}x
                        \end{array}
                    \end{equation}
                    Using \cite[Thm. B.11]{dutour:cuspDirichlet}, the function
                    \begin{equation}
                    \label{eq:studyAmuKNeq0SixthPart07}
                        \begin{array}{lll}
                            s & \longmapsto & - 2 \frac{\sin \left( \pi s \right)}{\pi} \left( 4 \mu + 1 \right)^{-s} \sum\limits_{j=0}^{+ \infty} \frac{\left( s \right)_j}{j! 4^j} \sum\limits_{k \geqslant 1}^{\left( \mathcal{R} \right)} f_{s,j} \left( k \right)
                        \end{array}
                    \end{equation}
                    is entire, and its derivative at $s=0$ is given by
                    \begin{equation}
                    \label{eq:studyAmuKNeq0SixthPart08}
                        \begin{array}{l}
                            - f_{0,0} \left( 1 \right) - 2 i \int_0^{+\infty} \frac{f_{0,0} \left( 1 + i t \right) - f_{0,0} \left( 1 - it \right)}{e^{2 \pi t} - 1} \mathrm{d}t
                        \end{array}
                    \end{equation}
                    It should be noted that the functions $f_{s,j}$ implicitely depend on $\mu$. Let us find the asymptotic expansion of \eqref{eq:studyAmuKNeq0SixthPart08} as $\mu$ goes to infinity. We have
                    \begin{equation}
                    \label{eq:studyAmuKNeq0SixthPart09}
                        \begin{array}{lll}
                            - f_{0,0} \left( 1 \right) & = & \scriptstyle - \log \left( \sqrt{4 \pi^2 \left( 1 + \alpha \right)^2 a^2 + \frac{1}{4} + \mu} + \sqrt{4 \pi^2 \left( 1 - \alpha \right)^2 a^2 + \frac{1}{4} + \mu} \right) \\[0.5em]
                            
                            & = & - \frac{1}{2} \log \mu - \log 2 + o \left( 1 \right)
                        \end{array}
                    \end{equation}
                    as $\mu$ goes to infinity. Let us now move on to the integral constituting the second part of \eqref{eq:studyAmuKNeq0SixthPart08}. Using once again the fundamental theorem of calculus, and estimates on the derivative $f_{0,0}'$, we can apply the dominated convergence theorem, to get
                    \begin{equation}
                    \label{eq:studyAmuKNeq0SixthPart10}
                        \begin{array}{lll}
                            \lim\limits_{\mu \rightarrow + \infty} \int_0^{+\infty} \frac{f_{0,0} \left( 1 + i t \right) - f_{0,0} \left( 1 - it \right)}{e^{2 \pi t} - 1} \mathrm{d}t & = & 0,
                        \end{array}
                    \end{equation}
                    thus giving the $\mu$-asymptotic expansion of \eqref{eq:studyAmuKNeq0SixthPart08}. We now turn our attention to the function induced by the second half of \eqref{eq:studyAmuKNeq0SixthPart06}, which is given by
                    \begin{equation}
                    \label{eq:studyAmuKNeq0SixthPart11}
                        \begin{array}{lll}
                            s & \longmapsto & \scriptstyle - 2 \frac{\sin \left( \pi s \right)}{\pi} \left( 4 \mu + 1 \right)^{-s} \sum\limits_{j=0}^{+ \infty} \frac{\left( s \right)_j}{j! 4^j} \int_1^{+\infty} x^{-2\delta \left( s + j \right)} \log \left( \sqrt{4 \pi^2 \left( x + \alpha \right)^2 a^2 + \frac{1}{4} + \mu} \right. \\
                            
                            && \qquad \qquad \qquad \qquad \qquad \qquad \qquad \qquad \scriptstyle \left. + \sqrt{4 \pi^2 \left( x - \alpha \right)^2 a^2 + \frac{1}{4} + \mu} \right) \mathrm{d}x.
                        \end{array}
                    \end{equation}
                    We will compute each of the integrals appearing in \eqref{eq:studyAmuKNeq0SixthPart11} using hypergeometric functions. The first step is to perform an integration by parts. We have
                    \begin{equation}
                    \label{eq:studyAmuKNeq0SixthPart12}
                        \begin{array}{ll}
                            \multicolumn{2}{l}{\scriptstyle \int_1^{+\infty} x^{-2\delta \left( s + j \right)} \log \left( \sqrt{4 \pi^2 \left( x + \alpha \right)^2 a^2 + \frac{1}{4} + \mu} + \sqrt{4 \pi^2 \left( x - \alpha \right)^2 a^2 + \frac{1}{4} + \mu} \right) \mathrm{d}x} \\[1em]
                            
                            = \hspace{-4pt} & \scriptstyle \left[ - \frac{1}{2 \delta \left( s + j \right) - 1} \cdot x^{-2 \delta \left( s + j \right) + 1} \log \left( \sqrt{4 \pi^2 \left( x + \alpha \right)^2 a^2 + \frac{1}{4} + \mu} + \sqrt{4 \pi^2 \left( x - \alpha \right)^2 a^2 + \frac{1}{4} + \mu} \right) \mathrm{d}x \right]_1^{+ \infty} \hspace{-6pt} \\[1em]
                            
                            & \scriptstyle + \frac{4 \pi^2 a^2}{2 \delta \left( s + j \right) - 1} \int_1^{+ \infty} x^{-2 \delta \left( s + j \right) + 1} \cdot \frac{1}{\sqrt{4 \pi^2 \left( x + \alpha \right)^2 a^2 + \frac{1}{4} + \mu} + \sqrt{4 \pi^2 \left( x - \alpha \right)^2 a^2 + \frac{1}{4} + \mu}} \\[1em]
                            
                            & \qquad \qquad \qquad \qquad \qquad \quad \scriptstyle \cdot \left( \frac{x + \alpha}{\sqrt{4 \pi^2 \left( x + \alpha \right)^2 a^2 + \frac{1}{4} + \mu}} + \frac{x - \alpha}{\sqrt{4 \pi^2 \left( x - \alpha \right)^2 a^2 + \frac{1}{4} + \mu}} \right) \mathrm{d}x.
                        \end{array}
                    \end{equation}
                    Having $\alpha \neq 0$, we have, for any real number $x \geqslant 1$,
                    \begin{equation}
                    \label{eq:studyAmuKNeq0SixthPart13}
                        \begin{array}{lll}
                            \multicolumn{3}{l}{\scriptstyle \frac{1}{\sqrt{4 \pi^2 \left( x + \alpha \right)^2 a^2 + \frac{1}{4} + \mu} + \sqrt{4 \pi^2 \left( x - \alpha \right)^2 a^2 + \frac{1}{4} + \mu}}} \\[1em]
                            
                            \qquad \qquad \qquad & = & \scriptstyle \frac{1}{16 \pi^2 \alpha x a^2} \left( \sqrt{4 \pi^2 \left( x + \alpha \right)^2 a^2 + \frac{1}{4} + \mu} - \sqrt{4 \pi^2 \left( x - \alpha \right)^2 a^2 + \frac{1}{4} + \mu} \right)
                        \end{array}
                    \end{equation}
                    and we can use this to further the computation started in \eqref{eq:studyAmuKNeq0SixthPart12}. We have
                    \begin{equation}
                    \label{eq:studyAmuKNeq0SixthPart14}
                        \begin{array}{lll}
                            \multicolumn{3}{l}{\scriptstyle \int_1^{+\infty} x^{-2\delta \left( s + j \right)} \log \left( \sqrt{4 \pi^2 \left( x + \alpha \right)^2 a^2 + \frac{1}{4} + \mu} + \sqrt{4 \pi^2 \left( x - \alpha \right)^2 a^2 + \frac{1}{4} + \mu} \right) \mathrm{d}x} \\[1em]
                            
                            & = & \scriptstyle \frac{1}{2 \delta \left( s + j \right) - 1} \log \left( \sqrt{4 \pi^2 \left( 1 + \alpha \right)^2 a^2 + \frac{1}{4} + \mu} + \sqrt{4 \pi^2 \left( 1 - \alpha \right)^2 a^2 + \frac{1}{4} + \mu} \right) \\[1em]
                            
                            && \scriptstyle + \frac{1}{2 \delta \left( s + j \right) - 1} \cdot \frac{1}{4 \alpha} \int_1^{+ \infty} x^{-2 \delta \left( s + j \right)} \left( 2 \alpha + \left( x - \alpha \right) \sqrt{\frac{4 \pi^2 \left( x + \alpha \right)^2 a^2 + 1/4 + \mu}{4 \pi^2 \left( x - \alpha \right)^2 a^2 + 1/4 + \mu}} \right. \\[1em]
                            
                            && \qquad \qquad \qquad \qquad \qquad \qquad \qquad \scriptstyle \quad\left. - \left( x + \alpha \right) \sqrt{\frac{4 \pi^2 \left( x - \alpha \right)^2 a^2 + 1/4 + \mu}{4 \pi^2 \left( x + \alpha \right)^2 a^2 + 1/4 + \mu}} \right) \mathrm{d}x \\[1em]
                            
                            & = & \scriptstyle \frac{1}{2 \delta \left( s + j \right) - 1} \log \left( \sqrt{4 \pi^2 \left( 1 + \alpha \right)^2 a^2 + \frac{1}{4} + \mu} + \sqrt{4 \pi^2 \left( 1 - \alpha \right)^2 a^2 + \frac{1}{4} + \mu} \right) + \frac{1}{2} \cdot \frac{1}{\left( 2 \delta \left( s + j \right) - 1 \right)^2} \\[1em]
                            
                            && \scriptstyle \qquad + \frac{1}{4 \alpha} \cdot \frac{1}{2 \delta \left( s + j \right) - 1} \int_1^{+ \infty} x^{- 2 \delta \left( s + j \right) + 1} \sqrt{\frac{4 \pi^2 \left( x + \alpha \right)^2 a^2 + 1/4 + \mu}{4 \pi^2 \left( x - \alpha \right)^2 a^2 + 1/4 + \mu}} \mathrm{d}x \\
\end{array}
                    \end{equation}
                    We will deal with each of the term appearing on the right-hand side of \eqref{eq:studyAmuKNeq0SixthPart14} separately. First, since we never have $2 \delta j = 1$, we note that the function
                    \begin{equation}
                    \label{eq:studyAmuKNeq0SixthPart15}
                        \begin{array}{lll}
                            s & \longmapsto & \scriptstyle - 2 \frac{\sin \left( \pi s \right)}{\pi} \left( 4 \mu + 1 \right)^{-s} \sum\limits_{j=0}^{+ \infty} \frac{\left( s \right)_j}{j! 4^j} \frac{1}{2 \delta \left( s + j \right) - 1} \log \left( \sqrt{4 \pi^2 \left( 1 + \alpha \right)^2 a^2 + \frac{1}{4} + \mu} \right. \\
                            
                            && \qquad \qquad \qquad \qquad \qquad \qquad \qquad \qquad \scriptstyle \left. + \sqrt{4 \pi^2 \left( 1 - \alpha \right)^2 a^2 + \frac{1}{4} + \mu} \right)
                        \end{array}
                    \end{equation}
                    is holomorphic around $0$, vanishes at $0$, and its derivative at $s=0$ is given by
                    \begin{equation}
                    \label{eq:studyAmuKNeq0SixthPart16}
                        \begin{array}{lll}
                            \scriptstyle -2 \log \left( \sqrt{4 \pi^2 \left( 1 + \alpha \right)^2 a^2 + \frac{1}{4} + \mu} + \sqrt{4 \pi^2 \left( 1 - \alpha \right)^2 a^2 + \frac{1}{4} + \mu} \right) & = & \scriptstyle \log \mu + 2 \log 2 + o \left( 1 \right)
                        \end{array}
                    \end{equation}
                    as $\mu$ goes to infinity. The function induced by the second part of the right-hand side of \eqref{eq:studyAmuKNeq0SixthPart14}, given by
                    \begin{equation}
                    \label{eq:studyAmuKNeq0SixthPart17}
                        \begin{array}{lll}
                            s & \longmapsto & \scriptstyle - \frac{\sin \left( \pi s \right)}{\pi} \left( 4 \mu + 1 \right)^{-s} \sum\limits_{j=0}^{+ \infty} \frac{\left( s \right)_j}{j! 4^j} \cdot \frac{1}{\left( 2 \delta \left( s + j \right) - 1 \right)^2},
                        \end{array}
                    \end{equation}
                    is also holomorphic around $0$, vanishes at $0$, and its derivative at $s=0$ satisfies
                    \begin{equation}
                    \label{eq:studyAmuKNeq0SixthPart18}
                        \begin{array}{lll}
                            \scriptstyle \frac{\partial}{\partial s}_{\vert s=0} \left[ - \frac{\sin \left( \pi s \right)}{\pi} \left( 4 \mu + 1 \right)^{-s} \sum\limits_{j=0}^{+ \infty} \frac{\left( s \right)_j}{j! 4^j} \cdot \frac{1}{\left( 2 \delta \left( s + j \right) - 1 \right)^2} \right] & = & - 1.
                        \end{array}
                    \end{equation}
                    In order to deal with the functions induced by the remaining four terms of the right-hand side of \eqref{eq:studyAmuKNeq0SixthPart14}, we note that we have
                    \begin{equation}
                    \label{eq:studyAmuKNeq0SixthPart19PowerSeriesExpSqrt1}
                        \begin{array}{lll}
                            \scriptstyle \sqrt{\frac{4 \pi^2 \left( x - \alpha \right)^2 a^2 + 1/4 + \mu}{4 \pi^2 \left( x + \alpha \right)^2 a^2 + 1/4 + \mu}} & = & \scriptstyle \left( 1 - \frac{16 \pi^2 a^2 \alpha x}{4 \pi^2 \left( x + \alpha \right)^2 a^2 + 1/4 + \mu} \right)^{1/2} \\[1em]
                            
                            & = & \scriptstyle - \frac{1}{2 \sqrt{\pi}} \sum\limits_{n=0}^{+\infty} \frac{1}{n!} \Gamma \left( n - \frac{1}{2} \right) \left( \frac{16 \pi^2 a^2 \alpha x}{4 \pi^2 \left( x + \alpha \right)^2 a^2 + 1/4 + \mu} \right)^n
                        \end{array}
                    \end{equation}
                    for any real number $x \geqslant 1$, using the binomial formula, recalled in \cite[Prop. C.26]{dutour:cuspDirichlet}. Furthermore, we have
                    \begin{equation}
                    \label{eq:studyAmuKNeq0SixthPart20PowerSeriesExpSqrt2}
                        \begin{array}{lll}
                            \scriptstyle \sqrt{\frac{4 \pi^2 \left( x + \alpha \right)^2 a^2 + 1/4 + \mu}{4 \pi^2 \left( x - \alpha \right)^2 a^2 + 1/4 + \mu}} & = & \scriptstyle \left( 1 - \frac{16 \pi^2 a^2 \alpha x}{4 \pi^2 \left( x + \alpha \right)^2 a^2 + 1/4 + \mu} \right)^{-1/2} \\[1em]
                            
                            & = & \scriptstyle \frac{1}{\sqrt{\pi}} \sum\limits_{n=0}^{+\infty} \frac{1}{n!} \Gamma \left( n + \frac{1}{2} \right) \left( \frac{16 \pi^2 a^2 \alpha x}{4 \pi^2 \left( x + \alpha \right)^2 a^2 + 1/4 + \mu} \right)^n.
                        \end{array}
                    \end{equation}
                    Let us begin the computation of the last four integrals in \eqref{eq:studyAmuKNeq0SixthPart14}. We have
                    \begin{equation}
                    \label{eq:studyAmuKNeq0SixthPart21}
                        \begin{array}{ll}
                            \multicolumn{2}{l}{\scriptstyle \frac{1}{4 \alpha} \cdot \frac{1}{2 \delta \left( s + j \right) - 1} \int_1^{+ \infty} x^{- 2 \delta \left( s + j \right) + 1} \left(\sqrt{\frac{4 \pi^2 \left( x + \alpha \right)^2 a^2 + 1/4 + \mu}{4 \pi^2 \left( x - \alpha \right)^2 a^2 + 1/4 + \mu}} - \sqrt{\frac{4 \pi^2 \left( x - \alpha \right)^2 a^2 + 1/4 + \mu}{4 \pi^2 \left( x + \alpha \right)^2 a^2 + 1/4 + \mu}} \right) \mathrm{d}x} \\[1em]
                            
                            = & \scriptstyle \frac{1}{4 \alpha \sqrt{\pi}} \cdot \frac{1}{2 \delta \left( s + j \right) - 1} \sum\limits_{n=1}^{+ \infty} \frac{1}{\left( n - 1 \right)!} \Gamma \left( n - \frac{1}{2} \right) \int_1^{+ \infty} x^{- 2 \delta \left( s + j \right) + 1} \left( \frac{16 \pi^2 \alpha^2 \alpha x}{4 \pi^2 \left( x + \alpha \right)^2 a^2 + 1/4 + \mu} \right)^n \mathrm{d}x
                        \end{array}
                    \end{equation}
                    and the last integral above can be computed using hypergeometric functions, and more precisely the integration formula \cite[Cor. C.32]{dutour:cuspDirichlet}. We have
                    \begin{equation}
                    \label{eq:studyAmuKNeq0SixthPart22}
                        \begin{array}{lll}
                            \multicolumn{3}{l}{\scriptstyle \int_1^{+ \infty} x^{- 2 \delta \left( s + j \right) + 1} \left( \frac{16 \pi^2 \alpha^2 \alpha x}{4 \pi^2 \left( x + \alpha \right)^2 a^2 + 1/4 + \mu} \right)^n \mathrm{d}x} \\[1em]
                            
                            & = & \scriptstyle \left( 4 \pi a \right)^{2n} \alpha^n \int_1^{+ \infty} x^{- 2 \delta \left( s + j \right) + 1 + n} \left( 4 \pi^2 \left( x + \alpha \right)^2 a^2 + \frac{1}{4} + \mu \right)^{-n} \mathrm{d}x \\[0.5em]
                            
                            & = & \scriptstyle 2^{2n-1} \left( 2 \pi a \right)^{2 \delta \left( s + j \right) + n - 2} \left( \frac{1}{4} + \mu \right)^{- \delta \left( s + j \right) + 1 - n/2} \alpha^n  \\[0.5em]
                            
                            && \qquad \qquad \quad \scriptstyle \cdot \int_0^{\frac{1/4+\mu}{4 \pi^2 \left( 1 + \alpha \right)^2 a^2}} \left( 1 - \frac{2 \pi \alpha a}{\sqrt{1/4 + \mu}} t^{1/2} \right)^{- 2 \delta \left( s + j \right) + n + 1} \frac{t^{\delta \left( s + j \right) - 2 + n/2}}{\left( 1 + t \right)^n} \mathrm{d}t,
                        \end{array}
                    \end{equation}
                    by setting $t = \frac{1/4+\mu}{4 \pi^2 \left( x + \alpha \right)^2 a^2}$. On the interval of integration, we have
                    \begin{equation}
                    \label{eq:studyAmuKNeq0SixthPart23}
                        \begin{array}{lllllllll}
                            \scriptstyle 0 & \leqslant & \scriptstyle \frac{2 \pi \alpha a}{\sqrt{1/4 + \mu}} t^{1/2} & \leqslant & \scriptstyle \frac{2 \pi a \alpha}{\sqrt{1/4 + \mu}} \cdot \frac{\sqrt{1/4 + \mu}}{2 \pi \left( 1 + \alpha \right) a} & = & \scriptstyle \frac{\alpha}{1 + \alpha} & < & \scriptstyle 1,
                        \end{array}
                    \end{equation}
                    which allows us to use the binomial formula to expand the complex power in \eqref{eq:studyAmuKNeq0SixthPart22}, as well as interchange the sum and the integral. We have
                    \begin{equation}
                    \label{eq:studyAmuKNeq0SixthPart24}
                        \begin{array}{lll}
                            \scriptstyle \left( 1 - \frac{2 \pi \alpha a}{\sqrt{1/4 + \mu}} t^{1/2} \right)^{- 2 \delta \left( s + j \right) + n + 1} & = & \scriptstyle \sum\limits_{m=0}^{+ \infty} \frac{\Gamma \left( 2 \delta \left( s + j \right) - 1 - n + m \right)}{\Gamma \left( 2 \delta \left( s + j \right) - 1 - n \right)} \cdot \frac{1}{m!} \left( \frac{2 \pi \alpha a}{\sqrt{1/4 + \mu}} t^{1/2} \right)^m
                        \end{array}
                    \end{equation}
                    and we can plug \eqref{eq:studyAmuKNeq0SixthPart24} into \eqref{eq:studyAmuKNeq0SixthPart22}, yielding
                    \begin{equation}
                    \label{eq:studyAmuKNeq0SixthPart25}
                        \begin{array}{lll}
                            \multicolumn{3}{l}{\scriptstyle \int_1^{+ \infty} x^{- 2 \delta \left( s + j \right) + 1} \left( \frac{16 \pi^2 \alpha^2 \alpha x}{4 \pi^2 \left( x + \alpha \right)^2 a^2 + 1/4 + \mu} \right)^n \mathrm{d}x} \\[1em]
                            
                            & = & \scriptstyle 2^{2n-1} \left( 2 \pi a \right)^{2 \delta \left( s + j \right) + n - 2} \left( \frac{1}{4} + \mu \right)^{- \delta \left( s + j \right) + 1 - n/2} \alpha^n  \sum\limits_{m=0}^{+ \infty} \Big[ \frac{\Gamma \left( 2 \delta \left( s + j \right) - 1 - n + m \right)}{\Gamma \left( 2 \delta \left( s + j \right) - 1 - n \right)} \\[0.5em]
                            
                            && \scriptstyle \qquad \cdot \frac{1}{m!} \left( \frac{1}{4} + \mu \right)^{-m/2} \alpha^m \left( 2 \pi a \right)^m \int_0^{\frac{1/4+\mu}{4 \pi^2 \left( 1 + \alpha \right)^2 a^2}} \frac{1}{\left( 1 + t \right)^n} \cdot t^{\delta \left( s + j \right) - 2 + \left( n+m \right)/2} \mathrm{d}t \Big].
                        \end{array}
                    \end{equation}
                    The last integral in \eqref{eq:studyAmuKNeq0SixthPart25} can be computed using \cite[Cor. C.32]{dutour:cuspDirichlet}, and we have
                    \begin{equation}
                    \label{eq:studyAmuKNeq0SixthPart26}
                        \begin{array}{lll}
                            \multicolumn{3}{l}{\scriptstyle \int_0^{\frac{1/4+\mu}{4 \pi^2 \left( 1 + \alpha \right)^2 a^2}} \frac{1}{\left( 1 + t \right)^n} \cdot t^{\delta \left( s + j \right) - 2 + \left( n+m \right)/2} \mathrm{d}t} \\[1em]
                            
                            \qquad & = & \scriptstyle \frac{2}{2 \delta \left( s + j \right) + n + m - 2} \left( \frac{1/4 + \mu}{4 \pi^2 \left( 1 + \alpha \right)^2 a^2} \right)^{\delta \left( s + j \right) - 1 + \left( n+m \right)/2} \left( 1 + \frac{1/4 + \mu}{4 \pi^2 \left( 1 + \alpha \right)^2 a^2} \right)^{-n} \\[1em]
                            
                            && \scriptstyle \qquad \qquad \qquad \qquad \cdot F \left( n, \; 1; \; \delta \left( s + j \right) + \frac{n+m}{2}; \; \frac{1/4 + \mu}{4 \pi^2 \left( 1 + \alpha \right)^2 a^2 + 1/4 + \mu} \right).
                        \end{array}
                    \end{equation}
                    Using this computation together with \eqref{eq:studyAmuKNeq0SixthPart25}, we get
                    \begin{equation}
                    \label{eq:studyAmuKNeq0SixthPart27}
                        \begin{array}{ll}
                            \multicolumn{2}{l}{\scriptstyle \int_1^{+ \infty} x^{- 2 \delta \left( s + j \right) + 1} \left( \frac{16 \pi^2 \alpha^2 \alpha x}{4 \pi^2 \left( x + \alpha \right)^2 a^2 + 1/4 + \mu} \right)^n \mathrm{d}x} \\[1em]
                            
                            = & \scriptstyle \frac{1}{\left( 4 \pi^2 \left( 1 + \alpha \right)^2 a^2 + 1/4 + \mu \right)^n} \left( 1 + \alpha \right)^{-2 \delta \left( s + j \right) + n + 2} \sum\limits_{m=0}^{+ \infty} \left[ \frac{1}{m!} \cdot \frac{\Gamma \left( 2 \delta \left( s + j \right) - 1 - n + m \right)}{\Gamma \left( 2 \delta \left( s + j \right) - 1 - n \right)} \left( \frac{\alpha}{1 + \alpha} \right)^m \right. \\[1em]

                            & \scriptstyle \left. \qquad \cdot \frac{1}{2 \delta \left( s + j \right) + n + m - 2} F \left( n, \; 1; \; \delta \left( s + j \right) + \frac{n+m}{2}; \; \frac{1/4 + \mu}{4 \pi^2 \left( 1 + \alpha \right)^2 a^2 + 1/4 + \mu} \right) \right].
                        \end{array}
                    \end{equation}
                    Inserting \eqref{eq:studyAmuKNeq0SixthPart27} into \eqref{eq:studyAmuKNeq0SixthPart21}, we get
                    \begin{equation}
                    \label{eq:studyAmuKNeq0SixthPart28}
                        \begin{array}{lll}
                            \multicolumn{3}{l}{\scriptstyle \frac{1}{4 \alpha} \cdot \frac{1}{2 \delta \left( s + j \right) - 1} \int_1^{+ \infty} x^{- 2 \delta \left( s + j \right) + 1} \left( \sqrt{\frac{4 \pi^2 \left( x + \alpha \right)^2 a^2 + 1/4 + \mu}{4 \pi^2 \left( x - \alpha \right)^2 a^2 + 1/4 + \mu}} - \sqrt{\frac{4 \pi^2 \left( x - \alpha \right)^2 a^2 + 1/4 + \mu}{4 \pi^2 \left( x + \alpha \right)^2 a^2 + 1/4 + \mu}} \right) \mathrm{d}x} \\[0.5em]
                            
                            & = & \scriptstyle \frac{1}{4 \alpha \sqrt{\pi}} \cdot \frac{1}{2 \delta \left( s + j \right) - 1} \sum\limits_{n=1}^{+ \infty} \left[ \frac{1}{\left( n - 1 \right)!} \cdot \Gamma \left( n - \frac{1}{2} \right) \left( 4 \pi a \right)^{2n} \alpha^n \cdot \frac{1}{\left( 4 \pi^2 \left( 1 + \alpha \right)^2 a^2 + 1/4 + \mu \right)^n} \right. \\[1em]
                            
                            && \scriptstyle \quad \cdot \left( 1 + \alpha \right)^{- 2 \delta \left( s + j \right) + n + 2} \sum\limits_{m=0}^{+ \infty} \frac{1}{m!} \cdot \frac{\Gamma \left( 2 \delta \left( s + j \right) - 1 - n + m \right)}{\Gamma \left( 2 \delta \left( s + j \right) - 1 - n \right)} \left( \frac{\alpha}{1 + \alpha} \right)^m \cdot \frac{1}{2 \delta \left( s + j \right) + n + m - 2} \\[1em]
                            
                            && \scriptstyle \left. \qquad \qquad \quad \qquad \qquad \cdot F \left( n, \; 1; \; \delta \left( s + j \right) + \frac{n+m}{2}; \; \frac{1/4 + \mu}{4 \pi^2 \left( 1 + \alpha \right)^2 a^2 + 1/4 + \mu} \right) \right].
                        \end{array}
                    \end{equation}
                    Similarly, we have
                    \begin{equation}
                    \label{eq:studyAmuKNeq0SixthPart29}
                        \begin{array}{lll}
                            \multicolumn{3}{l}{\scriptstyle - \frac{1}{4} \cdot \frac{1}{2 \delta \left( s + j \right) - 1} \int_1^{+ \infty} x^{- 2 \delta \left( s + j \right)} \left( \sqrt{\frac{4 \pi^2 \left( x + \alpha \right)^2 a^2 + 1/4 + \mu}{4 \pi^2 \left( x - \alpha \right)^2 a^2 + 1/4 + \mu}} + \sqrt{\frac{4 \pi^2 \left( x - \alpha \right)^2 a^2 + 1/4 + \mu}{4 \pi^2 \left( x + \alpha \right)^2 a^2 + 1/4 + \mu}} \right) \mathrm{d}x} \\[0.5em]
                            
                            & = & \scriptstyle - \frac{1}{4 \sqrt{\pi}} \cdot \frac{1}{2 \delta \left( s + j \right) - 1} \sum\limits_{n=0}^{+ \infty} \left[ \frac{n-1}{n!} \Gamma \left( n - \frac{1}{2} \right) \left( 4 \pi a \right)^{2n} \alpha^n \cdot \frac{1}{\left( 4 \pi^2 \left( 1 + \alpha \right)^2 a^2 + 1/4 + \mu \right)^n} \right. \\[1em]
                            
                            && \scriptstyle \quad \cdot \left( 1 + \alpha \right)^{-2 \delta \left( s + j \right) + n + 1} \sum\limits_{m=0}^{+ \infty} \frac{1}{m!} \cdot \frac{\Gamma \left( 2 \delta \left( s + j \right) - n + m \right)}{\Gamma \left( 2 \delta \left( s + j \right) - n \right)} \left( \frac{\alpha}{1 + \alpha} \right)^m \cdot \frac{1}{2 \delta \left( s + j \right) + n + m - 1} \\[1em]
                            
                            && \scriptstyle \left. \qquad \qquad \quad \qquad \qquad \cdot F \left( n, \; 1; \; \delta \left( s + j \right) + \frac{n+m+1}{2}; \; \frac{1/4 + \mu}{4 \pi^2 \left( 1 + \alpha \right)^2 a^2 + 1/4 + \mu} \right) \right].
                        \end{array}
                    \end{equation}
                    We will study \eqref{eq:studyAmuKNeq0SixthPart28} and \eqref{eq:studyAmuKNeq0SixthPart29} separately. The term induced by \eqref{eq:studyAmuKNeq0SixthPart28} is
                    \begin{equation}
                    \label{eq:studyAmuKNeq0SixthPart30}
                        \begin{array}{l}
                            \scriptstyle - \frac{1}{2 \alpha \sqrt{\pi}} \cdot \frac{\sin \left( \pi s \right)}{\pi} \left( 4 \mu + 1 \right)^{-s} \sum\limits_{j=0}^{+ \infty} \frac{\left( s \right)_j}{j! 4^j} \cdot \frac{1}{2 \delta \left( s + j \right) - 1} \sum\limits_{n=1}^{+ \infty} \left[ \frac{1}{\left( n - 1 \right)!} \Gamma \left( n - \frac{1}{2} \right) \left( 4 \pi a \right)^{2n} \alpha^n \right. \\[1em]
                            
                            \scriptstyle \qquad \cdot \frac{1}{\left( 4 \pi^2 \left( 1 + \alpha \right)^2 a^2 + 1/4 + \mu \right)^n} \cdot \left( 1 + \alpha \right)^{- 2 \delta \left( s + j \right) + n + 2} \sum\limits_{m=0}^{+ \infty} \frac{1}{m!} \cdot \frac{\Gamma \left( 2 \delta \left( s + j \right) - 1 - n + m \right)}{\Gamma \left( 2 \delta \left( s + j \right) - 1 - n \right)} \left( \frac{\alpha}{1 + \alpha} \right)^m \\[1em]
                            
                            \scriptstyle \left. \qquad \qquad \qquad \cdot \frac{1}{2 \delta \left( s + j \right) + n + m - 2} F \left( n, \; 1; \; \delta \left( s + j \right) + \frac{n+m}{2}; \; \frac{1/4 + \mu}{4 \pi^2 \left( 1 + \alpha \right)^2 a^2 + 1/4 + \mu} \right) \right].
                        \end{array}
                    \end{equation}
                    We will first study the sum over $j \geqslant 1$ within \eqref{eq:studyAmuKNeq0SixthPart30}. Using the contiguous functions formula
                    \begin{equation}
                    \label{eq:studyAmuKNeq0SixthPart31}
                        \begin{array}{lll}
                            \scriptstyle c \left( 1 - t \right) F \left( a, \; b; \; c; \; t \right) - c F \left( a, \; b-1; \; c; \; t \right) + \left( c - a \right) t F \left( a, \; b; \; c+1; \; t \right) & = & \scriptstyle 0,
                        \end{array}
                    \end{equation}
                    which one may check directly for $0 \leqslant t < 1$ using the hypergeometric series, we get
                    \begin{equation}
                    \label{eq:studyAmuKNeq0SixthPart32}
                        \begin{array}{lll}
                            \multicolumn{3}{l}{\scriptstyle F \left( n, \; 1; \; \delta \left( s + j \right) + \frac{n+m}{2}; \; \frac{1/4 + \mu}{4 \pi^2 \left( 1 + \alpha \right)^2 a^2 + 1/4 + \mu} \right)} \\[0.5em]
                            
                            \qquad & = & \scriptstyle 1 + \frac{1/4 + \mu}{4 \pi^2 \left( 1 + \alpha \right)^2 a^2} - \frac{2 \delta \left( s + j \right) + m - n}{2 \delta \left( s + j \right) + m + n} \cdot \frac{1/4 + \mu}{4 \pi^2 \left( 1 + \alpha \right)^2 a^2} \\[0.5em]
                            
                            && \scriptstyle \qquad \qquad \qquad \qquad \quad \cdot F \left( n, \; 1; \; \delta \left( s + j \right) + \frac{n+m}{2} + 1; \; \frac{1/4 + \mu}{4 \pi^2 \left( 1 + \alpha \right)^2 a^2 + 1/4 + \mu} \right).
                        \end{array}
                    \end{equation}
                    The term induced by the first two parts of the right-hand side of \eqref{eq:studyAmuKNeq0SixthPart32} is
                    \begin{equation}
                    \label{eq:studyAmuKNeq0SixthPart33}
                        \begin{array}{l}
                            \scriptstyle - \frac{2}{\alpha \sqrt{\pi}} \cdot \frac{\sin \left( \pi s \right)}{\pi} \left( 4 \mu + 1 \right)^{-s} \sum\limits_{j=1}^{+ \infty} \frac{\left( s \right)_j}{j! 4^j} \cdot \frac{1}{2 \delta \left( s + j \right) - 1} \sum\limits_{n=1}^{+ \infty} \left[ \frac{1}{\left( n - 1 \right)!} \Gamma \left( n - \frac{1}{2} \right)  \alpha^n \left( 1 + \alpha \right)^{- 2 \delta \left( s + j \right) + n} \right. \\[1em]
                            
                            \scriptstyle \left. \quad \cdot \frac{\left( 4 \pi a \right)^{2 \left( n - 1 \right)}}{\left( 4 \pi^2 \left( 1 + \alpha \right)^2 a^2 + 1/4 + \mu \right)^{n-1}} \sum\limits_{m=0}^{+ \infty} \frac{1}{m!} \cdot \frac{\Gamma \left( 2 \delta \left( s + j \right) - 1 - n + m \right)}{\Gamma \left( 2 \delta \left( s + j \right) - 1 - n \right)} \left( \frac{\alpha}{1 + \alpha} \right)^m \cdot \frac{1}{2 \delta \left( s + j \right) + n + m - 2} \right].
                        \end{array}
                    \end{equation}
                    For $\left \vert s \right \vert \leqslant \varepsilon$ small enough, we have
                    \begin{equation}
                    \label{eq:studyAmuKNeq0SixthPart34}
                        \begin{array}{lll}
                            \multicolumn{3}{l}{\scriptstyle \left \vert \frac{1}{\left( n - 1 \right)!} \Gamma \left( n - \frac{1}{2} \right)  \alpha^n \left( 1 + \alpha \right)^{- 2 \delta \left( s + j \right) + n} \frac{\left( 4 \pi a \right)^{2 \left( n - 1 \right)}}{\left( 4 \pi^2 \left( 1 + \alpha \right)^2 a^2 + 1/4 + \mu \right)^{n-1}} \right.} \\[0.5em]
                            
                            \multicolumn{3}{l}{\scriptstyle \left. \qquad \qquad \qquad \qquad \qquad \qquad \cdot \frac{1}{m!} \cdot \frac{\Gamma \left( 2 \delta \left( s + j \right) - 1 - n + m \right)}{\Gamma \left( 2 \delta \left( s + j \right) - 1 - n \right)} \left( \frac{\alpha}{1 + \alpha} \right)^m \cdot \frac{1}{2 \delta \left( s + j \right) + n + m - 2} \right \vert} \\[1em]
                            
                            & \leqslant & \scriptstyle \frac{1}{\left( n - 1 \right)!} \Gamma \left( n - \frac{1}{2} \right)  \alpha^n \left( 1 + \alpha \right)^{- 2 \delta \left( j - \varepsilon \right)+n} \left( 4 \pi a \right)^{2 \left( n-1 \right)} \frac{1}{\left( 4 \pi^2 \left( 1 + \alpha \right)^2 a^2 + 1/4 + \mu \right)^{n-1}} \\[0.5em]
                            
                            && \qquad \scriptstyle \cdot \frac{1}{2 \delta \left( j - \varepsilon \right)} \cdot \frac{1}{m!} \cdot \frac{\Gamma \left( 2 \delta \left( \varepsilon + j \right) + 1 + n + m \right)}{\Gamma \left( 2 \delta \left( \varepsilon + j \right) + 1 + n \right)} \cdot \left( \frac{\alpha}{1 + \alpha} \right)^m
                        \end{array}
                    \end{equation}
                    if we have $n \geqslant 2$, or $m \geqslant 1$, or $j$ large enough. The estimate \eqref{eq:studyAmuKNeq0SixthPart34}, the binomial formula to find an upper-bound for the part that depends on $m$, and Stirling's formula prove that \eqref{eq:studyAmuKNeq0SixthPart33} induces, for $\mu$ large enough, a holomorphic function around $0$, whose derivative at $s=0$ vanishes because of the Pochhammer symbol. The term induced by the last part of the right-hand side of \eqref{eq:studyAmuKNeq0SixthPart32} is
                    \begin{equation}
                    \label{eq:studyAmuKNeq0SixthPart35}
                        \begin{array}{lll}
                            \scriptstyle \frac{1}{2 \alpha \sqrt{\pi}} \cdot \frac{\sin \left( \pi s \right)}{\pi} \left( 4 \mu + 1 \right)^{-s+1} \sum\limits_{j=1}^{+ \infty} \frac{\left( s \right)_j}{j! 4^j} \cdot \frac{1}{2 \delta \left( s + j \right) - 1} \sum\limits_{n=1}^{+ \infty} \left[ \frac{1}{\left( n - 1 \right)!} \Gamma \left( n - \frac{1}{2} \right) \left( 4 \pi a \right)^{2 \left( n - 1 \right)} \alpha^n \right. \\[1em]
                            
                            \scriptstyle \quad \cdot \left( 1 + \alpha \right)^{- 2 \delta \left( s + j \right) + n} \frac{1}{\left( 4 \pi^2 \left( 1 + \alpha \right)^2 a^2 + 1/4 + \mu \right)^n} \cdot \sum\limits_{m=0}^{+ \infty} \frac{1}{m!} \cdot \frac{\Gamma \left( 2 \delta \left( s + j \right) - 1 - n + m \right)}{\Gamma \left( 2 \delta \left( s + j \right) - 1 - n \right)} \left( \frac{\alpha}{1 + \alpha} \right)^m \\[1em]
                            
                            \scriptstyle \left. \quad \quad \cdot \frac{1}{2 \delta \left( s + j \right) + n + m - 2} \cdot \frac{2 \delta \left( s + j \right) + m - n}{2 \delta \left( s + j \right) + m + n} F \left( n, \; 1; \; \delta \left( s + j \right) + \frac{n+m}{2} + 1; \; \frac{1/4 + \mu}{4 \pi^2 \left( 1 + \alpha \right)^2 a^2 + 1/4 + \mu} \right) \right].
                        \end{array}
                    \end{equation}
                    The Euler integral formula, recalled in \cite[Prop. C.30]{dutour:cuspDirichlet}, yields, for any $t \in \left[ 0, 1 \right[$,
                    \begin{equation}
                    \label{eq:studyAmuKNeq0SixthPart36}
                        \begin{array}{lll}
                            \multicolumn{3}{l}{\scriptstyle F \left( n, \; 1; \; \delta \left( s + j \right) + \frac{n+m}{2} + 1; \; t \right)} \\[0.5em]
                            
                            \qquad \qquad \qquad & = & \scriptstyle \frac{1}{2} \left( 2 \delta \left( s + j \right) + n + m \right) \int_0^1 \left( 1 - x \right)^{\delta \left( s + j \right) + \left( n + m \right)/2 - 1} \left( 1 - tx \right)^{-n} \mathrm{d}x
                        \end{array}
                    \end{equation}
                    which in turn gives the following estimate
                    \begin{equation}
                    \label{eq:studyAmuKNeq0SixthPart37}
                        \begin{array}{lll}
                            \scriptstyle \left \vert \frac{1}{2 \delta \left( s + j \right) + n + m} F \left( n, \; 1; \; \delta \left( s + j \right) + \frac{n+m}{2} + 1; \; t \right) \right \vert & \leqslant & \scriptstyle \frac{1}{2 \delta \left( j - \varepsilon \right)} \left( 1 - t \right)^{-n/2}
                        \end{array}
                    \end{equation}
                    where $\left \vert s \right \vert \leqslant \varepsilon$ is small enough. We get, for $n \geqslant 1$,
                    \begin{equation}
                    \label{eq:studyAmuKNeq0SixthPart38}
                        \begin{array}{lll}
                            \multicolumn{3}{l}{\scriptstyle \left \vert \frac{1}{\left( n - 1 \right)!} \Gamma \left( n - \frac{1}{2} \right) \left( 4 \pi a \right)^{2 \left( n - 1 \right)} \alpha^n \left( 1 + \alpha \right)^{- 2 \delta \left( s + j \right) + n} \frac{1}{\left( 4 \pi^2 \left( 1 + \alpha \right)^2 a^2 + 1/4 + \mu \right)^n} \right.} \\[0.5em]
                            
                            \multicolumn{3}{l}{\scriptstyle \qquad \qquad \qquad \cdot \frac{1}{m!} \cdot \frac{\Gamma \left( 2 \delta \left( s + j \right) - 1 - n + m \right)}{\Gamma \left( 2 \delta \left( s + j \right) - 1 - n \right)} \left( \frac{\alpha}{1 + \alpha} \right)^m \cdot \frac{1}{2 \delta \left( s + j \right) + n + m - 2} \cdot \frac{2 \delta \left( s + j \right) + m - n}{2 \delta \left( s + j \right) + m + n}} \\[0.5em]
                            
                            \multicolumn{3}{l}{\scriptstyle \qquad \qquad \qquad \qquad \qquad \qquad \left. F \left( n, \; 1; \; \delta \left( s + j \right) + \frac{n+m}{2} + 1; \; \frac{1/4 + \mu}{4 \pi^2 \left( 1 + \alpha \right)^2 a^2 + 1/4 + \mu} \right) \right \vert} \\[1em]
                            
                            & \leqslant & \scriptstyle \frac{1}{\left( n - 1 \right)!} \Gamma \left( n - \frac{1}{2} \right) \left( 4 \pi a \right)^{2 \left( n - 1 \right)} \alpha^n \left( 1 + \alpha \right)^{- 2 \delta \left( j - \varepsilon \right) + n} \frac{1}{\left( 4 \pi^2 \left( 1 + \alpha \right)^2 a^2 + 1/4 + \mu \right)^n} \\
\end{array}
                    \end{equation}
                    for $\left \vert s \right \vert \leqslant \varepsilon$ small enough, if we have $n \geqslant 2$, or $m \geqslant 1$, or $j$ large enough. The estimate \eqref{eq:studyAmuKNeq0SixthPart38} and the binomial formula thus prove that \eqref{eq:studyAmuKNeq0SixthPart35} induces, for $\mu$ large enough, a holomorphic function around $0$, whose derivative at $s=0$ vanishes because of the Pochhammer symbol. With respect to \eqref{eq:studyAmuKNeq0SixthPart30}, only the term $j=0$ remains. It is given by
                    \begin{equation}
                    \label{eq:studyAmuKNeq0SixthPart39}
                        \begin{array}{l}
                            \scriptstyle - \frac{1}{2 \alpha \sqrt{\pi}} \cdot \frac{\sin \left( \pi s \right)}{\pi} \left( 4 \mu + 1 \right)^{-s} \frac{1}{2 \delta s - 1} \sum\limits_{n=1}^{+ \infty} \left[ \frac{1}{\left( n - 1 \right)!} \Gamma \left( n - \frac{1}{2} \right) \left( 4 \pi a \right)^{2n} \alpha^n \left( 1 + \alpha \right)^{- 2 \delta s + n + 2} \right. \\[0.5em]
                            
                            \scriptstyle \qquad \qquad \cdot \frac{1}{\left( 4 \pi^2 \left( 1 + \alpha \right)^2 a^2 + 1/4 + \mu \right)^n} \sum\limits_{m=0}^{+ \infty} \frac{1}{m!} \cdot \frac{\Gamma \left( 2 \delta s - 1 - n + m \right)}{\Gamma \left( 2 \delta s - 1 - n \right)} \left( \frac{\alpha}{1 + \alpha} \right)^m \\[0.5em]
                            
                            \scriptstyle \qquad \qquad \qquad \qquad \qquad\left. \cdot \frac{1}{2 \delta s + n + m - 2} F \left( n, \; 1; \; \delta s + \frac{n+m}{2}; \; \frac{1/4 + \mu}{4 \pi^2 \left( 1 + \alpha \right)^2 a^2 + 1/4 + \mu} \right) \right].
                        \end{array}
                    \end{equation}
                    Using arguments similar to those presented for the sum over $j \geqslant 1$, the sum ranging over $m \geqslant 2$ included in \eqref{eq:studyAmuKNeq0SixthPart39} is seen to induce a holomorphic function around $0$, whose derivative at $s=0$ equals
                    \begin{equation}
                    \label{eq:studyAmuKNeq0SixthPart40}
                        \begin{array}{l}
                            \scriptstyle \frac{1}{2 \alpha \sqrt{\pi}} \sum\limits_{n=1}^{+ \infty} \left[ \frac{1}{\left( n - 1 \right)!} \Gamma \left( n - \frac{1}{2} \right) \left( 4 \pi a \right)^{2n} \left( 1 + \alpha \right)^{n+2} \frac{1}{\left( 4 \pi^2 \left( 1 + \alpha \right)^2 a^2 + 1/4 + \mu \right)^n} \alpha^n \right. \\[0.5em]
                            
                            \scriptstyle \; \; \left. \cdot \sum\limits_{m=2}^{n+1} \frac{1}{m!} \Fp_{s=0} \left( \frac{\Gamma \left( 2 \delta s - n - 1 + m \right)}{\Gamma \left( 2 \delta s - n - 1 \right)} \right) \left( \frac{\alpha}{1 + \alpha} \right)^m \cdot \frac{1}{n + m - 2} F \left( n, \; 1; \; \frac{n+m}{2}; \; \frac{1/4 + \mu}{4 \pi^2 \left( 1 + \alpha \right)^2 a^2 + 1/4 + \mu} \right) \right]
                        \end{array}
                    \end{equation}
                    since the quotient of Gamma functions vanishes for other integers $m$. We have
                    \begin{equation}
                    \label{eq:studyAmuKNeq0SixthPart41}
                        \begin{array}{lll}
                            \Fp_{s=0} \left( \frac{\Gamma \left( 2 \delta s - 1 - n + m \right)}{\Gamma \left( 2 \delta s - 1 - n \right)} \right) & = & \left( -1 \right)^m \cdot \frac{\left( n + 1 \right)!}{\left( n + 1 - m \right)!}
                        \end{array}
                    \end{equation}
                    using the Laurent series expansion of the Gamma function. We now need to modify the hypergeometric function in \eqref{eq:studyAmuKNeq0SixthPart40} in order to make its $\mu$-asymptotic study possible. To that effect, we use an analytic extension of \cite[Prop C.29]{dutour:cuspDirichlet}, which reads
                    \begin{equation}
                    \label{eq:studyAmuKNeq0SixthPart42AnalyticExtHyperGeoT1MinusT}
                        \begin{array}{lll}
                            \scriptstyle F \left( a, \; 1; \; c; \; t \right) & = & \scriptstyle \frac{1}{\Gamma \left( a \right)} \Gamma \left( c \right) \Gamma \left( a+1-c \right) \left( 1 - t \right)^{c-a-1} t^{1-c} + \frac{c-1}{c-a-1} F \left( a, \; 1; \; 2+a-c; \; 1-t \right)
                        \end{array}
                    \end{equation}
                    for any $t \in \left] 0,1 \right[$. Thus, we have, for any integers $n \geqslant 2$ and $1 \leqslant m \leqslant n+1$,
                    \begin{equation}
                    \label{eq:studyAmuKNeq0SixthPart43}
                        \begin{array}{lll}
                            \multicolumn{3}{l}{\scriptstyle F \left( n, \; 1; \; \frac{n+m}{2}; \; \frac{1/4 + \mu}{4 \pi^2 \left( 1 + \alpha \right)^2 a^2 + 1/4 + \mu} \right)} \\[1em]
                            
                            & = & \scriptstyle \frac{1}{\left( n-1 \right)!} \Gamma \left( \frac{n+m}{2} \right) \Gamma \left( \frac{n-m}{2} + 1 \right) \left( 2 \pi a \right)^{m-n-2} \left( 1 + \alpha \right)^{m-n-2} \\[0.5em]
                            
                            && \scriptstyle \qquad \qquad \qquad \qquad \qquad \qquad \cdot \left( 4 \pi^2 \left( 1 + \alpha \right)^2 a^2 + 1/4 + \mu \right)^n \left( \frac{1}{4} + \mu \right)^{1 - \frac{n+m}{2}} \\[0.5em]
                            
                            && \scriptstyle \qquad \qquad \qquad - \frac{n+m-2}{n+2-m} F \left( n, \; 1; \; \frac{n-m}{2} + 2; \; \frac{4 \pi^2 \left( 1 + \alpha \right)^2 a^2}{4 \pi^2 \left( 1 + \alpha \right)^2 a^2 + 1/4 + \mu} \right).
                        \end{array}
                    \end{equation}
                    The term induced by the first part of \eqref{eq:studyAmuKNeq0SixthPart43} is
                    \begin{equation}
                    \label{eq:studyAmuKNeq0SixthPart44}
                        %
                    \end{equation}
                    \noindent Using either an upper-bound or the explicit special value of the Gamma function at half-integers for the last two Gamma functions, one sees that \eqref{eq:studyAmuKNeq0SixthPart44} vanishes as $\mu$ goes to infinity. The term induced by the second part of \eqref{eq:studyAmuKNeq0SixthPart43} is
                    \begin{equation}
                    \label{eq:studyAmuKNeq0SixthPart45}
                        %
                    \end{equation}
                    Using the Euler integral formula, recalled in \cite[Prop. C.30]{dutour:cuspDirichlet}, we have
                    \begin{equation}
                    \label{eq:studyAmuKNeq0SixthPart46}
                        %
                    \end{equation}
                    which gives the estimate
                    \begin{equation}
                    \label{eq:studyAmuKNeq0SixthPart47}
                        %
                    \end{equation}
                    for any real number $t \in \left[ 0, 1 \right[$. We thus get
                    \begin{equation}
                    \label{eq:studyAmuKNeq0SixthPart48}
                        %
                    \end{equation}
                    Using estimate \eqref{eq:studyAmuKNeq0SixthPart48}, we can prove that the term \eqref{eq:studyAmuKNeq0SixthPart45} converges for $\mu$ large enough, and vanishes as $\mu$ goes to infinity. In \eqref{eq:studyAmuKNeq0SixthPart39}, only the terms $m=0,1$ remain. For $m=1$, we have
                    \begin{equation}
                    \label{eq:studyAmuKNeq0SixthPart49}
                        %
                    \end{equation}
                    The term induced by the first part of the right-hand side of \eqref{eq:studyAmuKNeq0SixthPart50} is given by
                    \begin{equation}
                    \label{eq:studyAmuKNeq0SixthPart51}
                        %
                    \end{equation}
                    Therefore, the sum over $n \geqslant 2$ in \eqref{eq:studyAmuKNeq0SixthPart51} induces, for $\mu$ large enough, a holomorphic function around $0$, and its derivative at $s=0$ vanishes as $\mu$ goes to infinity. The term $n=1$ in \eqref{eq:studyAmuKNeq0SixthPart51} is
                    \begin{equation}
                    \label{eq:studyAmuKNeq0SixthPart52}
                        %
                    \end{equation}
                    which induces a holomorphic function around $0$. Its derivative at $s=0$ can be computed using Laurent series expansion. Namely, we can write \eqref{eq:studyAmuKNeq0SixthPart52} as
                    \begin{equation}
                    \label{eq:studyAmuKNeq0SixthPart53}
                        %
                    \end{equation}
                    and the derivative at $s=0$ of \eqref{eq:studyAmuKNeq0SixthPart52} equals
                    \begin{equation}
                    \label{eq:studyAmuKNeq0SixthPart54}
                        %
                    \end{equation}
                    as $\mu$ goes to infinity. The term induced by the second part of \eqref{eq:studyAmuKNeq0SixthPart50} is given by
                    \begin{equation}
                    \label{eq:studyAmuKNeq0SixthPart55}
                        %
                    \end{equation}
                    Using the Euler integral formula, the term \eqref{eq:studyAmuKNeq0SixthPart55} induces, for $\mu$ large enough, a holomorphic function around $0$, whose derivative at $s=0$ vanishes as $\mu$ goes to infinity. In \eqref{eq:studyAmuKNeq0SixthPart39}, the term $m=0$ is
                    \begin{equation}
                    \label{eq:studyAmuKNeq0SixthPart56}
                        %
                    \end{equation}
                    The term induced by the first part of \eqref{eq:studyAmuKNeq0SixthPart57} is
                    \begin{equation}
                    \label{eq:studyAmuKNeq0SixthPart58}
                        %
                    \end{equation}
                    The sum over $n \geqslant 3$ in \eqref{eq:studyAmuKNeq0SixthPart58} thus induces, for $\mu$ large enough, a holomorphic function around $0$, whose derivative at $s=0$ vanishes as $\mu$ goes to infinity. The term $n=1$ in \eqref{eq:studyAmuKNeq0SixthPart58}, given by
                    \begin{equation}
                    \label{eq:studyAmuKNeq0SixthPart59}
                        %
                    \end{equation}
                    induces a holomorphic function around $0$, whose derivative at $s=0$ equals
                    \begin{equation}
                    \label{eq:studyAmuKNeq0SixthPart60}
                        %
                    \end{equation}
                    as $\mu$ goes to infinity. The term $n=2$ in \eqref{eq:studyAmuKNeq0SixthPart58}, given by
                    \begin{equation}
                    \label{eq:studyAmuKNeq0SixthPart61}
                        %
                    \end{equation}
                    also induces a holomorphic function around $0$. Its derivative at $s=0$ can be computed using Laurent series expansions. Namely, we can rewrite \eqref{eq:studyAmuKNeq0SixthPart61} as
                    \begin{equation}
                    \label{eq:studyAmuKNeq0SixthPart62}
                        %
                    \end{equation}
                    and the derivative at $s=0$ of \eqref{eq:studyAmuKNeq0SixthPart61} equals
                    \begin{equation}
                    \label{eq:studyAmuKNeq0SixthPart63}
                        %
                    \end{equation}
                    as $\mu$ goes to infinity. The term induced by the second part of \eqref{eq:studyAmuKNeq0SixthPart57} is
                    \begin{equation}
                    \label{eq:studyAmuKNeq0SixthPart64}
                        %
                    \end{equation}
                    Using the Euler integral formula, the term \eqref{eq:studyAmuKNeq0SixthPart64} induces, for $\mu$ large enough, a holomorphic function around $0$, whose derivative at $s=0$ vanishes as $\mu$ goes to infinity. Having finished with the study of \eqref{eq:studyAmuKNeq0SixthPart30}, we turn our attention to the term induced by \eqref{eq:studyAmuKNeq0SixthPart29}, given by
                    \begin{equation}
                    \label{eq:studyAmuKNeq0SixthPart65}
                        %
                    \end{equation}
                    We will begin by studying the sum over $j \geqslant 1$ in \eqref{eq:studyAmuKNeq0SixthPart65}. Using \eqref{eq:studyAmuKNeq0SixthPart31}, we have
                    \begin{equation}
                    \label{eq:studyAmuKNeq0SixthPart66}
                        %
                    \end{equation}
                    The term induced by the first two parts of \eqref{eq:studyAmuKNeq0SixthPart66} is
                    \begin{equation}
                    \label{eq:studyAmuKNeq0SixthPart67}
                        %
                    \end{equation}
                    We have, for $\left \vert s \right \vert \leqslant \varepsilon$ small enough,
                    \begin{equation}
                    \label{eq:studyAmuKNeq0SixthPart68}
                        %
                    \end{equation}
                    if we have $n \geqslant 1$, or $m \geqslant 1$, or $j$ large enough.  The estimate \eqref{eq:studyAmuKNeq0SixthPart68}, the binomial formula, and Stirling's formula prove that \eqref{eq:studyAmuKNeq0SixthPart67} induces, for $\mu$ large enough, a holomorphic function around $0$, whose derivative at $s=0$ vanishes because of the Pochhammer symbol, having assumed that we never have $2 \delta j = 1$. The term induced by the last part of \eqref{eq:studyAmuKNeq0SixthPart66} is
                    \begin{equation}
                    \label{eq:studyAmuKNeq0SixthPart69}
                        %
                    \end{equation}
                    Using the Euler integral formula, we have, for any $t \in \left[ 0, 1 \right[$,
                    \begin{equation}
                    \label{eq:studyAmuKNeq0SixthPart70}
                        %
                    \end{equation}
                    and we can use \eqref{eq:studyAmuKNeq0SixthPart70} to get, for $\left \vert s \right \vert \leqslant \varepsilon$ small enough,
                    \begin{equation}
                    \label{eq:studyAmuKNeq0SixthPart71}
                        %
                    \end{equation}
                    if we have $n \geqslant 1$, or $m \geqslant 1$, or $j$ large enough. The estimate \eqref{eq:studyAmuKNeq0SixthPart71}, the binomial formula, and Stirling's formula prove that \eqref{eq:studyAmuKNeq0SixthPart69} induces, for $\mu$ large enough, a holomorphic function around $0$, whose derivative at $s=0$ vanishes because of the Pochhammer symbol, having assumed that we never have $2 \delta j = 1$. In \eqref{eq:studyAmuKNeq0SixthPart65}, only the term $j=0$ remains. It is given by
                    \begin{equation}
                    \label{eq:studyAmuKNeq0SixthPart72}
                        %
                    \end{equation}
                    Using arguments similar to those presented for the sum over $j \geqslant 1$, the sum ranging over $m \geqslant 2$ included in \eqref{eq:studyAmuKNeq0SixthPart72} is seen to induce, for $\mu$ large enough, a holomorphic function around $0$, whose derivative at $s=0$ equals
                    \begin{equation}
                    \label{eq:studyAmuKNeq0SixthPart73}
                        %
                    \end{equation}
                    as the quotient of Gamma functions vanishes for other integers $m$ and $n$. We have
                    \begin{equation}
                    \label{eq:studyAmuKNeq0SixthPart74}
                        %
                    \end{equation}
                    using the Laurent series expansion of the Gamma function. Using \eqref{eq:studyAmuKNeq0SixthPart42AnalyticExtHyperGeoT1MinusT}, we have
                    \begin{equation}
                    \label{eq:studyAmuKNeq0SixthPart75}
                        %
                    \end{equation}
                    The term induced by the first part of \eqref{eq:studyAmuKNeq0SixthPart75} is
                    \begin{equation}
                    \label{eq:studyAmuKNeq0SixthPart76}
                        %
                    \end{equation}
                    For $\mu$ large enough, the term \eqref{eq:studyAmuKNeq0SixthPart76} converges, and actually vanishes as $\mu$ goes to infinity. The term induced by the last part of \eqref{eq:studyAmuKNeq0SixthPart75} is
                    \begin{equation}
                    \label{eq:studyAmuKNeq0SixthPart77}
                        %
                    \end{equation}
                    Using the Euler integral formula, we get the estimate
                    \begin{equation}
                    \label{eq:studyAmuKNeq0SixthPart78}
                        %
                    \end{equation}
                    for any integers $n \geqslant 2$ and $2 \leqslant m \leqslant n$. Using \eqref{eq:studyAmuKNeq0SixthPart78}, the term \eqref{eq:studyAmuKNeq0SixthPart77} is seen to converge, for $\mu$ large enough, and actually vanishes as $\mu$ goes to infinity. In \eqref{eq:studyAmuKNeq0SixthPart72}, only the terms $m=0,1$ remain. The term $m=1$ is
                    \begin{equation}
                    \label{eq:studyAmuKNeq0SixthPart79}
                        %
                    \end{equation}
                    In \eqref{eq:studyAmuKNeq0SixthPart79}, we note that there is no term associated with $n=1$. Furthermore, we will leave the term $n=0$ for later. Using \eqref{eq:studyAmuKNeq0SixthPart42AnalyticExtHyperGeoT1MinusT}, we have, for every integer $n \geqslant 2$,
                    \begin{equation}
                    \label{eq:studyAmuKNeq0SixthPart80}
                        %
                    \end{equation}
                    The term induced by the first part of the right-hand side of \eqref{eq:studyAmuKNeq0SixthPart80} is given by
                    \begin{equation}
                    \label{eq:studyAmuKNeq0SixthPart81}
                        %
                    \end{equation}
                    For $\mu$ large enough, the term \eqref{eq:studyAmuKNeq0SixthPart81} induces a holomorphic function around $0$, whose derivative at $s=0$ vanishes as $\mu$ goes to infinity. The term associated with the second part of \eqref{eq:studyAmuKNeq0SixthPart80} is
                    \begin{equation}
                    \label{eq:studyAmuKNeq0SixthPart82}
                        %
                    \end{equation}
                    The Euler integral formula yields the estimate
                    \begin{equation}
                    \label{eq:studyAmuKNeq0SixthPart83}
                        %
                    \end{equation}
                    and we can use \eqref{eq:studyAmuKNeq0SixthPart83} to prove that \eqref{eq:studyAmuKNeq0SixthPart82} induces, for $\mu$ large enough, a holomorphic function around $0$, whose derivative at $s=0$ vanishes as $\mu$ goes to infinity. In \eqref{eq:studyAmuKNeq0SixthPart79}, only the term $n=0$ remains. It is given by
                    \begin{equation}
                    \label{eq:studyAmuKNeq0SixthPart84}
                        %
                    \end{equation}
                    It induces a holomorphic function around $0$, whose derivative at $s=0$ is
                    \begin{equation}
                    \label{eq:studyAmuKNeq0SixthPart85}
                        %
                    \end{equation}
                    The term $m=0$ in \eqref{eq:studyAmuKNeq0SixthPart72} is given by
                    \begin{equation}
                    \label{eq:studyAmuKNeq0SixthPart86}
                        %
                    \end{equation}
                    In \eqref{eq:studyAmuKNeq0SixthPart86}, we note that there is no term corresponding to $n=1$. Once again, we will leave the term $n=0$ aside temporarily. Using \eqref{eq:studyAmuKNeq0SixthPart42AnalyticExtHyperGeoT1MinusT}, we have, for any $n \geqslant 2$
                    \begin{equation}
                    \label{eq:studyAmuKNeq0SixthPart87}
                        %
                    \end{equation}
                    For $\mu$ large enough, the term \eqref{eq:studyAmuKNeq0SixthPart88} induces a holomorphic function around $0$, whose derivative at $s=0$ vanishes as $\mu$ goes to infinity. The term induced by the second part of \eqref{eq:studyAmuKNeq0SixthPart87} is
                    \begin{equation}
                    \label{eq:studyAmuKNeq0SixthPart89}
                        %
                    \end{equation}
                    Using the Euler integral formula, we get
                    \begin{equation}
                    \label{eq:studyAmuKNeq0SixthPart90}
                        %
                    \end{equation}
                    for $s$ near $0$, and any integer $n \geqslant 2$. The estimate \eqref{eq:studyAmuKNeq0SixthPart90} proves that, for $\mu$ large enough, the term \eqref{eq:studyAmuKNeq0SixthPart89} induces a holomorphic function around $0$, whose derivative at $s=0$ vanishes as $\mu$ goes to infinity. The term $n=0$ in \eqref{eq:studyAmuKNeq0SixthPart86} is given by
                    \begin{equation}
                    \label{eq:studyAmuKNeq0SixthPart91}
                        %
                    \end{equation}
                    It induced a holomorphic function around $0$, whose derivative at $s=0$ is given by
                    \begin{equation}
                    \label{eq:studyAmuKNeq0SixthPart92}
                        %
                    \end{equation}
                    The proof of the proposition is thus complete.
                \end{proof}

            \subsubsection*{\textbf{Seventh part}}
            \label{AMuKNeq0:part7}
                
                Let us turn our attention to the next term which requires the use of the Ramanujan summation. It is induced by the second term in the right-hand side of \eqref{eq:logExpansionGkLargeOrder}, in lemma \ref{lem:logExpansionGkLargeOrder}.
                
                \begin{proposition}
                \label{prop:studyAMuKNeq0SeventhPart}
                    Assume that $1/\left( 2 \delta \right)$ is not an integer. The function
                    \begin{equation}
                        \label{eq:functionAMukSeventhPart}
                        %
                    \end{equation}
                    which is holomorphic on the half-plane $\Re s > 0$, has a holomorphic continuation near $0$, whose value at $0$ is given by
                    \begin{equation}
                    \label{eq:studyAMukSeventhPartValueAt0}
                        %
                    \end{equation}
                    and whose derivative at $s=0$ equals
                    \begin{equation}
                    \label{eq:studyAMukSeventhPartDerAt0}
                        %
                    \end{equation}
                    as $\mu$ goes to infinity.
                \end{proposition}
                
                \begin{proof}
                    The argument is similar as the one used in the proof of proposition \ref{prop:studyAMuKNeq0SixthPart}.
                \end{proof}
                
                Before we move on to the next term, let us note that the sum of the continuations of the functions from propositions \ref{prop:studyAMuKNeq0FourthPart} and \ref{prop:studyAMuKNeq0SeventhPart} vanish at $s=0$, since we have
                \begin{equation}
                \label{eq:studyAMukFourthAndSeventhPartValueAt0}
                    %
                \end{equation}
                This can be proved using explicit Laurent series expansions. Note that \eqref{eq:studyAMukFourthAndSeventhPartValueAt0} is the sum of \eqref{eq:studyAMukFourthPartValueAt0} and \eqref{eq:studyAMukSeventhPartValueAt0}.

            \subsubsection*{\textbf{Eighth part}}
            \label{AMuKNeq0:part8}
                
                Let us finally deal with the term induced by the last part in the right-hand side of \eqref{eq:logExpansionGkLargeOrder}, in lemma \ref{lem:logExpansionGkLargeOrder}.
                
                \begin{proposition}
                \label{prop:studyAMuKNeq0EighthPart}
                    The function
                    \begin{equation}
                        \label{eq:functionAMukEighthPart}
                        \begin{array}{lll}
                            s & \hspace{-6pt} \mapsto \hspace{-6pt} & \scriptstyle \frac{\sin \left( \pi s \right)}{8\pi} \left( 4 \mu + 1 \right)^{-s+1} \hspace{-2pt} \sum\limits_{k \neq 0} \left( \left \vert k \right \vert^{2\delta} - \frac{1}{4} \right)^{-s} \hspace{-2pt} \frac{1}{\sqrt{4 \pi^2 \left( k + \alpha \right)^2 a^2 + 1/4 + \mu} + \sqrt{4 \pi^2 \left( k - \alpha \right)^2 a^2 + 1/4 + \mu}} \\[1em]
                            
                            && \scriptstyle \qquad \qquad \qquad \qquad \qquad \cdot \left( \frac{1}{4 \pi^2 \left( k + \alpha \right)^2 a^2 + 1/4 + \mu} + \frac{1}{4 \pi^2 \left( k - \alpha \right)^2 a^2 + 1/4 + \mu} \right)
                        \end{array}
                    \end{equation}
                    is holomorphic around $0$, vanishes at $0$, and its derivative at $s=0$ equals
                    \begin{equation}
                    \label{eq:studyAMukEighthPartDerAt0}
                        \begin{array}{c}
                            \frac{1}{2 \pi a} + o \left( 1 \right)
                        \end{array}
                    \end{equation}
                    as $\mu$ goes to infinity.
                \end{proposition}
                
                \begin{proof}
                    The proof follows the same lines as that of proposition \ref{prop:studyAMuKNeq0SixthPart}, the argument being simplified by the fact that \eqref{eq:functionAMukEighthPart} can be directly differentiated at $s=0$, without needing a continuation.
                \end{proof}

            \subsubsection*{\textbf{Ninth part}}
            \label{AMuKNeq0:part9}
                
                The last function we need to look at to complete the study of the terms $A_{\mu,k}$ with $k \neq 0$ is induced by the last part of the right-hand side of \eqref{eq:recallDefHMuK}.
                
                \begin{proposition}
                \label{prop:studyAMuKNeq0NinthPart}
                    The function
                    \begin{equation}
                        \begin{array}{lll}
                            s & \longmapsto & \scriptstyle - i \frac{\sin \left( \pi s \right)}{\pi} \cdot \left( 4 \mu + 1 \right)^{-s+1/2} \sum\limits_{k \neq 0} \left( \left \vert k \right \vert^{2\delta} - \frac{1}{4} \right)^{-s+1} \frac{\partial}{\partial t}_{\vert t = i\sqrt{\frac{1}{4} + \mu}} \log \left \vert g_k \left( t \right) \right \vert
                        \end{array}
                    \end{equation}
                    has a holomorphic continuation near $0$, which vanishes at $s=0$
                \end{proposition}
                
                \begin{proof}
                    This results from propositions \ref{prop:globalStudyAMuKNeq0}, \ref{prop:studyAMuKNeq0FirstPart}, \ref{prop:studyAMuKNeq0SecondPart}, \ref{prop:studyAMuKNeq0ThirdPart}, \ref{prop:studyAMuKNeq0FourthPart}, \ref{prop:studyAMuKNeq0FifthPart}, \ref{prop:studyAMuKNeq0SixthPart}, \ref{prop:studyAMuKNeq0SeventhPart}, and \ref{prop:studyAMuKNeq0EighthPart}.
                \end{proof}

        \subsection{\texorpdfstring{Study of the integrals $M_{\mu,k}$}{Study of the integrals MMuk}}
        \label{subsec:studyMMuk}
        
            Recall that, in subsection \ref{subsec:splittingIntervalIntegration}, the interval of integration used for the terms $I_{\mu,k}$ was split into two parts, yielding a decomposition
            \begin{equation}
                \begin{array}{lll}
                    I_{\mu,k} \left( s \right) & = & L_{\mu,k} \left( s \right) + M_{\mu,k} \left( s \right)
                \end{array}
            \end{equation}
            on the strip $1 < \Re s < 2$, for any integer $k$, with the exception of $k=0$ should $\alpha$ vanish. Having studied the terms $L_{\mu,k}$ in subsections \ref{subsec:generalStudyLMuK} and \ref{subsec:studyAMuK}, we now turn our attention to the terms $M_{\mu,k}$, which were defined in \eqref{eq:intMMuK0} and \eqref{eq:intMMuKNeq0} by
            \begin{equation}
            \label{eq:defMMu0MMuKNeq0}
                \begin{array}{llll}
                    M_{\mu,0} \left( s \right) & = & \frac{\sin \left( \pi s \right)}{\pi} \int_{2 \sqrt{\frac{1}{4} + \mu}}^{+ \infty} \left( t^2 - \left( \frac{1}{4} + \mu \right) \right)^{-s} h_{\mu,0} \left( t \right) \mathrm{d}t \\[1em]
                    
                    M_{\mu,k} \left( s \right) & = & \frac{\sin \left( \pi s \right)}{\pi} \int_{2 \left \vert k \right \vert^{\delta} \sqrt{\frac{1}{4} + \mu}}^{+ \infty} \left( t^2 - \left( \frac{1}{4} + \mu \right) \right)^{-s} h_{\mu,k} \left( t \right) \mathrm{d}t & \text{for } k \neq 0
                \end{array}
            \end{equation}
            the function $h_{\mu,k}$ having been defined in \eqref{eq:functionHMuK} by
            \begin{equation}
            \label{eq:recall2DefHMuK}
                \begin{array}{lll}
                    h_{\mu,k} \left( t \right) & = & i \left( \frac{\partial}{\partial t} \log g_k \right) \left( it \right) - \frac{2it}{\sqrt{4 \mu + 1}} \frac{\partial}{\partial t}_{\vert t = i\sqrt{\frac{1}{4} + \mu}} \log g_k \left( t \right)
                \end{array}
            \end{equation}
            and where the function $g_k$ was introduced in \eqref{eq:functiong} as
            \begin{equation}
            \label{eq:recallDefGK}
                \begin{array}{lll}
                    \scriptstyle g_k \left( \nu \right) & = & \scriptstyle 1 + 4 \pi \alpha a + 2 \pi \left \vert k + \alpha \right \vert a \frac{K_{i\nu}' \left( 2 \pi \left \vert k + \alpha \right \vert a \right)}{K_{i\nu} \left( 2 \pi \left \vert k + \alpha \right \vert a \right)} + 2 \pi \left \vert k - \alpha \right \vert a \frac{K_{i\nu}' \left( 2 \pi \left \vert k - \alpha \right \vert a \right)}{K_{i\nu} \left( 2 \pi \left \vert k - \alpha \right \vert a \right)}.
                \end{array}
            \end{equation}
            Let us now study these terms, separating the cases $k = 0$ and $k \neq 0$.

            \subsubsection{The case $k = 0$}
            \label{subsubsec:studyMMuK0}
                
                We first need to study the term $M_{\mu,0}$ for any $\mu \geqslant 0$. Similarly to what was done in \cite[Sec. 3.5]{dutour:cuspDirichlet} and in \cite[Sec. 6.3]{freixas-vonPippich:RROrbifold}, we will begin by breaking apart $M_{\mu,0}$ into two parts. Let us assume that we have $\alpha \neq 0$.
                
                \begin{definition}
                \label{def:MTildeMu0RMu0}
                    For any real number $\mu \geqslant 0$, we set
                    \begin{equation}
                    \label{eq:defMTildeMu0}
                        \begin{array}{lll}
                            \widetilde{M}_{\mu, 0} \left( s \right) & = & i \frac{\sin \left( \pi s \right)}{\pi} \int_{2 \sqrt{\frac{1}{4} + \mu}}^{+ \infty} \left( t^2 - \left( \frac{1}{4} + \mu \right) \right)^{-s} \left( \frac{\partial}{\partial t} \log g_0 \right) \left( it \right) \mathrm{d}t
                        \end{array}
                    \end{equation}
                    on the strip $1 < \Re s < 2$, where we also set
                    \begin{equation}
                    \label{eq:defRMu0}
                        \begin{array}{lll}
                            R_{\mu,0} \left( s \right) & = & \scriptstyle - \frac{2i}{\sqrt{4 \mu + 1}} \cdot \frac{\sin \left( \pi s \right)}{\pi} \left( \int_{2 \sqrt{\frac{1}{4} + \mu}}^{+ \infty} t \left( t^2 - \left( \frac{1}{4} + \mu \right) \right)^{-s} \mathrm{d}t \right) \frac{\partial}{\partial t}_{\vert t = i\sqrt{\frac{1}{4} + \mu}} \log g_0 \left( t \right).
                        \end{array}
                    \end{equation}
                \end{definition}
                
                We have $M_{\mu,0} \left( s \right) = \widetilde{M}_{\mu,0} \left( s \right) + R_{\mu,0} \left( s \right)$ on the strip $1 < \Re s < 2$. Furthermore, note that the integral in \eqref{eq:defRMu0} can be computed exactly, though it was not done in definition \ref{def:MTildeMu0RMu0} to make the splitting of $M_{\mu,0}$ easier to see.
                
                \begin{proposition}
                \label{prop:studyRMuK0}
                    Assume we have $\mu \geqslant 0$. The function
                    \begin{equation}
                        \begin{array}{lll}
                            s & \longmapsto & R_{\mu,0} \left( s \right)
                        \end{array}
                    \end{equation}
                    has a holomorphic continuation near $0$, whose derivative at $s=0$ is given by
                    \begin{equation}
                    \label{eq:RMu0DerAt0}
                        \begin{array}{lll}
                            R_{\mu,0}' \left( 0 \right) & = & \frac{3}{2} i \sqrt{\frac{1}{4} + \mu} \; \frac{\partial}{\partial t}_{\vert t = i\sqrt{\frac{1}{4} + \mu}} \log \left \vert g_0 \left( t \right) \right \vert.
                        \end{array}
                    \end{equation}
                    In particular, we have
                    \begin{equation}
                    \label{eq:RMu0DerAt0AsMuGoesTo0}
                        \begin{array}{lll}
                            \lim\limits_{\; \; \mu \rightarrow 0^+} R_{\mu,0}' \left( 0 \right) & = & \frac{3}{4} i \frac{\partial}{\partial t}_{\vert t = i/2} \log \left \vert g_0 \left( t \right) \right \vert.
                        \end{array}
                    \end{equation}
                \end{proposition}
                
                \begin{proof}
                    This is a consequence of the fact that we have
                    \begin{equation}
                        \begin{array}{lll}
                            \int_{2 \sqrt{\frac{1}{4} + \mu}}^{+ \infty} t \left( t^2 - \left( \frac{1}{4} + \mu \right) \right)^{-s} \mathrm{d}t & = & - \frac{1}{2} \cdot \frac{1}{1-s} \cdot 3^{-s+1} \left( \frac{1}{4} + \mu \right)^{-s+1}
                        \end{array}
                    \end{equation}
                    for any complex number $s$ with $\Re s > 1$.
                \end{proof}
                
                Let us now deal with the term \eqref{eq:defMTildeMu0}. To that effect, recall that we have
                \begin{equation}
                \label{eq:recallAsymptG0LargeT}
                    \begin{array}{lll}
                        \log \left \vert g_0 \left( i t \right) \right \vert & = & \log 2 + \frac{1}{2} \log \left( 4 \pi^2 \alpha^2 a^2 + t^2 \right) - \frac{2 \pi \alpha a}{\sqrt{4 \pi^2 \alpha^2 a^2 + t^2}} \\[0.5em]
                        
                        && \qquad \qquad - \frac{1}{2} \cdot \frac{1}{\sqrt{4 \pi^2 \alpha^2 a^2 + t^2}} \cdot \frac{t^2}{4 \pi^2 \alpha^2 a^2 + t^2} + \eta_{2,3} \left( t \right),
                    \end{array}
                \end{equation}
                for $t$ large enough, as stated in lemma \eqref{lem:asymptExpLargeOrderFunctionLogGK0}, with the estimate
                \begin{equation}
                    \begin{array}{lll}
                        \left \vert \eta_{2,3} \left( t \right) \right \vert & \leqslant & \frac{C}{4 \pi^2 \alpha^2 a^2 + t^2},
                    \end{array}
                \end{equation}
                with $C > 0$ being a constant.
                
                \begin{proposition}
                \label{prop:studyMTildeMuK0}
                    Assume we have $\mu \geqslant 0$. The function
                    \begin{equation}
                        \begin{array}{lll}
                            s & \longmapsto & \widetilde{M}_{\mu,0} \left( s \right)
                        \end{array}
                    \end{equation}
                    has a holomorphic continuation near $0$, whose derivative at $s=0$ is given by
                    \begin{equation}
                    \label{eq:MTildeMu0DerAt0}
                        \begin{array}{lll}
                            \widetilde{M}_{\mu,0}' \left( 0 \right) & \hspace{-3pt} = \hspace{-3pt} & \left \{ \hspace{-1pt} \begin{array}{ll} \scriptstyle - \log \left \vert g_0 \left( 2 i \sqrt{\frac{1}{4} + \mu} \right) \right \vert + \log 2 + o \left( 1 \right) & \text{as } \mu \rightarrow + \infty \\[0.5em] \scriptstyle - \log \left \vert g_0 \left( i \right) \right \vert + \log 2 + o \left( 1 \right) & \text{as } a \rightarrow + \infty, \text{ with } \mu = 0. \end{array} \right.
                        \end{array}
                    \end{equation}
                \end{proposition}
                
                \begin{proof}
                    For any complex number $s$ with $\Re s > 1$, we have
                    \begin{equation}
                    \label{eq:MTildeMuK001}
                        \begin{array}{lll}
                            \scriptstyle \widetilde{M}_{\mu, 0} \left( s \right) & = & \scriptstyle - \frac{\sin \left( \pi s \right)}{\pi} 3^{-s} \left( \frac{1}{4} + \mu \right)^{-s} \left( \log \left \vert g_0 \left( 2 i \sqrt{\frac{1}{4} + \mu} \right) \right \vert - \log 2 \right. \\[0.5em]
                            
                            && \scriptstyle \qquad \qquad \qquad \qquad \left. - \frac{1}{2} \log \left( 4 \pi^2 \alpha^2 a^2 + 4 \mu + 1 \right) + \frac{2 \pi \alpha a}{\sqrt{4 \pi^2 \alpha^2 a^2 + 4 \mu + 1}} \right) \\[0.5em]
                            
                            && \scriptstyle - 2s \frac{\sin \left( \pi s \right)}{\pi} \int_{2 \sqrt{\frac{1}{4} + \mu}}^{+ \infty} t \left( t^2 - \left( \frac{1}{4} + \mu \right) \right)^{-s-1} \left( \log \left \vert g_0 \left( it \right) \right \vert - \log 2 \right. \\[0.5em]
                            
                            && \scriptstyle \qquad \qquad \qquad \qquad \left. - \frac{1}{2} \log \left( 4 \pi^2 \alpha^2 a^2 + t^2 \right) + \frac{2 \pi \alpha a}{\sqrt{4 \pi^2 \alpha^2 a^2 + t^2}} \right) \mathrm{d}t \\[1em]
                            
                            && \scriptstyle + \frac{1}{2} \cdot \frac{\sin \left( \pi s \right)}{\pi} \int_{2 \sqrt{\frac{1}{4} + \mu}}^{+ \infty} \left( t^2 - \left( \frac{1}{4} + \mu \right) \right)^{-s} \frac{\partial}{\partial t} \log \left( 4 \pi^2 \alpha^2 a^2 + t^2 \right) \mathrm{d}t \\[1em]
                            
                            && \scriptstyle - 2 \pi \alpha a \frac{\sin \left( \pi s \right)}{\pi} \int_{2 \sqrt{\frac{1}{4} + \mu}}^{+ \infty} \left( t^2 - \left( \frac{1}{4} + \mu \right) \right)^{-s} \frac{\partial}{\partial t} \frac{1}{\sqrt{4 \pi^2 \alpha^2 a^2 + t^2}} \mathrm{d}t.
                        \end{array}
                    \end{equation}
                    The second term on the right-hand side of \eqref{eq:MTildeMuK001} induces a holomorphic function around $0$ whose derivative at $s=0$ vanishes, which we see using \eqref{eq:recallAsymptG0LargeT}. The first term of the right-hand side of \eqref{eq:MTildeMuK001} also yields a holomorphic function near $0$, whose derivative at $s=0$ equals
                    \begin{equation}
                    \label{eq:MTildeMuK002}
                        \begin{array}{lll}
                            \multicolumn{3}{l}{\scriptstyle - \log \left \vert g_0 \left( 2 i \sqrt{\frac{1}{4} + \mu} \right) \right \vert + \log 2 + \frac{1}{2} \log \left( 4 \pi^2 \alpha^2 a^2 + 4 \mu + 1 \right) - \frac{2 \pi \alpha a}{\sqrt{4 \pi^2 \alpha^2 a^2 + 4 \mu + 1}}} \\[0.5em]
                            
                            \qquad \qquad \qquad \qquad \qquad \qquad & = & \scriptstyle - \log \left \vert g_0 \left( 2 i \sqrt{\frac{1}{4} + \mu} \right) \right \vert + 2 \log 2 + \frac{1}{2} \log \mu + o \left( 1 \right)
                        \end{array}
                    \end{equation}
                    as $\mu$ goes to infinity. Furthermore, the same derivative, taken with $\mu = 0$, equals
                    \begin{equation}
                    \label{eq:MTildeMuK003}
                        \begin{array}{lll}
                            \multicolumn{3}{l}{\scriptstyle - \log \left \vert g_0 \left( i \right) \right \vert + \log 2 + \frac{1}{2} \log \left( 4 \pi^2 \alpha^2 a^2 + 1 \right) - \frac{2 \pi \alpha a}{\sqrt{4 \pi^2 \alpha^2 a^2 + 1}}} \\[0.5em]
                            
                            \qquad \qquad \qquad \qquad \qquad \qquad \qquad & = & \scriptstyle - \log \left \vert g_0 \left( i \right) \right \vert + \log a + \log \left( 4 \pi \alpha \right) - 1 + o \left( 1 \right)
                        \end{array}
                    \end{equation}
                    as $a$ goes to infinity. The third term of the right-hand side of \eqref{eq:MTildeMuK001}, given by
                    \begin{equation}
                    \label{eq:MTildeMuK004}
                        \begin{array}{lll}
                            s & \longmapsto & \scriptstyle \frac{1}{2} \cdot \frac{\sin \left( \pi s \right)}{\pi} \int_{2 \sqrt{\frac{1}{4} + \mu}}^{+ \infty} \left( t^2 - \left( \frac{1}{4} + \mu \right) \right)^{-s} \frac{\partial}{\partial t} \log \left( 4 \pi^2 \alpha^2 a^2 + t^2 \right) \mathrm{d}t
                        \end{array}
                    \end{equation}
                    has a holomorphic continuation near $0$, whose derivative at $s=0$ is given by
                    \begin{equation}
                    \label{eq:MTildeMuK005}
                        \begin{array}{lll}
                            \multicolumn{3}{l}{\scriptstyle \frac{\partial}{\partial s}_{\vert s = 0} \left[ \frac{1}{2} \cdot \frac{\sin \left( \pi s \right)}{\pi} \int_{2 \sqrt{\frac{1}{4} + \mu}}^{+ \infty} \left( t^2 - \left( \frac{1}{4} + \mu \right) \right)^{-s} \frac{\partial}{\partial t} \log \left( 4 \pi^2 \alpha^2 a^2 + t^2 \right) \mathrm{d}t \right]} \\[1em]
                            
                            \qquad \qquad \qquad \qquad \qquad & = & \left \{ \begin{array}{ll} \scriptstyle - \frac{1}{2} \log \mu - \log 2 + o \left( 1 \right) & \scriptstyle \text{as } \mu \rightarrow + \infty \\[0.5em] \scriptstyle - \log \left( 2 \pi \alpha a \right) + o \left( 1 \right) & \scriptstyle \text{as } a \rightarrow + \infty, \text{ with } \mu = 0 \end{array} \right.
                        \end{array}
                    \end{equation}
                    using \cite[Prop. 3.71]{dutour:cuspDirichlet}. The last term of the right-hand side of \eqref{eq:MTildeMuK001}, given by
                    \begin{equation}
                        \begin{array}{lll}
                            s & \longmapsto & \scriptstyle - 2 \pi \alpha a \frac{\sin \left( \pi s \right)}{\pi} \int_{2 \sqrt{\frac{1}{4} + \mu}}^{+ \infty} \left( t^2 - \left( \frac{1}{4} + \mu \right) \right)^{-s} \frac{\partial}{\partial t} \frac{1}{\sqrt{4 \pi^2 \alpha^2 a^2 + t^2}} \mathrm{d}t
                        \end{array}
                    \end{equation}
                    has a holomorphic continuation near $0$, which can be proved using computations similar to those from \cite[Prop. 3.70]{dutour:cuspDirichlet}, and its derivative at $s=0$ equals
                    \begin{equation}
                        \begin{array}{lll}
                            \multicolumn{3}{l}{\scriptstyle \frac{\partial}{\partial s}_{\vert s = 0} \left[ - 2 \pi \alpha a \frac{\sin \left( \pi s \right)}{\pi} \int_{2 \sqrt{\frac{1}{4} + \mu}}^{+ \infty} \left( t^2 - \left( \frac{1}{4} + \mu \right) \right)^{-s} \frac{\partial}{\partial t} \frac{1}{\sqrt{4 \pi^2 \alpha^2 a^2 + t^2}} \mathrm{d}t \right]} \\[1em]
                            
                            \qquad \qquad \qquad \qquad \qquad \qquad & = & \left \{ \begin{array}{ll} \scriptstyle o \left( 1 \right) & \scriptstyle \text{as } \mu \rightarrow + \infty \\[0.5em] \scriptstyle o \left( 1 \right) & \scriptstyle \text{as } a \rightarrow + \infty, \text{ with } \mu = 0. \end{array} \right.
                        \end{array}
                    \end{equation}
                    This concludes the proof of the proposition.
                \end{proof}

            \subsubsection{The case $k \neq 0$}
            \label{subsubsec:studyMMukSplit}
                
                Let us now study the sum over all integers $k \neq 0$ of the terms $M_{\mu,k}$, for any $\mu \geqslant 0$. As in \cite[Sec. 3.5]{dutour:cuspDirichlet} and \cite[Sec. 6.3]{freixas-vonPippich:RROrbifold}, we will begin by breaking apart $M_{\mu,k}$ into two parts. Note that $\alpha$ may vanish in this paragraph.
                
                \begin{definition}
                    For any real number $\mu \geqslant 0$ and any integer $k \neq 0$, we set
                    \begin{equation}
                    \label{eq:defMTildeMuKNeq0}
                        \begin{array}{lll}
                            \widetilde{M}_{\mu, k} \left( s \right) & = & i \frac{\sin \left( \pi s \right)}{\pi} \int_{2 \left \vert k \right \vert^{\delta} \sqrt{\frac{1}{4} + \mu}}^{+ \infty} \left( t^2 - \left( \frac{1}{4} + \mu \right) \right)^{-s} \left( \frac{\partial}{\partial t} \log g_k \right) \left( it \right) \mathrm{d}t
                        \end{array}
                    \end{equation}
                    on the strip $1 < \Re s < 2$, where we also set
                    \begin{equation}
                    \label{eq:defRMuKNeq0}
                        \begin{array}{lll}
                            R_{\mu,k} \left( s \right) & \hspace{-6pt} = \hspace{-6pt} & \scriptstyle - \frac{2i}{\sqrt{4 \mu + 1}} \cdot \frac{\sin \left( \pi s \right)}{\pi} \left( \int_{2 \left \vert k \right \vert^{\delta} \sqrt{\frac{1}{4} + \mu}}^{+ \infty} t \left( t^2 - \left( \frac{1}{4} + \mu \right) \right)^{-s} \mathrm{d}t \right) \frac{\partial}{\partial t}_{\vert t = i\sqrt{\frac{1}{4} + \mu}} \log g_k \left( t \right) \\[1em]
                            
                            & \hspace{-6pt} = \hspace{-6pt} & \scriptstyle \frac{i}{1-s} \cdot \frac{\sin \left( \pi s \right)}{\pi} \cdot \left( 4 \mu + 1 \right)^{-s+1/2} \left( \left \vert k \right \vert^{2 \delta} - \frac{1}{4} \right)^{-s+1} \frac{\partial}{\partial t}_{\vert t = i\sqrt{\frac{1}{4} + \mu}} \log g_k \left( t \right).
                        \end{array}
                    \end{equation}
                \end{definition}
                
                We have $M_{\mu,k} \left( s \right) = \widetilde{M}_{\mu,k} \left( s \right) + R_{\mu,k} \left( s \right)$ on the strip $1 < \Re s < 2$.
                
                \begin{proposition}
                \label{prop:studyRMuKNeq0}
                    Assume we have $\mu \geqslant 0$. The function
                    \begin{equation}
                        \begin{array}{lll}
                            s & \longmapsto & \sum\limits_{k \neq 0} R_{\mu,k} \left( s \right)
                        \end{array}
                    \end{equation}
                    has a holomorphic continuation near $0$, which vanishes at $0$, and whose derivative at $s=0$, taken with $\mu = 0$, equals
                    \begin{equation}
                    \label{eq:RMuKNeq0DerAt0AAsympt}
                        \begin{array}{lll}
                            \frac{\partial}{\partial s}_{\vert s = 0} \sum\limits_{k \neq 0} R_{0,k} \left( s \right) & = & o \left( 1 \right)
                        \end{array}
                    \end{equation}
                    as $a$ goes to infinity.
                \end{proposition}
                
                \begin{proof}
                    Let us note that the function we study here, as evidenced in \eqref{eq:defRMuKNeq0}, is related to the function from proposition \ref{prop:studyAMuKNeq0NinthPart} by a factor $1/\left( s-1 \right)$, thus yielding the holomorphic continuation and its vanishing at $0$. Finally, the computation of the derivative at $s=0$ as $a$ goes to infinity stems from the special values \eqref{eq:specialValueHalfModBesselFunction}.
                \end{proof}
                
                We will now deal with the sum over $k \neq 0$ of the terms \eqref{eq:defMTildeMuKNeq0}. Let us first recall what we obtained in lemma \ref{lem:logExpansionGkLargeOrder} and in its proof. We have
                \begin{equation}
                \label{eq:recallLogExpansionGkLargeOrder}

                \end{equation}
                for either $t \geqslant A$ or $a \geqslant B$ large enough. The remainder $\eta_{2,3}$ is then given by
                \begin{equation}
                \label{eq:recallRemainderLargeOrderAsymptLogGKNeq0}
                    %
                \end{equation}
                and in particular satisfies
                \begin{equation}
                \label{eq:recallEstimateRemainderAsymptLargeOrderGKNeq0}
                    %
                \end{equation}
                for some constant $C > 0$ and either $t \geqslant A$ large enough or $\left \vert k \right \vert \geqslant K_0$ large enough.

                \subsubsection*{\textbf{First part}}
                \label{MTildeMuKNeq0:part1}
                
                    We begin with the term induced by the remainder $\eta_{2,3}$.
                    
                    \begin{proposition}
                    \label{prop:studyMTildeMuKNeq0FirstPart}
                        The function
                        \begin{equation}
                        \label{eq:functionStudyMTildeMuKNeq0FirstPart}
                            %
                        \end{equation}
                        is holomorphic around $0$, and its derivative at $s=0$ equals
                        \begin{equation}
                            \left \{ %
 \right.
                        \end{equation}
                    \end{proposition}
                    
                    \begin{proof}
                        We first note that the holomorphy around $0$ of \eqref{eq:functionStudyMTildeMuKNeq0FirstPart} is a consequence of estimate \eqref{eq:recallEstimateRemainderAsymptLargeOrderGKNeq0}, which holds for any $\mu \geqslant 0$ and any $a > 0$ if $\left \vert k \right \vert$ is large enough, and an integration by parts. Taking the derivative at $s=0$, only the integrated term remains, which is given by
                        \begin{equation}
                            %
                        \end{equation}
                        The $\mu$-asymptotic computation of the derivative at $s=0$ is also a consequence of \eqref{eq:recallEstimateRemainderAsymptLargeOrderGKNeq0}. Let us now obtain the $a$-asymptotic expansion of the same derivative, taken with $\mu = 0$. The first term induced by \eqref{eq:recallRemainderLargeOrderAsymptLogGKNeq0}, given by
                        \begin{equation}
                        \label{eq:studyMTildeMuKNeq0FirstPart01}
                            %
                        \end{equation}
                        is holomorphic around $0$, and its derivative at $s=0$ vanishes as $a$ goes to infinity, using \eqref{eq:estimateRemainderEta22}. The term induced by the rest of \eqref{eq:recallRemainderLargeOrderAsymptLogGKNeq0} is given by
                        \begin{equation}
                        \label{eq:studyMTildeMuKNeq0FirstPart02}
                            %
                        \end{equation}
                        It is now seen that the study of \eqref{eq:studyMTildeMuKNeq0FirstPart02} is reduced to that of
                        \begin{equation}
                        \label{eq:studyMTildeMuKNeq0FirstPart03}
                            %
                        \end{equation}
                        by using Newton's binomial formula, since every other part it yields induces a holomorphic function around $0$, whose derivatives at $s=0$ vanish as $a$ goes to infinity. Let us therefore study \eqref{eq:studyMTildeMuKNeq0FirstPart03}. The dominated convergence theorem yields
                        \begin{equation}
                        \label{eq:studyMTildeMuKNeq0FirstPart04}
                            %
                        \end{equation}
                        the last equality being a consequence of the generating function formula
                        \begin{equation}
                        \label{eq:studyMTildeMuKNeq0FirstPart05}
                            %
                        \end{equation}
                        for any $0 \leqslant x < 1$, where $\psi$ is the Digamma function, \textit{i.e.} the logarithmic derivative of the Gamma function. This completes the proof of the proposition.
                    \end{proof}

                \subsubsection*{\textbf{Second part}}
                \label{MTildeMuKNeq0:part2}
                
                    Let us now study of the term induced by the logarithm from \eqref{eq:recallLogExpansionGkLargeOrder}.
                    
                    \begin{proposition}
                    \label{prop:studyMTildeMuKNeq0SecondPart}
                        Assume that we have $\alpha \neq 0$. The function
                        \begin{equation}
                        \label{eq:functionStudyMTildeMuKNeq0SecondPartAlphaNeq0}
                            %
                        \end{equation}
                        has a holomorphic continuation near $0$. Its derivative at $s=0$ equals
                        \begin{equation}
                        \label{eq:functionStudyMTildeMuKNeq0SecondPartAlphaNeq0DerAt0AsMuInfinity}
                            %
                        \end{equation}
                        as $\mu$ goes to infinity. The same derivative, taken with $\mu = 0$ further satisfies
                        \begin{equation}
                        \label{eq:functionStudyMTildeMuKNeq0SecondPartAlphaNeq0DerAt0AsAInfinity}
                            %
                        \end{equation}
                    \end{proposition}
                    
                    \begin{proof}
                        Let us first use the binomial formula. We have
                        \begin{equation}
                        \label{eq:studyMTildeMuKNeq0SecondPart01}
                            %
                        \end{equation}
                        and we can now study the sum over $j \geqslant 1$ by computing the derivative. We have
                        \begin{equation}
                        \label{eq:studyMTildeMuKNeq0SecondPart02}
                            %
                        \end{equation}
                        for any $t$ in the interval of integration and any integer $k \neq 0$. We have
                        \begin{equation}
                        \label{eq:studyMTildeMuKNeq0SecondPart03}
                            %
                        \end{equation}
                        Therefore, using \eqref{eq:studyMTildeMuKNeq0SecondPart03}, the term induced by the sum over $j \geqslant 1$, given by
                        \begin{equation}
                        \label{eq:studyMTildeMuKNeq0SecondPart04}
                            %
                        \end{equation}
                        induces a holomorphic function around $0$, whose derivative at $s=0$ vanishes, because of $\left( s \right)_j$ for $j \geqslant 1$. The term corresponding to $j=0$ is given by
                        \begin{equation}
                        \label{eq:studyMTildeMuKNeq0SecondPart05}
                            %
                        \end{equation}
                        Let us perform an integration by parts. We have
                        \begin{equation}
                        \label{eq:studyMTildeMuKNeq0SecondPart06}
                            %
                        \end{equation}
                        The term induced by the first part of the right-hand side of \eqref{eq:studyMTildeMuKNeq0SecondPart06} is given by
                        \begin{equation}
                        \label{eq:studyMTildeMuKNeq0SecondPart07}
                            %
                        \end{equation}
                        Using the computations performed in the proof of proposition \ref{prop:studyAMuKNeq0ThirdPart}, the function induced by \eqref{eq:studyMTildeMuKNeq0SecondPart07} has a holomorphic continuation near $0$. Its derivative at $s=0$ needs not be computed as $\mu$ goes to infinity, as it is in \eqref{eq:functionStudyMTildeMuKNeq0SecondPartAlphaNeq0DerAt0AsMuInfinity}. However, we need to study it as $a$ goes to infinity, with $\mu = 0$. For any integer $k \geqslant 1$, we have
                        \begin{equation}
                        \label{eq:studyMTildeMuKNeq0SecondPart08}
                            %
                        \end{equation}
                        We now note that we have
                        \begin{equation}
                        \label{eq:studyMTildeMuKNeq0SecondPart09}
                            %
                        \end{equation}
                        on the interval of integration. Thus, the function induced from \eqref{eq:studyMTildeMuKNeq0SecondPart07} by the last part of the right-hand side of \eqref{eq:studyMTildeMuKNeq0SecondPart08} is holomorphic around $0$, and its derivative at $s=0$ vanishes as $a$ goes to infinity. The remaining term, given by
                        \begin{equation}
                        \label{eq:studyMTildeMuKNeq0SecondPart10}
                            %
                        \end{equation}
                        induces a holomorphic function around $0$, whose derivative at $s=0$ equals
                        \begin{equation}
                        \label{eq:studyMTildeMuKNeq0SecondPart11}
                            %
                        \end{equation}
                        We now move on to the term induced by the last part of \eqref{eq:studyMTildeMuKNeq0SecondPart06}, given by
                        \begin{equation}
                        \label{eq:studyMTildeMuKNeq0SecondPart12}
                            %
                        \end{equation}
                        Let us begin with the second term of the right-hand side of \eqref{eq:studyMTildeMuKNeq0SecondPart12}. We have
                        \begin{equation}
                        \label{eq:studyMTildeMuKNeq0SecondPart13}
                            %
                        \end{equation}
                        The term induced by the constant part of \eqref{eq:studyMTildeMuKNeq0SecondPart13}, given by
                        \begin{equation}
                        \label{eq:studyMTildeMuKNeq0SecondPart14}
                            %
                        \end{equation}
                        induces a holomorphic function around $0$, whose derivative at $s=0$ equals
                        \begin{equation}
                        \label{eq:studyMTildeMuKNeq0SecondPart15}
                            %
                        \end{equation}
                        We must now deal with
                        \begin{equation}
                        \label{eq:studyMTildeMuKNeq0SecondPart16}
                            %
                        \end{equation}
                        Let us compute the integrals present in \eqref{eq:studyMTildeMuKNeq0SecondPart16}. We have
                        \begin{equation}
                        \label{eq:studyMTildeMuKNeq0SecondPart17}
                            %
                        \end{equation}
                        Hence, we have to deal with
                        \begin{equation}
                        \label{eq:studyMTildeMuKNeq0SecondPart18}
                            %
                        \end{equation}
                        Let us begin with the sum over $n \geqslant 3$. Using the contiguous functions relation
                        \begin{equation}
                        \label{eq:studyMTildeMuKNeq0SecondPart19ContiguousFunctions}
                            %
                        \end{equation}
                        for $t \in \left[ 0, 1 \right[$, which one may check by direct computation using the power series expansions of hypergeometric functions, we have
                        \begin{equation}
                        \label{eq:studyMTildeMuKNeq0SecondPart20}
                            %
                        \end{equation}
                        for any $t \in \left[ 0,1 \right[$. The first part of the right-hand side of \eqref{eq:studyMTildeMuKNeq0SecondPart20} yields the term
                        \begin{equation}
                        \label{eq:studyMTildeMuKNeq0SecondPart21}
                            %
                        \end{equation}
                        For any integer $n \geqslant 3$ and any $\left \vert s \right \vert \leqslant \varepsilon$ small enough, we have
                        \begin{equation}
                        \label{eq:studyMTildeMuKNeq0SecondPart22}
                            %
                        \end{equation}
                        Having $4 \alpha < \left( 1 + \alpha \right)^2$, the term \eqref{eq:studyMTildeMuKNeq0SecondPart21} induces a holomorphic function near $0$, whose derivative at $s=0$ vanishes, because of the extra factor $s$. Going back to decomposition \eqref{eq:studyMTildeMuKNeq0SecondPart20}, the second part of the right-hand side induces the term
                        \begin{equation}
                        \label{eq:studyMTildeMuKNeq0SecondPart23}
                            %
                        \end{equation}
                        We now note that, for any fixed real number $\varepsilon > 0$ small enough, the hypergeometric function $t \longmapsto F \left( n,\; 1; \; -\varepsilon+n+2; \; t \right)$ is increasing on $\left[0,1 \right[$, and we have
                        \begin{equation}
                        \label{eq:studyMTildeMuKNeq0SecondPart24}
                            %
                        \end{equation}
                        Therefore, we have, for any complex number $s$ with $\left \vert s \right \vert \leqslant \varepsilon$ small enough
                        \begin{equation}
                        \label{eq:studyMTildeMuKNeq0SecondPart25}
                            %
                        \end{equation}
                        Using an estimate similar to \eqref{eq:studyMTildeMuKNeq0SecondPart22} on the right-hand side of \eqref{eq:studyMTildeMuKNeq0SecondPart25} for $n \geqslant 3$, the term \eqref{eq:studyMTildeMuKNeq0SecondPart23} induces a holomorphic function around $0$, whose derivative at $s=0$ vanishes. The term corresponding to $n=2$ in \eqref{eq:studyMTildeMuKNeq0SecondPart18} is
                        \begin{equation}
                        \label{eq:studyMTildeMuKNeq0SecondPart26}
                            %
                        \end{equation}
                        Using the Euler integral formula, we have
                        \begin{equation}
                        \label{eq:studyMTildeMuKNeq0SecondPart27}
                            %
                        \end{equation}
                        for any $t \in \left[ 0,1 \right[$. We can plug this estimate into \eqref{eq:studyMTildeMuKNeq0SecondPart26} to get
                        \begin{equation}
                        \label{eq:studyMTildeMuKNeq0SecondPart28}
                            %
                        \end{equation}
                        Therefore, the term \eqref{eq:studyMTildeMuKNeq0SecondPart26} induces a holomorphic function around $0$, whose derivative at $s=0$ vanishes. The term corresponding to $n=1$ in \eqref{eq:studyMTildeMuKNeq0SecondPart18} is
                        \begin{equation}
                        \label{eq:studyMTildeMuKNeq0SecondPart29}
                            %
                        \end{equation}
                        Let us break down the hypergeometric function in \eqref{eq:studyMTildeMuKNeq0SecondPart29}, using \eqref{eq:studyAmuKNeq0SixthPart42AnalyticExtHyperGeoT1MinusT}. We have
                        \begin{equation}
                        \label{eq:studyMTildeMuKNeq0SecondPart30}
                            %
                        \end{equation}
                        The term induced by the first part of the right-hand side of \eqref{eq:studyMTildeMuKNeq0SecondPart30} is given by
                        \begin{equation}
                        \label{eq:studyMTildeMuKNeq0SecondPart31}
                            %
                        \end{equation}
                        In order to obtain \eqref{eq:studyMTildeMuKNeq0SecondPart31}, we had to use the special value $\Gamma \left( 3/2 \right) = \sqrt{\pi}/2$ and the reflection formula $\Gamma \left( 1-s \right) \Gamma \left( s \right) = \pi / \sin \left( \pi s \right)$. The term \eqref{eq:studyMTildeMuKNeq0SecondPart31} induces a holomorphic function around $0$, whose derivative at $s=0$ equals
                        \begin{equation}
                        \label{eq:studyMTildeMuKNeq0SecondPart32}
                            %
                        \end{equation}
                        where $\psi$ denotes the Digamma function, \textit{i.e.} the logarithmic derivative of the Gamma function. Note that, in order to get \eqref{eq:studyMTildeMuKNeq0SecondPart32}, one uses Laurent series expansions, and the special value $\zeta_H \left( 1 + \alpha, 2 \right) = \psi' \left( 1 + \alpha \right)$ of the Hurwitz zeta function. The term induced by the second part of the right-hand side of \eqref{eq:studyMTildeMuKNeq0SecondPart30} is
                        \begin{equation}
                        \label{eq:studyMTildeMuKNeq0SecondPart33}
                            %
                        \end{equation}
                        For any integer $k \geqslant 1$ and any real number $t \in \left[ 0, 1 \right[$, we have
                        \begin{equation}
                        \label{eq:studyMTildeMuKNeq0SecondPart34}
                            %
                        \end{equation}
                        We will now extract the first few terms of the hypergeometric function. We have
                        \begin{equation}
                        \label{eq:studyMTildeMuKNeq0SecondPart36}
                            %
                        \end{equation}
                        for any $t \in \left[ 0, 1 \right[$, where the last part of the right-hand side of \eqref{eq:studyMTildeMuKNeq0SecondPart36} is a generalized hypergeometric function, whose definition in terms of a power series can be found in \cite[Prop.-Def. C.36]{dutour:cuspDirichlet}. The term induced by the second part of the right-hand side of \eqref{eq:studyMTildeMuKNeq0SecondPart36} is
                        \begin{equation}
                        \label{eq:studyMTildeMuKNeq0SecondPart37}
                            %
                        \end{equation}
                        Now, note that the function $t \longmapsto F \left( -s, \; -s, \; 1; \; 1-s; \; 2; \; t \right)$ is uniformly bounded on $\left[ 0,1 \right[$ for $s$ around $0$, since we have $\Re \left( 1 - s - \left( -s \right) - \left( -s \right) \right) = 1 + \Re s > 0$ assuming we have $\Re s > -1$. This proves that \eqref{eq:studyMTildeMuKNeq0SecondPart37} induces a holomorphic function around $0$, whose derivative at $s=0$ vanishes. The term associated with the first part of the right-hand side of \eqref{eq:studyMTildeMuKNeq0SecondPart36} is
                        \begin{equation}
                        \label{eq:studyMTildeMuKNeq0SecondPart38}

                        \end{equation}
                        Using the binomial formula, we have
                        \begin{equation}
                        \label{eq:studyMTildeMuKNeq0SecondPart39}
                            %
                        \end{equation}
                        The term induced by the sum over $j \geqslant 1$ yields a holomorphic function around $0$, whose derivative at $s=0$ vanishes, because of the Pochhammer symbol. Therefore, we need only study the term $j = 0$, which is given by
                        \begin{equation}
                        \label{eq:studyMTildeMuKNeq0SecondPart40}
                            %
                        \end{equation}
                        The second term in the right-hand side of \eqref{eq:studyMTildeMuKNeq0SecondPart40} induces a holomorphic function around $0$, whose derivative at $s=0$ equals
                        \begin{equation}
                        \label{eq:studyMTildeMuKNeq0SecondPart41}
                            %
                        \end{equation}
                        The first term of \eqref{eq:studyMTildeMuKNeq0SecondPart40} also induces a holomorphic function around $0$, and we compute its derivative using Laurent series expansions. We have
                        \begin{equation}
                        \label{eq:studyMTildeMuKNeq0SecondPart42}
                            %
                        \end{equation}
                        We must now deal with the first term of the right-hand side of \eqref{eq:studyMTildeMuKNeq0SecondPart12}, given by
                        \begin{equation}
                        \label{eq:studyMTildeMuKNeq0SecondPart43}
                            %
                        \end{equation}
                        The first step towards dealing with this term is to perform an integration by parts which reverses a portion of the one effected in \eqref{eq:studyMTildeMuKNeq0SecondPart06}. We have
                        \begin{equation}
                        \label{eq:studyMTildeMuKNeq0SecondPart44}
                            %
                        \end{equation}
                        As we will see later, there is no need to study the first part of the right-hand side of \eqref{eq:studyMTildeMuKNeq0SecondPart44}. For any integer $k \geqslant 1$, we have
                        \begin{equation}
                        \label{eq:studyMTildeMuKNeq0SecondPart45}
                            %
                        \end{equation}
                        and we can compute the integrals above using \cite[Cor. C.32]{dutour:cuspDirichlet}. We thus have
                        \begin{equation}
                        \label{eq:studyMTildeMuKNeq0SecondPart46}
                            %
                        \end{equation}
                        We now use \eqref{eq:studyAmuKNeq0SixthPart42AnalyticExtHyperGeoT1MinusT} and the reflection formula on the Gamma function to get
                        \begin{equation}
                        \label{eq:studyMTildeMuKNeq0SecondPart47}
                            %
                        \end{equation}
                        for every real number $t \in \left] 0,1 \right[$. The term induced by the first part of the right-hand side of \eqref{eq:studyMTildeMuKNeq0SecondPart47} is given by
                        \begin{equation}
                        \label{eq:studyMTildeMuKNeq0SecondPart48}
                            %
                        \end{equation}
                        This term induces a holomorphic function near $0$, whose derivative at $s=0$ equals
                        \begin{equation}
                        \label{eq:studyMTildeMuKNeq0SecondPart49}
                            %
                        \end{equation}
                        The term induced by the second part of the right-hand side of \eqref{eq:studyMTildeMuKNeq0SecondPart47} is given by
                        \begin{equation}
                        \label{eq:studyMTildeMuKNeq0SecondPart50}
                            %
                        \end{equation}
                        Using the Euler integral formula and the binomial formula, we have
                        \begin{equation}
                        \label{eq:studyMTildeMuKNeq0SecondPart51}
                            %
                        \end{equation}
                        for any real number $t \in \left[ 0, 1 \right[$. The sum over $j \geqslant 1$ induced \eqref{eq:studyMTildeMuKNeq0SecondPart51} yields a holomorphic function around $0$, whose derivative at $s=0$ vanishes, because of the Pochhammer symbol. The term corresponding to $j = 0$ is given by
                        \begin{equation}
                        \label{eq:studyMTildeMuKNeq0SecondPart52}
                            %
                        \end{equation}
                        The first part of the right-hand side of \eqref{eq:studyMTildeMuKNeq0SecondPart52} needs not be studied, as it is canceled by the first part of \eqref{eq:studyMTildeMuKNeq0SecondPart44}. The second part of \eqref{eq:studyMTildeMuKNeq0SecondPart52}, given by
                        \begin{equation}
                        \label{eq:studyMTildeMuKNeq0SecondPart53}
                            %
                        \end{equation}
                        induces a holomorphic function near $0$, whose derivative at $s=0$ is given by
                        \begin{equation}
                        \label{eq:studyMTildeMuKNeq0SecondPart54}
                            %
                        \end{equation}
                        The last part of \eqref{eq:studyMTildeMuKNeq0SecondPart52}, given by
                        \begin{equation}
                        \label{eq:studyMTildeMuKNeq0SecondPart55}
                            %
                        \end{equation}
                        requires more care. On the half-plane $\Re s > 1$, we have
                        \begin{equation}
                        \label{eq:studyMTildeMuKNeq0SecondPart56}
                            %
                        \end{equation}
                        using the binomial formula. Thus, we can rewrite \eqref{eq:studyMTildeMuKNeq0SecondPart55} as
                        \begin{equation}
                        \label{eq:studyMTildeMuKNeq0SecondPart57}
                            %
                        \end{equation}
                        The last part of the right-hand side of \eqref{eq:studyMTildeMuKNeq0SecondPart57} induces a holomorphic function around $0$, whose derivative at $s=0$ vanishes, because of the Pochhammer symbol. The first part of the right-hand side of \eqref{eq:studyMTildeMuKNeq0SecondPart57} also induces a holomorphic function around $0$, and its derivative at $s=0$ is given by
                        \begin{equation}
                        \label{eq:studyMTildeMuKNeq0SecondPart58}
                            %
                        \end{equation}
                        Let us now deal with the second part of the right-hand side of \eqref{eq:studyMTildeMuKNeq0SecondPart57}. We have
                        \begin{equation}
                        \label{eq:studyMTildeMuKNeq0SecondPart59}
                            %
                        \end{equation}
                        The second part of \eqref{eq:studyMTildeMuKNeq0SecondPart59} induces a holomorphic function around $0$, whose derivative at $s=0$ vanishes. The first part of \eqref{eq:studyMTildeMuKNeq0SecondPart59} also induces a holomorphic function around $0$. Let us compute its derivative at $s=0$. We have
                        \begin{equation}
                        \label{eq:studyMTildeMuKNeq0SecondPart60}
                            %
                        \end{equation}
                        The proof of the proposition is thus complete.
                    \end{proof}
                    
                    Let us now take care of the case $\alpha = 0$.
                    
                    \begin{proposition}
                    \label{prop:studyMTildeMuKNeq0SecondPartAlpha0}
                        The function
                        \begin{equation}
                        \label{eq:functionStudyMTildeMuKNeq0SecondPartAlpha0}
                            %
                        \end{equation}
                        has a holomorphic continuation near $0$, whose derivative at $s=0$ equals
                        \begin{equation}
                            %
                        \end{equation}
                        This derivative, now taken with $\mu = 0$, further equals
                        \begin{equation}
                            %
                        \end{equation}
                        as $a$ goes to infinity.
                    \end{proposition}
                    
                    \begin{proof}
                        This is, up to a factor $-2$, done in \cite[Prop. 3.70]{dutour:cuspDirichlet}.
                    \end{proof}

                \subsubsection*{\textbf{Third part}}
                \label{MTildeMuKNeq0:part3}
                
                    We will now study of the term induced by the second part of the right-hand side of \eqref{eq:recallLogExpansionGkLargeOrder}.
                    
                    \begin{proposition}
                    \label{prop:studyMTildeMuKNeq0ThirdPart}
                        Assume we have $\alpha \neq 0$. The function
                        \begin{equation}
                        \label{eq:functionStudyMTildeMuKNeq0ThirdPart}
                            %
                        \end{equation}
                        has a holomorphic continuation near $0$. Its derivative at $s=0$ equals
                        \begin{equation}
                            %
                        \end{equation}
                        as $\mu$ goes to infinity. Furthermore, the same derivative, taken with $\mu = 0$, equals
                        \begin{equation}
                            %
                        \end{equation}
                    \end{proposition}
                    
                    \begin{proof}
                        The argument used here is similar to the one detailed in the proof of proposition \ref{prop:studyMTildeMuKNeq0SecondPart}. Note that we can use an analogue of \eqref{eq:studyAmuKNeq0FourthPart05TrickConjugateRoot}, having assumed $\alpha \neq 0$.
                    \end{proof}
                    
                    There is no counterpart to study when we have $\alpha = 0$, as the function \eqref{eq:functionStudyMTildeMuKNeq0SecondPartAlphaNeq0} vanishes identically in that case.

                \subsubsection*{\textbf{Fourth part}}
                \label{MTildeMuKNeq0:part4}
                
                    We now come to the last term induced by \eqref{eq:recallLogExpansionGkLargeOrder}.
                    
                    \begin{proposition}
                    \label{prop:studyMTildeMuKNeq0FourthPart}
                        Assume we have $\alpha \neq 0$. The function
                        \begin{equation}
                        \label{eq:functionStudyMTildeMuKNeq0FourthPart}
                            \begin{array}{lll}
                                \scriptstyle s & \longmapsto & \scriptstyle - \frac{1}{2} \frac{\sin \left( \pi s \right)}{\pi} \sum\limits_{k \neq 0} \int_{2 \left \vert k \right \vert^{\delta} \sqrt{\frac{1}{4} + \mu}}^{+ \infty} \left( t^2 - \left( \frac{1}{4} + \mu \right) \right)^{-s} \\[0.5em]
                                
                                && \qquad \qquad \qquad \qquad \scriptstyle \frac{\partial}{\partial t} \left( \frac{t^2}{\sqrt{4 \pi^2 \left( k + \alpha \right)^2 a^2 + t^2} + \sqrt{4 \pi^2 \left( k - \alpha \right)^2 a^2 + t^2}} \right. \\[0.5em]
                                
                                && \qquad \qquad \qquad \qquad \qquad \quad \scriptstyle \cdot  \left( \frac{1}{4 \pi^2 \left( k + \alpha \right)^2 a^2 + t^2} + \frac{1}{4 \pi^2 \left( k - \alpha \right)^2 a^2 + t^2} \right) \mathrm{d}t
                            \end{array}
                        \end{equation}
                        has a holomorphic continuation near $0$. Its derivative at $s=0$ is given by
                        \begin{equation}
                            \begin{array}{lll}
                                \scriptstyle \frac{\partial}{\partial s}_{\vert s = 0} \left[ \frac{\sin \left( \pi s \right)}{\pi} \left( 4 \mu + 1 \right)^{-s+1} \sum\limits_{k=1}^{+ \infty} \frac{1}{k^{2 \delta \left( s-1 \right)}} \right. \\[0.5em] 
                                
                                \scriptstyle \qquad \qquad \qquad \cdot \frac{1}{\sqrt{4 \pi^2 \left( k + \alpha \right)^2 a^2 + \left( 4 \mu + 1 \right) k^{2 \delta}} + \sqrt{4 \pi^2 \left( k - \alpha \right)^2 a^2 + \left( 4 \mu + 1 \right) k^{2 \delta}}} \\[0.5em] 
                                
                                \scriptstyle \qquad \qquad \qquad \qquad \left. \cdot \left( \frac{1}{4 \pi^2 \left( k + \alpha \right)^2 a^2 + \left( 4 \mu + 1 \right) k^{2 \delta}} + \frac{1}{4 \pi^2 \left( k - \alpha \right)^2 a^2 + \left( 4 \mu + 1 \right) k^{2 \delta}} \right) \right] - \frac{1}{2 \pi a}
                            \end{array}
                        \end{equation}
                        and the same derivative, taken with $\mu = 0$, vanishes as $a$ goes to infinity.
                    \end{proposition}
                    
                    \begin{proof}
                        This follows the same type of argument as the proof of proposition \ref{prop:studyMTildeMuKNeq0SecondPart}.
                    \end{proof}
                    
                    Let us now deal with the case $\alpha = 0$.
                    
                    \begin{proposition}
                    \label{prop:studyMTildeMuKNeq0FourthPartAlpha0}
                        The function
                        \begin{equation}
                            \begin{array}{lll}
                                \scriptstyle s & \longmapsto & \scriptstyle - \frac{1}{2} \frac{\sin \left( \pi s \right)}{\pi} \sum\limits_{k \neq 0} \int_{2 \left \vert k \right \vert^{\delta} \sqrt{\frac{1}{4} + \mu}}^{+ \infty} \left( t^2 - \left( \frac{1}{4} + \mu \right) \right)^{-s} \frac{\partial}{\partial t} \left( \frac{t^2}{\left( 4 \pi^2 k^2 a^2 + t^2 \right)^{3/2}} \right) \mathrm{d}t
                            \end{array}
                        \end{equation}
                        has a holomorphic continuation near $0$. Its derivative at $s=0$ is given by
                        \begin{equation}
                            \begin{array}{lll}
                                \scriptstyle \frac{\partial}{\partial s}_{\vert s = 0} \left[ \frac{\sin \left( \pi s \right)}{\pi} \left( 4 \mu + 1 \right)^{-s+1} \sum\limits_{k=1}^{+ \infty} \frac{1}{k^{2 \delta \left( s-1 \right)}} \cdot \frac{1}{\left( 4 \pi^2 \left( k + \alpha \right)^2 a^2 + \left( 4 \mu + 1 \right) k^{2 \delta} \right)^{3/2}} \right] - \frac{1}{2 \pi a}
                            \end{array}
                        \end{equation}
                        and the same derivative, taken with $\mu = 0$, vanishes as $a$ goes to infinity.
                    \end{proposition}
                    
                    \begin{proof}
                        Up to a factor, this is included in \cite[Prop. 3.72]{dutour:cuspDirichlet}, and follows a reasonning similar to proposition \ref{prop:studyMTildeMuKNeq0SecondPart}.
                    \end{proof}

\section{The determinant of the pseudo-Laplacian with Alvarez--Wentworth boundary conditions}

    In this final section, we will put together all the results obtained in this paper, to obtain the asymptotics studies of the zeta-determinant of the pseudo-Laplacian with Alvarez--Wentworth boundary conditions.
    
    \subsection{Results on the zeta functions}
    \label{subsec:resultsOnZetaFunctions}
    
        We begin by stating the relevant results on the partial spectral zeta functions and their derivatives at $s=0$.
    
        \subsubsection{The function $\zeta_{L,\mu}^{0}$}
        
            Let us first summarize the results obtained so far regarding the partial spectral zeta function which corresponds to the integer $k = 0$ in the Fourier decomposition \eqref{eq:fourierDecompositions}. Let us assume that we have $\alpha \neq 0$.
            
            \begin{theorem}
            \label{thm:asymptZetaFunction0}
                The function $\zeta_{L,\mu}^{0}$ has a holomorphic continuation near $0$, and its derivative at $s=0$ satisfies
                \begin{equation}
                    \begin{array}{lll}
                        \zeta_{L,\mu}^{0} \, ' \left( 0 \right) & = & \left \{ \begin{array}{ll} - \sqrt{\mu} \log \mu + 2 \left( \log \left( \pi \alpha a \right) - 1 \right) \sqrt{\mu} \\[0.5em]
                        
                            \qquad - \log \mu  + 2 \log 2 + o \left( 1 \right) & \text{as } \mu \rightarrow + \infty. \\[1em]
                            
                            2 \pi \alpha a - \frac{1}{2} \log a + \frac{3}{2} \log 2 & \text{as } a \rightarrow + \infty,\\[0.5em]
                            
                            \qquad - \frac{1}{2} \log \left( \pi \alpha \right) + o \left( 1 \right) & \text{ with } \mu = 0.
                        \end{array} \right.
                    \end{array}
                \end{equation}
            \end{theorem}
            
            \begin{proof}
                This is a consequence of theorem \ref{thm:asymptDerAt0DirichletContribution0} and propositions \ref{prop:studyBMuK0}, \ref{prop:studyAMu0AsMuGoesToInfinity}, \ref{prop:studyAMu0AsAGoesToInfinity}, \ref{prop:studyRMuK0}, and \ref{prop:studyMTildeMuK0}.
            \end{proof}
        
        \subsubsection{The zeta function $\zeta_{L,\mu}^{\pm}$}
            
            More care is required in order to collect the various results pertaining to the partial spectral zeta function which corresponds to the integers $k \neq 0$ in the Fourier decomposition \eqref{eq:fourierDecompositions}. We do not make any general assumption on $\alpha$.
            
            \begin{theorem}
            \label{thm:asymptZetaFunctionPMAlphaNeq0AsMuGoesToInfinity}
                The function $\zeta_{L,\mu}^{\pm}$ has a holomorphic continuation near $0$. Its derivative at $s=0$ satisfies
                \begin{equation}
                    \begin{array}{lll}
                        \zeta_{L,\mu}^{\pm} \, ' \left( 0 \right) & = & \frac{1}{2 \pi a} \mu \log \mu - \frac{1}{2 \pi a} \mu - 3 \sqrt{\mu} \log \mu + 2 \left[ - 3 \log 2 + \alpha \log \frac{1 + \alpha}{1 - \alpha} \right. \\[1em]
                        
                        && \quad + \frac{1}{2 a} - 2 \int_0^{+ \infty} \frac{1}{e^{2 \pi t} - 1} \left( \arctan \frac{t}{1 + \alpha} + \arctan \frac{t}{1 - \alpha} \right) \mathrm{d}t \\[1em]
                        
                        && \quad \quad \left. + 1 + \frac{3}{2} \log \left( 4 \pi^2 \left( 1 - \alpha^2 \right) a^2 \right) \right] \sqrt{\mu} - \alpha \log \mu + o \left( 1 \right)
                    \end{array}
                \end{equation}
                as $\mu$ goes to infinity.
            \end{theorem}
            
            \begin{proof}
                This is a consequence of theorem \ref{thm:asymptDerAt0DirichletContributionPM}, as well as of propositions \ref{prop:studyBMuKNeq0}, \ref{prop:studyAMuKNeq0FirstPart}, \ref{prop:studyAMuKNeq0SecondPart}, \ref{prop:studyAMuKNeq0ThirdPart}, \ref{prop:studyAMuKNeq0FourthPart}, \ref{prop:studyAMuKNeq0FifthPart}, \ref{prop:studyAMuKNeq0SixthPart}, \ref{prop:studyAMuKNeq0SeventhPart}, \ref{prop:studyAMuKNeq0EighthPart}, \ref{prop:studyAMuKNeq0NinthPart}, \ref{prop:studyRMuKNeq0}, \ref{prop:studyMTildeMuKNeq0FirstPart}, \ref{prop:studyMTildeMuKNeq0SecondPart}, \ref{prop:studyMTildeMuKNeq0SecondPartAlpha0}, \ref{prop:studyMTildeMuKNeq0ThirdPart}, \ref{prop:studyMTildeMuKNeq0FourthPart}, and \ref{prop:studyMTildeMuKNeq0FourthPartAlpha0}.
            \end{proof}
            
            We can now state the result for the $a$-asymptotic expansion.
            
            \begin{theorem}
            \label{thm:asymptZetaFunctionPMAsAGoesToInfinity}
                The function $\zeta_{L,\mu}^{\pm}$ has a holomorphic continuation near $0$. Its derivative at $s=0$ satisfies
                \begin{equation}
                    \begin{array}{lll}
                        \zeta_{L,0}^{\pm} \, ' \left( 0 \right) & = & - \frac{2 \pi}{3} a - 4 \pi \alpha^2 a - 2 \alpha \log a + \log \frac{\sin \left( \pi \alpha \right)}{\pi \alpha} + 2 \log \Gamma \left( 1 - \alpha \right) \\[1em]
                        
                        && \qquad + 2 \gamma \alpha - 2 \alpha \log \left( 4 \pi \right) + o \left( 1 \right)
                    \end{array}
                \end{equation}
                as $a$ goes to infinity.
            \end{theorem}
            
            \begin{proof}
                This is a consequence of theorem \ref{thm:asymptDerAt0DirichletContributionPM}, as well as of propositions \ref{prop:studyBMuKNeq0}, \ref{prop:globalStudyAMuKNeq0}, \ref{prop:studyRMuKNeq0}, \ref{prop:studyMTildeMuKNeq0FirstPart}, \ref{prop:studyMTildeMuKNeq0SecondPart}, \ref{prop:studyMTildeMuKNeq0SecondPartAlpha0}, \ref{prop:studyMTildeMuKNeq0ThirdPart}, \ref{prop:studyMTildeMuKNeq0FourthPart}, \ref{prop:studyMTildeMuKNeq0FourthPartAlpha0}.
            \end{proof}
    
    \subsection{Results on the determinants}
    
        We will now state the theorems from subsection \ref{subsec:resultsOnZetaFunctions} as results on the (logarithms of the) determinant of the pseudo-Laplacian with Alvarez--Wentworth boundary conditions.
        
        \subsubsection{\texorpdfstring{As $\mu \rightarrow + \infty$}{As mu goes to infinity}}
        
            We begin with the $\mu$-asymptotic expansion.
            
            \begin{theorem}
            \label{thm:asymptMuDeterminantAlphaNeq0}
                Let us assume we have $\alpha \neq 0$ and $a > 1/\left( 4 \pi \left( 1 - \alpha \right) \right)$. We have
                \begin{equation}
                    \begin{array}{lll}
                        \multicolumn{3}{l}{\log \det \left( \Delta_{L,0} + \mu \right)} \\[0.5em]
                        
                        \quad & = &  \frac{1}{2 \pi a} \mu \log \mu + \frac{1}{2 \pi a} \mu + 4 \sqrt{\mu} \log \mu - 2 \left[ \frac{1}{2 a} + \alpha \log \frac{1 + \alpha}{1 - \alpha} + \log \alpha \right. \\[1em]
                        
                        && \quad + \frac{3}{2} \log \left( 1 - \alpha^2 \right) - 2 \int_0^{+ \infty} \frac{1}{e^{2 \pi t} - 1} \left( \arctan \frac{t}{1 + \alpha} + \arctan \frac{t}{1 - \alpha} \right) \mathrm{d}t \\[1em]
                        
                        && \quad \quad \left. + 4 \log \left( \pi a \right) \right] \sqrt{\mu} + \left( 1 + \alpha \right) \log \mu - 2 \log 2 + o \left( 1 \right)
                    \end{array}
                \end{equation}
                as $\mu$ goes to infinity.
            \end{theorem}
            
            \begin{proof}
                This is a direct consequence of theorems \ref{thm:asymptZetaFunction0} and \ref{thm:asymptZetaFunctionPMAlphaNeq0AsMuGoesToInfinity}.
            \end{proof}
            
            \begin{theorem}
            \label{thm:asymptMuDeterminantAlpha0}
                Let us assume we have $\alpha = 0$ and $a > 1/\left( 4 \pi \right)$. We have
                \begin{equation}
                    \begin{array}{lll}
                        \multicolumn{3}{l}{\log \det \left( \Delta_{L,0} + \mu \right)} \\[0.5em] 
                        
                        \qquad & = & - \frac{1}{2 \pi a} \mu \log \mu + \frac{1}{2 \pi a} \mu + 3 \sqrt{\mu} \log \mu - 2 \left[ \frac{1}{2a} + 1 + 3 \log \left( \pi a \right) \right. \\[1em]
                        
                        && \qquad \left. - 4 \int_0^{+ \infty} \frac{1}{e^{2 \pi t} - 1} \arctan \left( t \right) \mathrm{d}t \right] \sqrt{\mu} + o \left( 1 \right)
                    \end{array}
                \end{equation}
                as $\mu$ goes to infinity.
            \end{theorem}
            
            \begin{proof}
                This is a direct corollary to theorem \ref{thm:asymptZetaFunctionPMAlphaNeq0AsMuGoesToInfinity}.
            \end{proof}
            
            \begin{remark}
                Note that the assumption $a > 1/\left( 4 \pi \left( 1 - \alpha \right) \right)$ was necessary in proposition \ref{prop:zerosfk} to localize the eigenvalues of the pseudo-Laplacian. It must therefore be kept in the final result.
            \end{remark}

        \subsubsection{\texorpdfstring{As $a \rightarrow + \infty$}{As a goes to infinity}}
    
            We end with the $a$-asymptotic expansion.
            
            \begin{theorem}
            \label{thm:asymptADeterminantAlphaNeq0}
                Let us assume we have $\alpha \neq 0$. We have
                \begin{equation}
                    \begin{array}{lll}
                        \log \det' \Delta_{L,0} & = & \frac{2 \pi}{3} a + 4 \pi \alpha^2 a - 2 \pi \alpha a + \left( \frac{1}{2} + 2 \alpha \right) \log a - \log \frac{\sin \left( \pi \alpha \right)}{\pi \alpha} \\[1em]
                        
                        && \qquad - 2 \log \Gamma \left( 1 - \alpha \right)- 2 \alpha \gamma + 2 \alpha \log \left( 4 \pi \right) - \frac{3}{2} \log 2 \\[1em]
                        
                        && \qquad \qquad + \frac{1}{2} \log \left( \pi \alpha \right) + o \left( 1 \right)
                    \end{array}
                \end{equation}
                as $a$ goes to infinity.
            \end{theorem}
            
            \begin{proof}
                This is a direct consequence of theorems \ref{thm:asymptZetaFunction0} and \ref{thm:asymptZetaFunctionPMAsAGoesToInfinity}.
            \end{proof}
            
            \begin{theorem}
            \label{thm:asymptADeterminantAlpha0}
                Let us assume we have $\alpha = 0$. We have
                \begin{equation}
                    \begin{array}{lll}
                        \log \det' \Delta_{L,0} & = & \frac{2 \pi}{3} a + o \left( 1 \right)
                    \end{array}
                \end{equation}
                as $a$ goes to infinity.
            \end{theorem}
            
            \begin{proof}
                This is a direct corollary to theorem \ref{thm:asymptZetaFunctionPMAsAGoesToInfinity}.
            \end{proof}

\appendix

\bibliographystyle{amsplain}

\bibliography{bibliography}

\end{document}